\documentclass[aop,seceqn,citesort,MSNbibl,dvips]{arximspdf}

\doi{10.1214/10-AOP571}
\volume{39}
\issue{3}
\pubyear{2011}
\firstpage{779}
\lastpage{856}

\makeatletter

\newtheorem{theo}{Theorem}[section]
\newtheorem{lemma}{Lemma}[section]
\newtheorem{prop}{Proposition}[section]
\newtheorem{corollary}{Corollary}[section]

\newproclaim{rem}{Remark}[section]
\newproclaim{defi}{Definition}[section]

\newcommand{\ho}{\mathrm{hom}}
\newcommand{\R}{\mathbb{R}}
\newcommand{\Z}{\mathbb{Z}}
\newcommand{\N}{\mathbb{N}}
\newcommand{\Id}{\operatorname{Id}}
\newcommand{\e}{\varepsilon}
\newcommand{\calA}{\mathcal{A}_{\alpha\beta}}
\newcommand{\ee}{\mathbf{e}}
\newcommand{\ener}{\mathcal{E}}

\newcommand{\llangle}{\langle\!\langle}
\newcommand{\rrangle}{\rangle\!\rangle}

\makeatother

\begin{document}
\begin{frontmatter}

\title{An optimal variance estimate in stochastic homogenization
of discrete elliptic equations}
\runtitle{Variance estimate for effective diffusion}

\begin{aug}
\author[A]{\fnms{Antoine} \snm{Gloria}\corref{}\ead
[label=e1]{antoine.gloria@inria.fr}} and
\author[B]{\fnms{Felix} \snm{Otto}\ead[label=e2]{otto@mis.mpg.de}}
\runauthor{A. Gloria and F. Otto}
\affiliation{INRIA and Max Planck Institute}
\address[A]{Project-team SIMPAF\\
INRIA Lille-Nord Europe\\
France\\
\printead{e1}}
\address[B]{Max Planck Institute\\
for Mathematics in the Sciences\\
Leipzig\\
Germany\\
\printead{e2}}
\end{aug}

\received{\smonth{5} \syear{2009}}
\revised{\smonth{5} \syear{2010}}

%
\begin{abstract}
We consider a discrete elliptic equation on the $d$-dimensional lattice
$\mathbb{Z}^d$ with random coefficients $A$ of the simplest type: they
are identically distributed and independent from edge to edge. On
scales large w.r.t. the lattice spacing (i.e., unity), the solution
operator is known to behave like the solution operator of a
(continuous) elliptic equation with constant deterministic
coefficients. This symmetric ``homogenized'' matrix $A_\ho=a_{\ho}\Id$
is characterized by $\xi\cdot A_{\ho}\xi=\langle(\xi+\nabla\phi
)\cdot
A(\xi+\nabla\phi)\rangle$ for any direction $\xi\in\mathbb
{R}^d$, where
the random field $\phi$ (the ``corrector'') is the unique solution of
$-\nabla^*\cdot A(\xi+\nabla\phi) = 0$ such that $\phi(0)=0$,
$\nabla\phi$ is stationary and $\langle\nabla\phi\rangle=0$,
$\langle\cdot\rangle$~denoting the ensemble average (or expectation).

It is known (``by ergodicity'') that the above ensemble average of the
energy density $\ener=(\xi+\nabla\phi)\cdot A(\xi+\nabla\phi)$, which
is a stationary random field, can be recovered by a system average. We
quantify this by proving that the variance of a spatial average of
$\ener$ on length scales $L$ satisfies the optimal estimate, that is,
$\operatorname{var}[\sum\ener\eta_L] \lesssim L^{-d}$, where the
averaging function
[i.e., $\sum\eta_L=1$, $\operatorname{supp}(\eta_L)\subset\{|x|\le
L\}$] has to be
smooth in the sense that $|\nabla\eta_L|\lesssim L^{-1-d}$. In two
space dimensions (i.e., $d=2$), there is a logarithmic correction. This
estimate is optimal since it shows that smooth averages of the energy
density $\ener$ decay in $L$ as if $\ener$ would be independent from
edge to edge (which it is not for $d>1$).

This result is of practical significance, since it allows to
estimate the dominant error when numerically computing $a_{\ho}$.
\end{abstract}

%

\end{frontmatter}

\section{Introduction}

\subsection{Motivation, informal statement and optimality of the result}

We study discrete elliptic equations. More precisely,
we consider real functions $u$ of the sites $x$
in a $d$-dimensional Cartesian lattice $\mathbb{Z}^d$.
Every edge $e$ of the lattice is endowed with a ``conductivity''
$a(e)>0$. This defines a discrete elliptic differential operator
$-\nabla^*\cdot A\nabla$ via
\[
-\nabla^*\cdot(A\nabla u)(x) := \sum_{y\in\Z^d, |x-y|=1} a(e)\bigl(u(x)-u(y)\bigr),
\]
where the sum is over the $2d$ sites $y$ which are connected by
an edge $e=[x,y]$ to the site $x$. It is sometimes more convenient
to think in terms of the associated Dirichlet form, that is,
\begin{eqnarray*}
\sum_{x\in\Z^d}(\nabla v\cdot A\nabla u)(x)
:\!&=&
\sum_{x\in\mathbb{Z}^d}v(x) \bigl(-\nabla^*\cdot(A\nabla u)(x)\bigr)\\
&=&
\sum_{e}\bigl(v(x)-v(y)\bigr)a(e)\bigl(u(x)-u(y)\bigr),
\end{eqnarray*}
where the last sum is over all edges $e$, and $(x,y)$ denotes the two sites
connected by $e$, that is, $e=[x,y]=[y,x]$ (with the convention that an
edge is not oriented).
We assume the conductivities
$a$ to be uniformly elliptic in the sense of
\[
\alpha\le a(e) \le\beta\qquad\mbox{for all edges } e
\]
for some fixed constants $0<\alpha\le\beta<\infty$.

We are interested in random coefficients. To fix ideas,
we consider the simplest situation possible:
\[
\{a(e)\}_{e}\qquad \mbox{are independently and identically distributed
(i.i.d.)}.
\]
Hence, the statistics are described by a distribution on the finite
interval $[\alpha,\beta]$. We would like to see this discrete elliptic operator
with random coefficients as a good model problem for
continuum elliptic operators with random coefficients of correlation
length unity.

The first results in stochastic homogenization of linear
elliptic equations in the
continuous setting are due to Kozlov \cite{Kozlov-79} and
Papanicolaou and Varadhan \cite{Papanicolaou-Varadhan-79},
essentially using compensated compactness.
The adaptation of these results to discrete elliptic equations
in quite more general situations
than the one considered above (i.e., under general ergodic assumptions)
is due to K\"unnemann \cite{Kunnemann-83} following the approach by
Papanicolaou and Varadhan for the continuous case, and also to Kozlov
\cite{Kozlov-87} (where more
general discrete elliptic operators are considered).
Note that the discrete elliptic operator $-\nabla^*\cdot A\nabla$ is
the infinitesimal generator of a random walk in
a random environment, whence the rephrasing of the homogenization
result in \cite{Kunnemann-83}
as the diffusion limit for reversible jump processes in $\Z^d$ with
random bond conductivities.
With the same point of view, it is also worth mentioning the seminal paper
by Kipnis and Varadhan \cite{Kipnis-Varadhan-86} using central limit
theorems for martingales.

The general homogenization result proved in these articles states
that there exist \textit{homogeneous and deterministic} coefficients
$A_{\ho}$
such that the solution operator of the continuum differential operator
$-\nabla\cdot A_{\ho}\nabla$
describes the large scale behavior of the solution operator
of the discrete differential operator $-\nabla^*\cdot A\nabla$.
As a by product of this homogenization result,
one obtains a characterization of the homogenized
coefficients $A_{\ho}$: it is shown that for
every direction $\xi\in\mathbb{R}^d$, there exists a unique
scalar field $\phi$ such that $\nabla\phi$ is stationary
[stationarity means that the fields
$\nabla\phi(\cdot)$ and $\nabla\phi(\cdot+z)$ have the same statistics
for all shifts $z\in\mathbb{Z}^d$] and $\langle\nabla\phi\rangle
=0$, solving
the equation
%
%
\begin{equation}\label{PV2}
-\nabla^*\cdot\bigl(A(\xi+\nabla\phi)\bigr) = 0 \qquad\mbox{in } \mathbb{Z}^d,
\end{equation}
and normalized by $\phi(0)=0$. As in periodic homogenization,
the function $\mathbb{Z}^d\ni x\mapsto\xi\cdot x+\phi(x)$
can be seen as the $A$-harmonic function which macroscopically
behaves as the affine function $\mathbb{Z}^d\ni x\mapsto\xi\cdot x$.
With this
``corrector'' $\phi$, the homogenized coefficients $A_{\ho}$
(which in general form a symmetric matrix and for our simple statistics
in fact a multiple of the identity: $A_{\ho}=a_{\ho}\Id$)
can be characterized as follows:
%
%
\begin{equation}\label{PV1}
\xi\cdot A_{\ho}\xi=
\langle(\xi+\nabla\phi)\cdot A(\xi+\nabla\phi)\rangle.
\end{equation}
Since
the scalar field
$(\xi+\nabla\phi)\cdot A(\xi+\nabla\phi)$ is stationary, it does
not matter (in terms of the distribution) at which site $x$ it
is evaluated in the formula (\ref{PV1}), so that we suppress
the argument $x$ in our notation.

The representation (\ref{PV1}) is of no immediate practical use,
since the equation (\ref{PV2}) has to be solved:
\begin{itemize}
\item for \textit{every realization} of the coefficients
$\{a(e)\}_{e}$ and
\item in the \textit{whole space} $\mathbb{Z}^d$.
\end{itemize}
In order to overcome the first difficulty, it is natural
to appeal to ergodicity (in the sense that ensemble averages
are equal to system averages), which suggests to replace
(\ref{PV1}) by
%
%
\begin{equation}\label{PV5}
\xi\cdot A_{\ho}\xi\leadsto
\sum(\xi+\nabla\phi)\cdot A(\xi+\nabla\phi)\eta_L,
\end{equation}
where $\eta_L$ is a suitable averaging function of length scale $L\gg1$,
that is,
%
%
\begin{equation}\label{PV13}
\operatorname{supp}(\eta_L)\subset\{|x|\le L\},\qquad
|\eta_L|\lesssim L^{-d},\qquad
\sum\eta_L=1.
\end{equation}
In fact, on expects the
energy density $(\xi+\nabla\phi)\cdot A(\xi+\nabla\phi)$, which is
a stationary random field, to display a decay of correlations
over large distances, so that (\ref{PV5}) seems a good approximation
for $L\gg1$.

However, one still has to solve (\ref{PV2}) on the whole space
$\mathbb{Z}^d$, albeit for a single realization of the coefficients.
In order to overcome this second difficulty, we start with
the following observation:
since $\phi$ on the ball $\{|x|\le L\}$ is expected to be little correlated
to $\phi$ outside the ball $\{|x|\ge R\}$ for
$R-L\gg1$, it seems natural to replace $\phi$ in (\ref{PV5})
by $\phi_R$:
%
%
\begin{equation}\label{PV9bis}
\sum(\xi+\nabla\phi)\cdot A(\xi+\nabla\phi)\eta_L
\leadsto
\sum(\xi+\nabla\phi_R)\cdot A(\xi+\nabla\phi_R)\eta_L,
\end{equation}
where $\phi_R$ is the solution
of an equation on a domain (say, a ball) of size $R$
with homogeneous boundary conditions (say, Dirichlet):
%
%
\begin{eqnarray}\label{PV3}
-\nabla^*\cdot\bigl(A(\xi+\nabla\phi_R)\bigr)&=&0
\qquad\mbox{in } \mathbb{Z}^d\cap\{|x|<R\},\nonumber\\[-8pt]\\[-8pt]
\phi_R&=&0\qquad\mbox{in } \mathbb{Z}^d\cap\{|x|\ge R\},\nonumber
\end{eqnarray}
so that the right-hand side of (\ref{PV9bis}) is indeed computable.

However, $\nabla\phi_R$ defined by (\ref{PV3}) is not
statistically stationary, which is a handicap for the
error analysis. It is therefore common in the analysis of the
error from spatial cut-off\vadjust{\goodbreak}
to introduce an intermediate step which consists in
replacing equation (\ref{PV2}) by
%
%
\begin{equation}\label{PV4}
T^{-1}\phi_T-\nabla^*\cdot\bigl(A(\xi+\nabla\phi_T)\bigr) = 0
\qquad\mbox{in } \mathbb{Z}^d.
\end{equation}
Clearly, the zero order term in (\ref{PV4}) introduces a characteristic
length scale $\sqrt{T}$ (the notation $T$ that alludes to time is used because
$T^{-1}$ corresponds to the death rate in the random walker
interpretation of the operator $T^{-1}-\nabla^*\cdot A\nabla$).
In a second step, (\ref{PV4}) is then replaced by
\begin{eqnarray*}
T^{-1}\phi_T-\nabla^*\cdot\bigl(A(\xi+\nabla\phi_{T,R})\bigr)&=&0
\qquad\mbox{in } \mathbb{Z}^d\cap\{|x|<R\},\\
\phi_{T,R}&=&0 \qquad\mbox{in } \mathbb{Z}^d\cap\{|x|\ge R\}.
\end{eqnarray*}
The Green's function $G_T(x,y)$ of
the operator $T^{-1}-\nabla^*\cdot A\nabla$ is known to decay faster
than any power in $\frac{\sqrt{T}}{|x-y|}\ll1$ uniformly
in the realization of the coefficients [see, in particular,
Lemma \ref{L11a}(iii)]. Therefore,
one expects that $\phi_T$ and $\phi_{T,R}$ agree on
the ball $\{|x|\le L\}$ up to an error which is of \textit{infinite} order
in $\e=\frac{\sqrt{T}}{R-L}$ ($\e$ is the inverse of the distance of
the ball $\{|x|\le L\}$ to the Dirichlet boundary $\{|x|=R\}$
measured in units of $\sqrt{T}$, see, e.g.,
\cite{Bourgeat-04}, Section 3, for related arguments).
Hence, we shall consider
$\sum(\xi+\nabla\phi_T)\cdot A(\xi+\nabla\phi_T)\eta_L$
as a very good proxy to the practically computable
$\sum(\xi+\nabla\phi_{T,R})\cdot A(\xi+\nabla\phi_{T,R})\eta_L$:
\[
\sum(\xi+\nabla\phi_T)\cdot A(\xi+\nabla\phi_T)\eta_L
\approx
\sum(\xi+\nabla\phi_{T,R})\cdot A(\xi+\nabla\phi_{T,R})\eta_L.
\]

In view of this remark, we restrict our attention to
the error we make when replacing
\[
\xi\cdot A_{\ho}\xi\leadsto
\sum(\xi+\nabla\phi_T)\cdot A(\xi+\nabla\phi_T)\eta_L.
\]
It is natural to measure this error in terms of the expected value of its
square. This error splits into two parts, the first arising from
the finiteness of the averaging length\vspace*{1pt} scale $L$
and the other arising from the finiteness of the cut-off length
scale $\sqrt{T}$:
%
%
\begin{eqnarray}\qquad
&&\Bigl\langle\Bigl|
\sum(\xi+\nabla\phi_T)\cdot A(\xi+\nabla\phi_T)\eta_L-\xi\cdot
A_{\ho
}\xi
\Bigr|^2\Bigr\rangle\nonumber\\
&&\qquad\stackrel{\mbox{\fontsize{8.36pt}{10.36pt}\selectfont{(\ref{PV1})}}}{=}
\Bigl\langle\Bigl|\sum(\xi+\nabla\phi_T)\cdot A(\xi+\nabla\phi_T)\eta_L
-\langle(\xi+\nabla\phi)\cdot A(\xi+\nabla\phi)\rangle
\Bigr|^2\Bigr\rangle\nonumber\\
\label{PV6}
&&\qquad\hspace*{3.84pt}=
\operatorname{var}\Bigl[\sum(\xi+\nabla\phi_T)\cdot A(\xi+\nabla\phi
_T)\eta_L\Bigr]
\\
\label{PV8bis}
&&\qquad\quad\hspace*{3.84pt}{} +
\Bigl|\Bigl\langle\sum(\xi+\nabla\phi_T)\cdot A(\xi+\nabla\phi_T)
\eta_L\Bigr\rangle
-\langle(\xi+\nabla\phi)\cdot A(\xi+\nabla\phi)\rangle
\Bigr|^2.
\end{eqnarray}
In view of the stationarity of $(\xi+\nabla\phi_T)\cdot A(\xi
+\nabla\phi_T)$, of (\ref{PV13}) and of (\ref{PV2}), the second part
(\ref{PV8bis}) of the error can be rewritten as
%
%
\begin{eqnarray}\label{PV8}
&&\Bigl|\Bigl\langle\sum(\xi+\nabla\phi_T)\cdot A(\xi
+\nabla\phi_T)
\eta_L\Bigr\rangle
-\langle(\xi+\nabla\phi)\cdot A(\xi+\nabla\phi)\rangle
\Bigr|^2\nonumber\\
&&\qquad=|\langle(\xi+\nabla\phi_T)\cdot A(\xi+\nabla\phi_T) -(\xi
+\nabla\phi)\cdot A(\xi+\nabla\phi) \rangle|^2\\
&&\qquad=
\langle(\nabla\phi_T-\nabla\phi)\cdot A(\nabla\phi_T-\nabla\phi
)\rangle
^2.\nonumber
\end{eqnarray}

What scaling can we expect for the two error terms
(\ref{PV6}) and (\ref{PV8})? A heuristic prediction
can be easily inferred from the regime of small ellipticity contrast,
that is, $1-\frac{\alpha}{\beta}\ll1$ (and $\alpha=1$ w.l.o.g.).
In this regime, to leading order, the two error terms
(\ref{PV6}) and (\ref{PV8}) behave like
\[
\operatorname{var}\Bigl[\sum\bigl( \xi\cdot(A-\langle A\rangle)\xi+2\xi
\cdot\nabla\bar\phi\bigr)\eta_L\Bigr]
\quad\mbox{and}\quad
\langle|\nabla\bar\phi_T-\nabla\bar\phi|^2\rangle^2,
\]
where $\bar\phi$ and $\bar\phi_T$ are defined via
%
%
\begin{eqnarray}
\label{PV9}
-\triangle\bar\phi
&=&
\nabla^*\cdot\bigl((A-\langle A\rangle)\xi\bigr),
\\
\label{PV10}
T^{-1}\bar\phi_T-\triangle\bar\phi_T
&=&
\nabla^*\cdot\bigl((A-\langle A\rangle)\xi\bigr),
\end{eqnarray}
respectively. In the first error term, we have replaced $\bar\phi_T$ by
$\bar\phi$ for simplicity of the exposition.

These error terms can be computed in a straightforward
manner. Indeed,
as shown in the \hyperref[app]{Appendix}, they scale for any direction
$|\xi|=1$ as:
%
%
\begin{eqnarray}\qquad
\label{pv10bis}
\operatorname{var}\Bigl[\sum\bigl(\xi\cdot(A-\langle A \rangle)\xi+2\xi
\cdot\nabla\bar\phi\bigr)\eta_L\Bigr]
&\sim&L^{-d}\\
\label{PV11}
\langle|\nabla\bar\phi_T-\nabla\bar\phi|^2\rangle^2
&\sim&
\cases{
T^{-d}, &\quad for $d<4$,\cr
T^{-4}\ln^2 T, &\quad for $d=4$,\cr
T^{-4}, &\quad for $d>4$.}
\end{eqnarray}

We now argue that the first error term (\ref{pv10bis})
is the dominant one (in dimensions $d<8$).
In order to do so, we argue that the choice of $L\sim\sqrt{T}$
is natural [for which (\ref{pv10bis}) dominates (\ref{PV11})
in dimensions $d<8$]. Indeed, we recall that in the ball $\{|x|\le L\}$,
$\phi_T$ is a proxy for the computable $\phi_{T,R}$
(defined on the larger ball $\{|x|\le R\}$). The error
is of \textit{infinite} order in the distance between the two balls,
measured in the length scale $\sqrt{T}$, that is,
in $\e:=\sqrt{T}/(R-L)\ll1$.
Hence, for the sake of discussing rates, we may indeed think
of $L\sim\sqrt{T}\sim R$.

In this paper, we therefore focus on the error
term (\ref{PV6}) coming from the finite range $L$ of
the spatial average. In Theorem \ref{th:main} (see also Remark \ref{Remark1}),
we shall establish
that (\ref{pv10bis}) holds as an estimate also for its
nonlinear counterpart (\ref{PV6}), that is,
%
%
\begin{equation}\label{PV12}
\operatorname{var}\Bigl[\sum(\xi+\nabla\phi_T)\cdot A(\xi+\nabla\phi
_T)\eta_L\Bigr]
\lesssim
L^{-d},
\end{equation}
with two minor restrictions:
\begin{itemize}
\item In dimension $d=2$, the prefactor depends logarithmically
on $T$ (whereas for $d\not=2$, the prefactor depends only on the
ellipticity constants).
\item The spatial\vspace*{1pt} averaging function $\eta_L$ has to be smooth
in the sense that $|\nabla\eta_L|\lesssim L^{-d-1}$ in addition
to (\ref{PV13}).
\end{itemize}
The estimate for the higher order term (\ref{PV8bis}) will be the
object of
a subsequent work.


\subsection{Discussion of the works of Yurinskii and
of Naddaf and Spencer}

In this subsection, we comment on two papers
on error estimates (in the sense of the previous subsection)
which from our perspective are the essential ones. We also
explain how our work relates to these two papers.

Still unsurpassed is the first quantitative paper,
the inspiring 1986 work by Yurinskii \cite{Yurinskii-86}.
He essentially deals with the error (\ref{PV8bis}) arising from
the spatial cut-off $T$. In our discrete setting of
i.i.d. coefficients $a(e)$ and for dimension $d>2$,
his result translates into
%
%
\begin{equation}\label{Yu1}
\langle|\nabla\phi_T-\nabla\phi|^2\rangle
\lesssim
T^{({2-d})/({4+d})+\delta},
\end{equation}
for $T\gg1$ and some arbitrarily small $\delta>0$, see
\cite{Yurinskii-86}, Theorem 2.1 (and \cite{E-05}, Lemma A.5,
for this rephrasing of Yurinskii's result).

Yurinskii derives estimate (\ref{Yu1}) by fairly elementary
arguments from the following crucial \textit{variance estimate} of the
spatial averages $\sum\phi_T\eta_L$ of $\phi_T$ on length scales $L$:
%
%
\begin{equation}\label{Yu2}
\operatorname{var}\Bigl[\sum\phi_T\eta_L\Bigr]
\lesssim
T\biggl(\frac{T}{L^d}\biggr)^{1/2-\delta}
\end{equation}
for $1\ll T\ll L^d$ and some arbitrarily small $\delta>0$,
see \cite{Yurinskii-86}, Lemma 2.4. Let us comment a bit on
the proof of (\ref{Yu2}): by stationarity of $\phi_T$, the
variance can be reformulated as a covariance, that is,%
\[
\operatorname{var}\Bigl[\sum\phi_T\eta_L\Bigr]
=
\operatorname{cov}\Bigl[\sum\phi_T\tilde\eta_L;\phi_T(0)\Bigr],
\]
with a modified averaging function $\tilde\eta_L$.
The starting point for (\ref{Yu2}) is to control the covariance by:
\begin{enumerate}[(ii)]
\item[(i)]
An additive decomposition of $\phi_T(0)$ over all finite subsets $S$ of
the lattice $\mathbb{Z}^d$, that is, $\phi_T(0)
=\sum_{S\subset\mathbb{Z}^d}\phi_{T,S}(0)$,
where $\phi_{T,S}(0)$ only depends on $a_{|S}$, that is,
the coefficients $a$ restricted to the subset $S$.
\item[(ii)] An estimate on how sensitively $\sum\phi_T\tilde\eta_L$
depends on $a_{|S}$.
\end{enumerate}
The decomposition in (i) is based on the probability measure on path space
$[0,\infty)\ni t\mapsto\eta(t)\in\mathbb{Z}^d$
describing the random walk generated by the operator $-\nabla^*\cdot
A\nabla$
(for a fixed realization of $a$). Indeed,
this probability measure on path space allows for a well-known representation
of $\phi_{T}(0)$ in terms of paths starting in $0$
(via the expected value). Hence,
the splitting can be obtained from restricting the expected
value to all paths $\eta$ with image $S$ (up to some exit time larger
than $T$),
see \cite{Yurinskii-86}, Lemma 2.3.

The sensitivity estimate (ii) comes in form of the
deterministic energy-type estimate
\[
\Bigl|\sum\phi_T\tilde\eta_L-\sum\tilde\phi_T\tilde\eta_L\Bigr|^2
\lesssim
\frac{T}{L^d}\sum_{\mathrm{edges}\ e\
\mathrm{s.}\mathrm{t.}\ e\cap S\not=\varnothing
}
\bigl(1+|\nabla\phi_T(e)|^2\bigr),
\]
where $\tilde\phi_T$ is the solution of
$T^{-1}\tilde\phi_T-\nabla^*\cdot\tilde A(\xi+\nabla\tilde\phi_T)=0$
with coefficients $\tilde A$ which differ from $A$ only on the subset $S$,
see \cite{Yurinskii-86}, (1.17).

The third ingredient for (\ref{Yu2}) is an estimate of
the probability that a path $\eta$ starting in $0$ crosses a given
edge $e$.
This
probability\vspace*{1pt} can be estimated in terms of the \textit{Green's function}
$G_T(x,0)$ of the operator $T^{-1}-\nabla^*\cdot A\nabla$
(where $x$ is one of the two sites on the edge $e$).
Yurinskii then appeals to estimates on
$G_T(x,y)$ that only depend on the ellipticity bounds $\alpha\le a\le
\beta$
of $A$ (and therefore do not depend on the realization of $a$)
see \cite{Yurinskii-86}, Lemma 2.1. As is well known, these type of
estimates rely on the \textit{Harnack inequality}.

Our variance estimate (\ref{PV12}) also relies on these
deterministic estimates of the Green's function $G_T(x,y)$,
see Lemma \ref{L11a}. However, our strategy to estimate a variance
differs substantially from Yurinskii's strategy of (i) and (ii).
As a matter of fact, with our methods, we could derive the optimal
variance estimate
%
%
\begin{equation}\label{Yu4}
\operatorname{var}\Bigl[\sum\phi_T\eta_L\Bigr]
\lesssim L^{2-d}
\end{equation}
for $L\gg1$. Estimate (\ref{Yu4}) is optimal in the sense that we
obtain the above scaling in the regime of ``vanishing ellipticity
ratio'' $1-\frac{\alpha}{\beta}\ll1$ by the arguments\vspace*{1pt} in the previous
subsection. Still, the optimal estimate (\ref{Yu4}) would not yield the
optimal estimate (\ref{PV11}) by Yurinskii's argument to pass from
(\ref{Yu2}) to (\ref{Yu1}).

Our strategy of estimating a variance is inspired by
an unpublished paper by Naddaf and Spencer \cite{Naddaf-Spencer-98}.
They use a \textit{spectral gap} estimate to control the
variance of some function $X$ of the coefficients
$\{a(e)\}_{\mathrm{edges}\ e}$
(i.e., a random variable):
%
%
\begin{equation}\label{Yu5}
\operatorname{var}[X] \lesssim
\biggl\langle\sum_{\mathrm{edges}\ e}
\biggl(\frac{\partial X}{\partial a(e)}
\biggr)^2
\biggr\rangle,
\end{equation}
see \cite{Naddaf-Spencer-98}, page 4. This type of estimate
can be seen as a Poincar\'e estimate with mean value zero w.r.t.
the infinite product measure that describes the distribution
of the coefficients (and the optimal constant in this estimate
is given by the smallest nonzero eigenvalue of the corresponding elliptic
operator, whence ``spectral gap''). Naddaf and Spencer derive (\ref
{Yu5}) via
the Brascamp--Lieb inequality for a large class of statistics for
$\{a(e)\}_{\mathrm{edges}\ e}$, which however does not
include all
i.i.d. statistics of $\{a(e)\}_{\mathrm{edges}\ e}$
considered by us. We therefore
rely on a slight modification of (\ref{Yu5}), see Lemma \ref{lem:var-estim}.

We also follow Naddaf and Spencer in the sense that we treat
the variance of
an \textit{energy density}. However, they express their result not in terms
of the energy density of $\phi_T$ but of a generic solution $u$
with a \textit{compactly supported, deterministic} right-hand side $f$,
that is,
%
%
\begin{equation}\label{Yu6}
-\nabla^*\cdot A\nabla u = \nabla^*\cdot f.
\end{equation}
Using (\ref{Yu6}),\vspace*{1pt} they obtain the formula
$\frac{\partial}{\partial a(e)}\sum\nabla u\cdot A\nabla u
=-|\nabla u(e)|^2$
so that an application of (\ref{Yu5}) yields the following
estimate on the energy density $X=\sum\nabla u\cdot A\nabla u$:
%
%
\begin{equation}\label{Yu7}
\operatorname{var}\Big[\sum\nabla u\cdot A\nabla u\Bigr]
\lesssim
\Bigl\langle\sum|\nabla u|^4\Bigr\rangle,
\end{equation}
see \cite{Naddaf-Spencer-98}, Proposition 1.

Naddaf and Spencer also remark
that provided the ellipticity contrast $1-\frac{\alpha}{\beta}$ is
small enough, \textit{Meyer's estimate} holds which states that
%
%
\begin{equation}\label{Yu8}
\sum|\nabla u|^4 \lesssim\sum|f|^4,
\end{equation}
with a constant that only depends on $\alpha$, $\beta$. The combination
of (\ref{Yu7}) and (\ref{Yu8}) yields the a priori estimate
%
%
\begin{equation}\label{Yu7bis}
\operatorname{var}\Bigl[\sum\nabla u\cdot A\nabla u\Bigr]
\lesssim\sum|f|^4,
\end{equation}
see \cite{Naddaf-Spencer-98}, Theorem 1. Since the left-hand side of
(\ref{Yu7bis}) scales as (volume)$^2$, while the right-hand side only
scales as $\mbox{volume}$, this estimate reveals the optimal decay of
fluctuations on the macroscopic level, very much like
(\ref{PV12}).---There is a somewhat theatrical convention in the
homogenization literature to call the lattice spacing $\e$ instead of 1
which highlights this scaling. Following Naddaf and Spencer, we use
Meyer's estimate, albeit applied on the Green's function $G_T(x,y)$,
see Lem\-ma~\ref{lem:int-grad}.

We will make use of the following notation:
\begin{itemize}
\item$d\geq2$ is the dimension;
\item$\int_{\Z^d} \,dx$ denotes the sum over $x\in\Z^d$, and $\int
_{D}\,dx$ denotes the sum over $x\in\Z^d$ such that $x\in D$, $D$ open
subset of $\R^d$;
\item$\langle\cdot\rangle$ is the ensemble average, or
equivalently the
expectation in the underlying probability space;
\item$\operatorname{var}[\cdot]$ is the variance associated with the
ensemble average;
\item$\lesssim$ and $\gtrsim$ stand for $\leq$ and $\geq$ up to a
multiplicative constant which only depends on the dimension $d$ and the
constants $\alpha,\beta$ (see Definition \ref{def:alpha-beta} below) if
not otherwise stated;
\item when both $\lesssim$ and $\gtrsim$ hold, we simply write $\sim$;
\item we use $\gg$ instead of $\gtrsim$ when the multiplicative
constant is (much) larger than~$1$;
\item$(\ee_1,\ldots,\ee_d)$ denotes the canonical basis of $\Z^d$.
\end{itemize}

\section{Main results}

\subsection{General framework}

\begin{defi}\label{def:alpha-beta}
We say that $a\dvtx\Z^d\times\Z^d\to\R^+,(x,y)\mapsto a(x,y)$ is a
conductivity function on $\Z^d$ if there exist $0<\alpha\leq\beta
<\infty$ such that:
\begin{itemize}
\item$a(x,y)=0$ if $|x-y|\neq1$,
\item$a(x,y)=a(y,x)\in[\alpha,\beta]$ if $|x-y|=1$.
\end{itemize}
We denote by $\calA$ the set of such conductivity functions.
\end{defi}
\begin{defi}
The elliptic operator $L\dvtx L^2_{\mathrm{loc}}(\Z^d)\to L^2_{\mathrm
{loc}}(\Z^d),u\mapsto Lu$ associated with a conductivity function
$a\in
\calA$ is defined for all $x\in\Z^d$ by
%
%
\begin{equation}\label{eq:def-elliptic}
(Lu)(x)=-\nabla^*\cdot A(x)\nabla u(x),
\end{equation}
where
\[
\nabla u(x):=\left[
\matrix{u(x+\ee_1)-u(x) \cr
\vdots\cr
u(x+\ee_d)-u(x)}
\right],\qquad
\
\nabla^* u(x):=\left[
\matrix{
u(x)-u(x-\ee_1) \cr
\vdots\cr
u(x)-u(x-\ee_d)}
\right]
\]
and
\[
A(x):=\operatorname{diag}[a(x,x+\ee_1),\ldots,a(x,x+\ee_d)].
\]
In particular, it holds that
\[
(Lu)(x)=\sum_{y,|x-y|=1}a(x,y)\bigl(u(x)-u(y)\bigr).
\]
If $a(x,y)=1$ for $|x-y|=1$, then the associated elliptic operator $L$
is the discrete Laplace operator, and is denoted by $-\bigtriangleup$.
\end{defi}
\begin{defi}[(Discrete integration by parts)]\label{def:ibp}
Let $d\geq2$, $h\in L^2(\Z^d)$ and $g\in L^2(\Z^d,\R^d)$. Then the
discrete integration by parts reads
\[
\int_{\Z^d}h(x)\nabla^*\cdot g(x) \,dx = -\int_{\Z^d} \nabla h(x)
\cdot g(x) \,dx.
\]
\end{defi}

We now turn to the definition of the statistics of the
conductivity function.
\begin{defi}
A conductivity function is said to be independent and identically
distributed (i.i.d.) if the coefficients $a(x,y)$ for $|x-y|=1$ are
i.i.d. random variables.
\end{defi}
\begin{defi}
The conductivity matrix $A$ is obviously stationary in the sense that
for all $z\in\Z^d$,
$A(\cdot+z)$ and $A(\cdot)$ have the same statistics; and for all
$x,z\in\Z^d$,
\[
\langle A(x+z) \rangle=\langle A(x) \rangle.
\]
Therefore, any translation invariant function of $A$, such as the modified
corrector $\phi_T$ (see Lemma \ref{lem:app-corr}), is jointly
stationary with $A$.
In particular, not only are $\phi_T$ and its gradient $\nabla\phi_T$
stationary, but also
any function of $A$, $\phi_T$ and $\nabla\phi_T$. A useful such
example is the energy density
$(\xi+\nabla\phi_T)\cdot A (\xi+\nabla\phi_T)$, which is stationary
by joint stationarity of $A$ and $\nabla\phi_T$.

Another translation invariant function of $A$ is the Green
functions $G_T$ of Definition \ref{def:Green}.
In this case, stationarity means that $G_T(\cdot+z,\cdot+z)$ has the
same statistics as $G_T(\cdot,\cdot)$ for all
$z\in\Z^d$, so that in particular, for all $x,y,z\in\Z^d$,
\[
\langle G_T(x+z,y+z) \rangle=\langle G_T(x,y) \rangle.
\]
\end{defi}
\begin{lemma}[(Corrector (\cite{Kunnemann-83}, Theorem 3))]\label{lem:corr}
Let $a\in\calA$ be an i.i.d. conductivity function, then for all $\xi
\in
\R^d$, there exists a unique random function $\phi\dvtx\Z^d\to\R$ which
satisfies the corrector equation
%
%
\begin{equation}\label{eq:corr}
-\nabla^*\cdot A(x) \bigl( \nabla\phi(x) +\xi\bigr)=0
\qquad\mbox{in }\Z^d,
\end{equation}
and such that $\phi(0)=0$, $\nabla\phi$ is stationary and $\langle
\nabla\phi\rangle=0$.
In addition, $\langle|\nabla\phi|^2 \rangle\lesssim|\xi|^2$.
\end{lemma}

We also define an ``approximation'' of the corrector as follows.
\begin{lemma}[(Approximate corrector (\cite{Kunnemann-83}, proof of
Theorem 3))]\label{lem:app-corr}
Let $a\in\calA$ be an i.i.d. conductivity function, then for all $T>0$
and $\xi\in\R^d$, there exists a unique stationary random function
$\phi
_{T}\dvtx\Z^d\to\R$ which satisfies the ``approximate'' corrector equation
%
%
\begin{equation}\label{eq:app-corr}
T^{-1}\phi_{T}(x)-\nabla^*\cdot A(x) \bigl( \nabla\phi_{T}(x) +\xi
\bigr)=0 \qquad\mbox{in }\Z^d,
\end{equation}
and such that $\langle\phi_{T} \rangle=0$.
In addition, $T^{-1} \langle\phi_{T}^2 \rangle+\langle{|\nabla
\phi_{T}|}^2 \rangle\lesssim|\xi|^2$.
\end{lemma}

Note that $\phi_T$ is stationary, whereas $\phi$ is not.
\begin{defi}[(Homogenized coefficients)]
Let $a\in\calA$ be an i.i.d. conductivity function and let $\xi\in
\R
^d$ and $\phi$ be as in Lemma \ref{lem:corr}.
We define the homogenized $d\times d$-matrix $A_\ho$ as
%
%
\begin{equation}\label{eq:fo-homog}
\xi\cdot A_\ho\xi= \langle{(\xi+\nabla\phi)\cdot A(\xi+\nabla
\phi)(0)} \rangle.
\end{equation}
\end{defi}

Note that (\ref{eq:fo-homog}) fully characterizes $A_\ho$
since $A_\ho$ is a symmetric matrix (it is in particular of the form
$a_\ho\Id$ for an i.i.d. conductivity function).

\subsection{Statement of the main result}

Our main result shows that the energy density
$\ener:=T^{-1}\phi_T^2+(\nabla\phi_T+\xi)\cdot A(\nabla\phi
_T+\xi)$
of the approximate corrector $\phi_T$, which is
a stationary scalar field, decorrelates
sufficiently rapidly so that smooth spatial averages
(defined with help of $\eta_L$) fluctuate
as they would if $\ener$ would be independent from site to site
(as is the case for the tensor field $A$ of the coefficients).
The strength of fluctuation is expressed in terms of the variance.
In more than two space dimensions (i.e., $d>2$),
the estimate does \textit{not} depend on the cut-off scale $\sqrt{T}$
and thus carries over to the energy density of the
corrector $\phi$. In two space dimensions, we are not able
to rule out a weak
(i.e., logarithmic) dependence on the cut-off scale $\sqrt{T}$:
\begin{theo}\label{th:main}
Let $a\in\calA$ be an i.i.d. conductivity function, and let $\phi$ and
$\phi_T$ denote the corrector and approximate correctors associated
with the conductivity function $a$ and direction $\xi\in\R^d$, $|\xi|=
1$. We then define for all $L>0$ and $T \gg1$ the symmetric matrix
$A_{L,T}$ characterized by
\[
\xi\cdot A_{L,T}\xi:=\int_{\Z^d} \bigl(T^{-1}\phi_T(x)^2+\bigl(\nabla
\phi
_T(x)+\xi\bigr)\cdot A(x)\bigl(\nabla\phi_T(x)+\xi\bigr)\bigr)\eta_L(x) \,dx,
\]
where $x\mapsto\eta_L(x)$ is an averaging function on $(-L,L)^d$ such
that $\int_{\Z^d}\eta_L(x)\,dx=1$ and $\|\nabla\eta_L\|_{L^\infty}
\lesssim L^{-d-1}$.
Then, there exists an exponent $q>0$ depending only on $\alpha,\beta$
such that
%
%
\begin{eqnarray}\label{eq:estim-var-hom}
&&\mbox{for }d=2\qquad \operatorname{var}[\xi\cdot A_{L,T}\xi] \lesssim
L^{-2}(\ln T)^q, \nonumber\\[-8pt]\\[-8pt]
&&\mbox{for }d>2\qquad \operatorname{var}[\xi\cdot A_{L,T}\xi] \lesssim
L^{-d}.\nonumber
\end{eqnarray}
In particular, for $d>2$, the variance estimate (\ref{eq:estim-var-hom})
holds for the energy density of the corrector $\phi$ itself.
\end{theo}
\begin{rem}\label{Remark1}
While it is natural to include the zero-order term
$T^{-1}\langle\phi_T^2\rangle$ into the definition of
the energy density, it is not essential for our result.
Here comes the reason: by a simplified version
of the string of arguments which lead to Theorem~\ref{th:main}
we can show that the variance of the zero-order
term is estimated as
\[
\operatorname{var}\biggl[\int_{\mathbb{Z}^d}\phi_T(x)^2\eta_L(x)\,dx\biggr]
\lesssim
\cases{
(\ln T)^q, &\quad for $d=2$,\cr
L^{2-d}, &\quad for $d>2$.}
\]
Hence, this term is of lower order in the regime
(of interest) $L\lesssim{T}$.
\end{rem}

The main ingredient to the proof of Theorem \ref{th:main} is
of independent interest. It states that all finite stochastic moments
of the approximate corrector $\phi_T$ are bounded independently
of $T$ for $d>2$ and grow at most logarithmically in $T$ for
$d=2$.
\begin{prop}\label{prop:main}
Let $a\in\calA$ be an i.i.d. conductivity function, $\xi\in\R^d$ with
$|\xi|=1$ and let $\phi_T$ denote the approximate corrector associated
with the conductivity function $a$, and $\xi$.
Then there exists a continuous function $\gamma\dvtx\R^+\to\R^+$ such that
for all $q\in\R^+$, there exists a constant $C_q$ such that for all $T>0$,
%
%
\begin{eqnarray}\label{eq:estim-prop}
&&\mbox{for }d=2\qquad \langle|\phi_T(0)|^{q} \rangle\leq C_q (\ln
T)^{\gamma(q)},
\nonumber\\[-8pt]\\[-8pt]
&&\mbox{for }d>2\qquad \langle|\phi_T(0)|^q \rangle\leq C_q.\nonumber
\end{eqnarray}
In addition, $\gamma(2n)=n(n+1)$ for all $n=2^l$, $l\in\N$ large enough.
\end{prop}

Let us give a heuristic argument for the behavior of $\langle|\phi
_T(0)|^q \rangle$
for $d=1$. In this case, for $T=\infty$, the gradient of the corrector
associated with $\xi=1$
is explicitly given by
\[
\nabla\phi= \frac{1}{a\langle a^{-1} \rangle}-1.
\]
Hence, $\phi(x)\in\R$ behaves as a discrete Brownian motion in $x\in
\Z$
once we have
fixed its value at 0.
Usually, one imposes $\phi(0)=0$ almost surely, so that for $|x|\sim
\sqrt{T}$,
\[
\langle|\phi(x)|^q | \phi(0)=0 \rangle\sim\bigl(\sqrt{T}\bigr)^{q/2}.
\]
Yet, one may choose a nontrivial initial value. In particular, one may
also consider $\phi(0)=\phi_T(0)$ (which yields a corrector field
different from the one in Definition~\ref{lem:corr}).
With $\phi$ defined this way, $\phi_T(x)$ and $\phi(x)$ are expected
not to differ much provided $|x|\ll\sqrt{T}$.
On the one hand, from this we deduce that $\phi_T(x)$ behaves locally
as a discrete Brownian motion starting at $\phi_T(0)$, so that we have
as above
\[
\langle|\phi_T(x)-\phi_T(0)|^q \rangle\sim|x|^{q/2}
\]
for all $q>0$ and $|x|\ll\sqrt{T}$.
On the other hand, since $\phi_T$ is stationary,
\[
\langle|\phi_T(x)-\phi_T(0)|^q \rangle\lesssim\langle|\phi
_T(x)|^q \rangle+\langle|\phi_T(0)|^q \rangle= 2 \langle|\phi
_T(0)|^q \rangle.
\]
These two estimates indeed suggest that
\[
\langle|\phi_T(0)|^q\rangle\gtrsim\sqrt{T}{}^{q/2-},
\]
where the minus sign accounts for the fact that the argument only holds
for $|x|\ll\sqrt{T}$---we may for instance miss logarithmic
corrections. Hence, there is a transition between unboundedness and
boundedness in $T$ for some $d\in(1,3)$. The linearization of the
problem in the regime of vanishing ellipticity contrast, that is,
$1-\frac{\alpha}{\beta}\ll1$, suggests that $d=2$ is indeed the
critical dimension for Proposition~\ref{prop:main}, that is, the
dimension where a logarithmic behavior is to be expected. However,
there is no reason why $d=2$ should be critical for Theorem
\ref{th:main}. Indeed, in the case of $d=1$, the statement of Theorem
\ref{th:main} holds without a logarithm.

In view of our discussion of the case $d=1$ and the observations in case
of vanishing ellipticity contrast, it is not surprising that the
statement of bounded stochastic moments is harder to prove the closer
we are to $d=2$.
For the experts in homogenization, let us give a quick sketch of the
strategy of the proof of this result.
Independent of the dimension, the proof always starts from the variance
estimate (Lemma \ref{lem:var-estim}) applied to $\phi_T(0)^q$
and makes use of the representation of $\frac{\partial\phi
_T(0)}{\partial a(e)}$
with help of the gradient $\nabla_xG_T(x,0)$ (Lemma \ref{lem:diff-phi}).
\begin{itemize}
\item In the case of $d>4$, the uniform pointwise, but suboptimal,
decay $|\nabla_x G_T(x$, $y)|\lesssim|x-y|^{d-2}$, which can be easily
obtained from the
same pointwise decay of the Green's function itself, is sufficient.
\item In case $d=4$, it would be enough to appeal to the H\"older
estimate (with exponent $\gamma$ only depending on the ellipticity
contrast) in order to get the
somewhat better pointwise decay $|\nabla_x G_T(x,y)|\lesssim
|x-y|^{d-2-\gamma}$.
\item In $d=3$, we need (in addition) the \textit{optimal} decay $|\nabla
_x G_T(x,y)|\lesssim|x-y|^{d-1}$, which cannot be a pointwise control,
but only an average control on dyadic annuli.
In fact, we need the control of the \textit{square} average, which we
easily obtain from the Cacciopoli estimate.
\item For $d=2$, the square average is not sufficient anymore, we need
the average to some power $p>2$, as provided by Meyers' estimate
(Lemma \ref{lem:int-grad}). This forces us---somewhat
counterintuitively---to first estimate high moments of $\phi_T$,
so that the exponent we put on the gradient of the
Green's function can be chosen close to 2 (and thus below Meyers' exponent).
\end{itemize}
In this presentation, we only display the last strategy (although it is
an overkill for dimensions $d>2$).

As a corollary of Proposition \ref{prop:main}, we obtain the
following existence and
uniqueness result of stationary solutions to the corrector equation
(\ref{PV2}) for $d>2$, which settles
a long-standing open question.
\begin{corollary}\label{prop:corrector}
Let $a\in\calA$ be an i.i.d. conductivity function.
Then, for $d>2$ and for all $\xi\in\R^d$, there exists a unique
stationary random field $\phi$
such that $\langle\phi\rangle=0$ and
\[
-\nabla^*\cdot\bigl(A(\xi+\nabla\phi)\bigr) = 0 \qquad\mbox{in }\Z^d.
\]
In addition, $\langle\phi^2+|\nabla\phi|^2 \rangle\lesssim|\xi|^2$.
\end{corollary}

We will not prove Corollary \ref{prop:corrector} in detail.
Here comes the argument.
Proposition~\ref{prop:main} yields the a priori estimate $\langle
\phi_T^2 \rangle<C$ which is uniform in $T$.
This additional estimate allows us to pass to the limit in the
probability space for $\phi_T$,
as it is done for $\nabla\phi_T$ in \cite{Kunnemann-83}, proof of
Theorem 3.
Note that the corrector fields of Lem\-ma~\ref{lem:corr} and
Theorem \ref
{prop:corrector} do not
coincide (only their gradients coincide).
Uniqueness further requires the argument by Papanicolaou and Varadhan in
\cite{Papanicolaou-Varadhan-79}, which
does not appear in \cite{Kunnemann-83}.

Let us point out that Proposition \ref{prop:main},
Theorem \ref{th:main} and Corollary \ref{prop:corrector} hold true for
more general distributions, provided the variance estimate of
Lemma \ref{lem:var-estim} below holds.
In particular, the law of $a(x,x+\ee_i)$ may depend on the direction
$\ee_i$, which would give a general diagonal homogenized matrix (not
necessarily a multiple of the identity matrix). More generally,
$a(x,x')$ and $a(y,y')$ may also be slightly correlated. We do not
pursue this direction in this article.

\subsection{Structure of the proof and statement of the auxiliary results}

Not surprisingly, in order to control the variance of some function
$X$ of the coefficients $a$ (like the spatial average of the
energy density of the approximate corrector $\phi_T$),
one needs to control the gradient
of $X$ w.r.t. $a$.
As in \cite{Naddaf-Spencer-98}, this is quantified
by the following general variance estimate:
\begin{lemma}[(Variance estimate)]\label{lem:var-estim}
Let $a=\{a_i\}_{i\in\mathbb{N}}$ be a sequence of i.i.d. random
variables with range $[\alpha,\beta]$.
Let $X$ be a Borel measurable function of $a\in\mathbb{R}^\mathbb{N}$
(i.e., measurable
w.r.t. the smallest $\sigma$-algebra on $\mathbb{R}^\mathbb{N}$
for which all coordinate functions $\mathbb{R}^\mathbb{N}\ni a\mapsto
a_i\in\mathbb{R}$
are Borel measurable, cf. \cite{Klenke-08}, Definition 14.4).

Then we have
%
%
\begin{equation}\label{eq:var-estim}
\operatorname{var}[X] \le\Biggl\langle
\sum_{i=1}^\infty\sup_{a_i}
\biggl|\frac{\partial X}{\partial a_i}\biggr|^2\Biggr\rangle
\operatorname{var}[a_1],
\end{equation}
where $\sup_{a_i}|\frac{\partial X}{\partial a_i}|$ denotes
the supremum of the modulus of
the $i$th partial derivative
\[
\frac{\partial X}{\partial a_i}(a_1,\ldots,a_{i-1},a_i,a_{i+1},\ldots)
\]
of $X$
with respect to the variable $a_i\in[\alpha,\beta]$.
\end{lemma}
\begin{rem} Let us comment a bit on Lemma \ref{lem:var-estim}.
Estimate (\ref{eq:var-estim}) is a weakened version of a spectral
gap estimate
%
%
\begin{equation}\label{5.-1}
\operatorname{var}[X] \lesssim\Biggl\langle
\sum_{i=1}^\infty
\biggl|\frac{\partial X}{\partial a_i}\biggr|^2\Biggr\rangle,
\end{equation}
which already played a central role in Naddaf and Spencer's analysis of
stochastic homogenization \cite{Naddaf-Spencer-98}, Section 2.
We note that for i.i.d. random variables, such a spectral
gap estimate (\ref{5.-1}) follows ``by tensorization'' from the one-dimensional
spectral gap estimate
%
%
\begin{equation}\label{5.-2}
\langle X(a_1)^2\rangle-\langle X(a_1)\rangle^2 \lesssim
\biggl\langle\biggl|\frac{\partial X}{\partial a_1}\biggr|^2
\biggr\rangle,
\end{equation}
see, for instance, \cite{Ledoux-01}, Lemma 1.1. The one-dimensional
spectral gap estimate (\ref{5.-2}) holds under mild assumptions on the
distribution of $a_1$. Yet, (\ref{5.-2}) does not hold for atomic
measures like $\langle X(a_1)\rangle=\frac{1}{2}(X(1)+X(2))$. Since
Lemma \ref{lem:var-estim} covers the case of atomic measures, we only
obtain the weaker form (\ref {eq:var-estim}) of (\ref{5.-1}). Despite
this technical detail, the proof of Lemma \ref{lem:var-estim} is very
similar to the one in~\cite{Ledoux-01}, Lemma 1.1.
\end{rem}

As in \cite{Naddaf-Spencer-98}, in the proof of Theorem \ref {th:main},
we will make use of the fact that
$T^{-1}\phi_T^2+(\nabla\phi_T+\xi)\cdot A(\nabla\phi_T+\xi)$ is an
energy density, which yields the following elementary formula for the
partial derivative w.r.t. the value $a(e)$ of the coefficient in the
edge $e=[z,z+\ee_i]$:
%
%
\begin{eqnarray}\label{I.1}
&&\frac{\partial}{\partial a(e)}\int
\bigl(T^{-1}\phi_T^2+(\nabla\phi_T+\xi)\cdot A(\nabla\phi_T+\xi
)\bigr)(x)\eta
_L(x) \,dx
\nonumber\\
&&\qquad=
-2\int\biggl(\frac{\partial\phi_T}{\partial a(e)}
\nabla\eta_L\cdot A(\nabla\phi_T+\xi)\biggr)(x)\,dx
\\
&&\qquad\quad{}
+\bigl(\eta_L (\nabla_i\phi_T+\xi_i)\bigr)^2(z),\nonumber
\end{eqnarray}
up to minor modifications coming from the discrete Leibniz rule,
see Step 1 of the proof of Theorem \ref{th:main}.

This formula makes the gradient of the averaging function
$\eta_L$ appear;
in order to benefit from this, we assume that the averaging function
is smooth so that we get an extra power of $L^{-1}$.
The merit of (\ref{I.1}) is
that we need to control the partial derivative
$\frac{\partial\phi_T(x)}{\partial a(e)}$ of the approximate corrector
$\phi_T(x)$ (and not of its spatial derivatives).
Not surprisingly, this partial derivative involves the Green's function
$G_T(x,\cdot)$. More precisely, it involves the gradient
$\nabla_{z_i} G_T(x,z)$ of the
Green's function with singularity in $z$
[and not its second gradient $\nabla_{z_i} \nabla_x G_T(x,z)$, for
which we would \textit{not} have the optimal decay rate uniformly in $a$].

We define discrete Green's functions as follows.
\begin{defi}[(Discrete Green's function)]\label{def:Green}
Let $d\geq2$.
For all $T>0$, the Green's function $G_T\dvtx\calA\times\Z^d\times\Z
^d\to\Z^d,(a,x,y)\mapsto G_T(x,y;a)$ associated with the conductivity
function $a$ is defined for all $y\in\Z^d$ and $a\in\calA$ as the
unique solution in $L^2_x(\Z^d)$ to
%
%
\begin{eqnarray}\label{eq:disc-Green}\qquad
&&\int_{\Z^d}T^{-1}G_T(x,y;a)v(x) \,dx\nonumber\\[-8pt]\\[-8pt]
&&\qquad{}+\int_{\Z^d}\nabla v(x)\cdot
A(x)\nabla_x G_T(x,y;a) \,dx=v(y)\qquad \forall v\in
L^2(\Z^d),\nonumber
\end{eqnarray}
where $A$ is as in (\ref{eq:def-elliptic}).
\end{defi}

Note that the existence and uniqueness of discrete Green's
functions is a consequence of Riesz' representation theorem.
Throughout this paper, when no confusion occurs, we use the short-hand
notation $G_T(x,y)$ for $G_T(x,y;a)$.

The following lemma provides the elementary formula
relating the ``susceptibility'' $\frac{\partial\phi_T(x)}{\partial a(e)}$
of $\phi_T(x)$ to the Green's function $G_T(x,y)$.
\begin{lemma}\label{lem:diff-phi}
Let $a\in\calA$ be an i.i.d. conductivity function, and let $G_T$ and
$\phi_T$ be the associated Green's function and approximate corrector
for $T>0$ and $\xi\in\R^d$, $|\xi|=1$.
Then, for all $e=[z,z+\ee_i]$ and $x\in\Z^d$,
%
%
\begin{equation}\label{eq:diff-phi-1}
\frac{\partial\phi_T(x;a)}{\partial a(e)}=-\bigl(\nabla_i \phi
_T(z;a)+\xi
_i\bigr) \nabla_{z_i} G_T(z,x;a),
\end{equation}
and for all $n\in\N$,
%
%
\begin{eqnarray}\label{eq:diff-phi-2}
&&
\sup_{a(e)} \biggl|\frac{\partial}{\partial a(e)} [\phi
_T(x;a)^{n+1}]\biggr|
\nonumber\\
&&\qquad \lesssim |\phi_T(x;a)|^{n}\bigl(|\nabla_i \phi_T(z;a)|+1\bigr)|\nabla_{z_i}
G_T(z,x;a)|\\
&&\qquad\quad{} + \bigl(|\nabla_i \phi_T(z;a)|+1\bigr)^{n+1}|\nabla_{z_i}
G_T(z,x;a)|^{n+1}.\nonumber
\end{eqnarray}
In addition, it holds that
%
%
\begin{equation}\label{eq:diff-phi-3}
{\sup_{a(e)}}|\nabla_i \phi_T(z;a)| \lesssim|\nabla_i \phi_T(z;a)|+1.
\end{equation}
Note that the multiplicative constant in (\ref{eq:diff-phi-2}) depends
on $n$ next to $\alpha$, $\beta$ and $d$.
\end{lemma}

In addition,\vspace*{1pt} Lemma \ref{lem:diff-phi} provides uniform
estimates on
$\frac{\partial[\phi_T(x)^n]}{\partial a(e)}$ in $a(e)$
(the case $n>1$ is needed in Proposition \ref{prop:main}). In order
to obtain this uniform control in $a(e)$, we need
to control $\nabla_z G(z,x;a)$ uniformly in $a(e)$. Again, this comes
from considering $\frac{\partial\nabla_z G(z,x;a)}{\partial a(e)}$.
The following lemma provides the elementary formula for
$\frac{\partial\nabla_z G(z,x;a)}{\partial a(e)}$ and a uniform estimate
in $a(e)$.
\begin{lemma}\label{lem:diff-Green}
Let $G_T\dvtx\calA\times\Z^d\times\Z^d \to\R, (a,x,y)\mapsto
G_T(x,y;a)$ be the Green's function associated with the conductivity
function $a$ for $T>0$.
For all $e=[z,z+\ee_i]$ and for all $x,y\in\Z^d$, it holds that
%
%
\begin{equation}\label{eq:diff-Green}
\frac{\partial}{\partial a(e)}G_T(x,y;a) = -\nabla_{z_i}G_T(x,z;a)
\nabla_{z_i}G_T(z,y;a).
\end{equation}
As a by-product, we also have: for all $x\in\Z^d$
%
%
\begin{equation}\label{eq:bd-G(x,e)}
{\sup_{a(e)}}|\nabla_{z_i}G_T(z,x;a)| \lesssim
|\nabla
_{z_i}G_T(z,x;a)|.
\end{equation}
\end{lemma}

There is a technical difficulty arising from the fact that
$a$ has infinitely many components. In Lemma \ref{lem:var-estim}, this
technical difficulty is handled by the strong measurability assumptions
on $X$. The following lemma establishes these measurability properties for
$\phi_T$, so that we can apply Lemma \ref{lem:var-estim}.
\begin{lemma}\label{lem:depend-coeff}
Let $a\in{\mathcal A}_{\alpha\beta}$ be an i.i.d. conductivity
function, and let $G_T(\cdot,\cdot;a)$ and $\phi_T(\cdot;a)$
be the associated Green's function and approximate corrector for $\xi
\in
\R^d$, $d\ge2$, and $T>0$.
Then for fixed
$x,y\in\mathbb{Z}^d$,
$G_T(x,y,\cdot)$ and $\phi_T(x;\cdot)$ are continuous w.r.t. the product
topology of ${\mathcal A}_{\alpha\beta}$ (i.e.,\vspace*{1pt}
the smallest/coarsest topology on $\mathbb{R}^E$, where $E$ denotes the
set of edges,
such that the coordinate functions $\mathbb{R}^{E}\ni a\mapsto a_e\in
\mathbb{R}$
are continuous for all edges $e\in E$).

In particular, $G_T(x,y;\cdot)$ and $\phi_T(x;\cdot)$ are
Borel measurable
functions of $a\in{\mathcal A}_{\alpha\beta}$, so that one may apply
Lemma \ref{lem:var-estim}
to $\phi_T(x;\cdot)$ and nonlinear funtions thereof.
\end{lemma}

The proof of Theorem \ref{th:main} crucially relies on the
fact that
$\phi_T$ is almost bounded independently of $T$ (in $d>2$).
More precisely, it relies on the fact that any moment
$\langle\phi_T(0)^n\rangle$ is bounded independently of $T$
as stated in Proposition \ref{prop:main}. Starting point for
Proposition \ref{prop:main} is again Lemma \ref{lem:var-estim}, which
is iteratively applied
to $\phi_T(0)^m$ where $m$ increases dyadically. This is how
Lemma \ref{lem:diff-phi} comes in again. However, the crucial gain in
stochastic integrability is provided by the following lemma.
It can be interpreted as a Cacciopoli estimate in probability and relies
on the stationarity of $\phi_T$.
\begin{lemma}\label{lem:Lp-grad}
Let $a\in\calA$ be an i.i.d. conductivity function, and let $\phi_T$ be
the approximate corrector associated with the coefficients $a$ for $\xi
\in\R^d$, $|\xi|=1$.
Then for $d\geq2$ and for all $n \in2\mathbb{N}$,
%
%
\begin{equation}\label{eq:phi-grad-phi}
\bigl\langle|\phi_T(0)|^n\bigl(|\nabla\phi_T(0)|^2+|\nabla^* \phi_T(0)|^2\bigr)
\bigr\rangle
\lesssim\langle|\phi_T|^n(0) \rangle,
\end{equation}
where the multiplicative constant does depend on $n$ next to $\alpha$,
$\beta$, and $d$, but not on $T>0$.
\end{lemma}

In order to prove Proposition \ref{prop:main} via Lemma
\ref{lem:var-estim} [applied to $\phi_T(0)^n$] and
Lem\-ma~\ref{lem:diff-phi}, we need some \textit{weak version} of the
optimal decay of the gradient $\nabla_z G_T(x,z)$ of Green's function
in $|x-z|$, that is,
%
%
\begin{equation}\label{I.2}
|\nabla_z G_T(x,z;a)| \lesssim|x-z|^{1-d}
\qquad\mbox{uniformly in } a \mbox{ and } T.
\end{equation}
This decay is the best we can hope as can be checked on the Green
function for the Laplace equation.
The same decay property is needed to prove Theorem \ref{th:main}
via Lemma \ref{lem:var-estim} [applied to (\ref{I.1})] and Lemma
\ref{lem:diff-phi}.
Yet it is well known from the continuum case that there are
no \textit{pointwise} in $z$ bounds of the type (\ref{I.2})
which would hold uniformly in the ellipticity constants
$\alpha$, $\beta$. (An elementary argument shows that any bound on
$\nabla_x G(x,y)$ which would be
uniform in $a$ and in $1/2\leq|x-y|\leq1$ would yield that a bounded
$a$-harmonic
function has bounded gradient. However, for $d=2$ and for any $\gamma
>0$, there are examples
of $a$-harmonic functions from the theory of quasi-conformal mappings
that are not H\"older continuous with exponent~$\gamma$,
see \cite{Gilbarg-Trudinger-98}, Section 12.1.)
Nevertheless, (\ref{I.2}) holds in the
\textit{square averaged} sense on dyadic annuli, as can be seen by
a standard Cacciopoli argument based on the optimal decay of the
Green's function itself, that is,
%
%
\begin{equation}\label{I.3}
G_T(x,z) \lesssim|x-z|^{2-d} \qquad\mbox{uniformly in } a \mbox
{ and } T,
\end{equation}
in the case $d>2$. The pointwise estimate (\ref{I.3}) in $x$ and
$z$ is a classical
result \cite{Grueter-Widman-82}, Theorem 1.1, that relies on Harnack's
inequality.
It has been partially extended to discrete settings, see
in particular the Harnack inequality on graphs \cite{Delmotte-97}.
However, we
did not find a suitable reference for the BMO-type estimate
in the case of $d=2$. On the other hand, we do not require the \textit{pointwise}
version of (\ref{I.3}), but just an averaged version
on dyadic annuli. The statements
we need are collected in the following lemma.
\begin{lemma}\label{L11a}
Let $a\in\calA$, $T>0$ and $G_T$ be the associated Green's function.
For all $d\ge2$ and $q\geq1,r\geq0$,
\begin{enumerate}[(iii)]
\item[(i)] BMO and $L^q$ estimate: for all $R \gg1$,
%
%
\begin{eqnarray}
\label{11.16}
&&\mbox{for } d=2\qquad  \int_{|x-y|\le R}\bigl|G_T(x,y)-\bar G_T(\cdot,y)_{\{
|x-y|\le R\}}\bigr|^q \,dx
\lesssim R^2,\\
\label{11.16b}
&&\mbox{for }d>2\qquad  \int_{R\le|x-y|\le2R}G_T(x,y)^q \,dx \lesssim
R^d (R^{2-d})^q,
\end{eqnarray}
where $\bar G_T(\cdot,y)_{\{ |y-x|\le R\}}$ denotes the average of
$G_T(\cdot,y)$ over the ball $\{x\in\Z^d,|x-y|\le R\}$.
\item[(ii)] Behavior for $R\sim\sqrt{T}$ and $d=2$:
%
%
\begin{equation}\label{L11.19b}
R^{-2}\int_{|x-y|\leq R}G_T(x,y)^2 \,dx \lesssim1 .
\end{equation}
\item[(iii)] Decay at infinity: for all $R\geq\sqrt{T}$,
%
%
\begin{equation}\label{11.17}
\int_{R\le|x-y|\le2R}G_T(x,y)^q \,dx
\lesssim
R^d (R^{2-d})^q\bigl(\sqrt{T}R^{-1}\bigr)^r.
\end{equation}
\end{enumerate}
The multiplicative constants in (\ref{11.16}), (\ref{11.16b}) and
(\ref
{11.17}) depend on $q,r$ next to $\alpha$, $\beta$ and $d$.
\end{lemma}

We present a proof of Lemma \ref{L11a} which for $d=2$ is a
direct version
of the indirect argument developed in \cite
{Dolzmann-Hungerbuehler-Mueller-00}, Lemma
2.5, in case of
a nonlinear, continuum equation. For the convenience of the reader,
we also include the proof for $d>2$---anyway, it has the same
building blocks as the argument for $d=2$. This makes our paper
self-contained w.r.t. the properties of $G_T$.

However, it is not quite enough to know (\ref{I.2}) in
the \textit{square}-averaged sense on dyadic annuli. In order
to compensate for the fact that we only control \textit{finite}
stochastic moments of $\nabla\phi_T(0)$ via Proposition \ref
{prop:main}, we need to control
a $p$th power of the gradient $\nabla_z G_T(x,z)$ of Green's function
in the optimal way for some $p>2$. This slight increase in
integrability is provided by Meyers' estimate, which yields such
a $p>2$ as a function of the ellipticity bounds $\alpha$, $\beta$ only.
Meyers' estimate has already been crucially used in \cite{Naddaf-Spencer-98},
however in a somewhat different spirit. There it is used that for
sufficiently small ellipticity contrast,
$1-\frac{\alpha}{\beta}\ll1$, one has $p\ge4$. The following lemma
is the version of Meyers' estimate we need and will prove.
\begin{lemma}[(Higher integrability of gradients)]\label{lem:int-grad}
Let $a\in\calA$ be a conductivity function, and $G_T$ be its associated
Green's function. Then, for $d\geq2$, there exists ${p}>2$ depending
only on $\alpha,\beta$, and $d$ such that for all $T>0$, ${p}\geq
q\geq2$,
$k >0$ and $R\gg1$,
%
%
\begin{equation} \label{eq:int-grad}
\int_{R\leq|z|\leq2R}|\nabla_zG_T(z,0) |^q\,dz
\lesssim R^d(R^{1-d})^q
\min\bigl\{1,\sqrt{T}R^{-1}\bigr\}^{k}.
\end{equation}
\end{lemma}

For technical reasons, we need a \textit{pointwise} decay of
$G_T(x,y;a)$ in $|x-y|$ uniformly in $a$ (but not in $T$). The decay we
obtain is suboptimal and easily follows from Lemmas \ref{L11a} and
\ref
{lem:int-grad} using the discreteness.
\begin{corollary}\label{coro:unif-bound}
For all $d\ge2$ and $T>0$, there exists a bounded radially symmetric
function $h_T\in L^1(\Z^d)$ depending only on $d,\alpha,\beta$, and $T$
such that
\[
G_T(x,y;a) \leq h_T(x-y)
\]
for all $x,y\in\Z^d$ and $a\in\calA$.
\end{corollary}

Lemmas \ref{L11a} and \ref{lem:int-grad} only treat $G_T$ away
from the diagonal $x=y$---which is a consequence of the fact that the
scaling symmetry is broken by the discreteness. Using the discreteness,
the following corollary establishes a bound independent of $T$ and $a$.
\begin{corollary}\label{coro:unif-bound-grad}
For all $a\in\calA$, $T>0$ and $x,y\in\Z^d$,
\[
|\nabla G_T(x,y;a)| \lesssim1.
\]
\end{corollary}

Finally, for the proof of Theorem \ref{th:main}, we need to know
that also the \textit{convolution} of the gradients of the Green's
functions decays at the optimal rate, that is,
%
%
\begin{eqnarray}\label{I.4}
&&\int_{\Z^d}|\nabla_z G_T(x,z)||\nabla_z
G_T(x',z)|\,dz\nonumber\\[-8pt]\\[-8pt]
&&\qquad\lesssim|x-x'|^{2-d} \qquad\mbox{uniformly in } a \mbox{ and }
T.\nonumber
\end{eqnarray}
As for (\ref{I.2}), it is not necessary
to know (\ref{I.4}) \textit{pointwise} in $(x,x')$, but only in an
averaged sense on dyadic annuli. The following lemma shows that
(\ref{I.4}) for linear averages can be inferred from (\ref{I.2})
for quadratic averages.
\begin{lemma}\label{lem:hh}
Let $h_T\in L^2_{\mathrm{loc}}(\Z^d)$ be such that for all $R\gg1$
and $T>0$,
%
%
\begin{eqnarray}
\label{eq:assump-h-d=2}
&&\mbox{for }d=2\qquad \int_{R< |z| \leq2R}h_T^2(z)\,dz\lesssim\min\bigl\{
1,\sqrt
{T}R^{-1}\bigr\}^2,\\
\label{eq:assump-h}
&&\mbox{for }d>2\qquad \int_{R< |z|\leq2R}h_T^2(z)\,dz \lesssim R^{2-d},
\end{eqnarray}
and for $R\sim1$
%
%
\begin{equation}\label{eq:assump-h-R=1}
\mbox{for }d\geq2\qquad \int_{ |z|\leq R}h_T^2(z)\,dz \lesssim1.
\end{equation}
Then for $R\gg1$
%
%
\begin{eqnarray}
\label{eq:res-h-d=2}
&&\mbox{for }d=2\qquad \int_{|x|\leq R} \int_{\Z^d} h_T(z)h_T(z-x) \,dz
\,dx\nonumber\\[-8pt]\\[-8pt]
&&\hphantom{\mbox{for }d=2}\qquad\qquad\lesssim R^2 \max\bigl\{1,\ln\bigl(\sqrt{T}R^{-1}\bigr)\bigr\},\nonumber\\
\label{eq:res-h}
&&\mbox{for }d>2\qquad \int_{|x|\leq R} \int_{\Z^d} h_T(z)h_T(z-x) \,dz \,dx
\lesssim R^2.
\end{eqnarray}
\end{lemma}

We present the proof of Proposition \ref{prop:main} and
Theorem \ref{th:main} in Section \ref{sec:proof}.
We gather in Section \ref{sec:proofs-Green} the proofs of the decay
estimates for the discrete Green functions (i.e., Lemmas \ref{L11a}
and \ref{lem:int-grad},
and Corollaries \ref{coro:unif-bound} and \ref{coro:unif-bound-grad})
since they are needed at several
places in the paper, and may be of independent interest.
The proofs of the remaining auxiliary lemmas are the object of
Section \ref{sec:proofs-aux}.\\


\section{Proofs of the main results}\label{sec:proof}

\subsection{\texorpdfstring{Proof of Proposition \protect\ref{prop:main}}{Proof of Proposition 2.1.}}

Starting point are Lemmas \ref{lem:var-estim} and \ref
{lem:depend-coeff}, which yield
\[
\operatorname{var}[\phi_T(0)^m] \lesssim\sum_e \biggl\langle\sup_{a(e)}
\biggl|\frac{\partial\phi_T(0)^{m}}{\partial a(e)}\biggr|^2 \biggr\rangle,
\]
where $\sum_e$ denotes the sum over the edges.
Using now (\ref{eq:diff-phi-2}) in Lemma \ref{lem:diff-phi}, this
inequality turns into
\begin{eqnarray*}
\operatorname{var}[\phi_T(0)^m] &\lesssim& \int_{\Z^d}\sum
_{i=1}^d \bigl\langle
\phi
_T(0)^{2(m-1)}\bigl(|\nabla_i\phi_T(z)|+1\bigr)^2|\nabla_{z_i}G_T(z,0)|^{2} \\
&&\hspace*{63pt}{} +\bigl(|\nabla_i\phi_T(z)|+1\bigr)^{2m}|\nabla
_{z_i}G_T(z,0)|^{2m}\bigr\rangle \,dz ,
\end{eqnarray*}
where we have replaced the sum over edges $e$ by the sum over sites
$z\in\Z^d$ and directions $\ee_i$ for $i\in\{1,\ldots,d\}$
according to
$e=[z,z+\ee_i]$.
Simplifying further, we obtain
%
%
\begin{eqnarray}\label{eq:var-use}
\operatorname{var}[\phi_T(0)^m]&\lesssim& \int_{\Z^d}\bigl\langle\phi
_T(0)^{2(m-1)}\bigl(|\nabla\phi_T(z)|+1\bigr)^2|\nabla_{z}G_T(z,0)|^{2}
\nonumber\\[-8pt]\\[-8pt]
&&\hspace*{47.6pt}{} +\bigl(|\nabla\phi_T(z)|+1\bigr)^{2m}|\nabla
_{z}G_T(z,0)|^{2m}\bigr\rangle \,dz.\nonumber
\end{eqnarray}

We proceed in four steps. Assuming first that for $n$ big enough
and for all $m\leq n$ it holds that
%
%
\begin{eqnarray}\label{eq:key-prop-d=2}
&& \int_{\Z^d} \bigl\langle\phi_T(0)^{2(m-1)}\bigl(|\nabla\phi
_T(z)|+1\bigr)^2|\nabla_zG_T(z,0)|^{2}\nonumber\\
&&\qquad\hspace*{25.5pt}{} +\bigl(|\nabla\phi_T(z)|+1\bigr)^{2m}|\nabla
_zG_T(z,0)|^{2m} \bigr\rangle \,dz
\\
&&\qquad \lesssim\bigl(\langle\phi_T(0)^{2n} \rangle^{{m/n}-{1}/({n(n+1)})} +1\bigr)
\cases{
\ln T, &\quad for $d=2$, \cr
1, &\quad for $d>2$,}\nonumber
\end{eqnarray}
we prove the claim in the first step.
The last three steps are dedicated to the proof of (\ref
{eq:key-prop-d=2}) for $n$ large enough.

\textit{Step} 1. Proof that (\ref{eq:var-use}) and (\ref{eq:key-prop-d=2})
imply (\ref{eq:estim-prop}).

For notational convenience, we set $\mu_d(T)=1$ for $d>2$ and
$\mu_d(T)=\ln T$ for $d=2$.
Let $n=2^l$, $l\in\N^*$. Using the elementary fact that
\[
\langle\phi_T(0)^{2m} \rangle\leq\langle\phi_T(0)^m \rangle
^2+\operatorname{var}[\phi_T(0)^m],
\]
from the cascade of inequalities (\ref{eq:var-use}) and (\ref
{eq:key-prop-d=2}) for $m=2^{l-q}$, $q\in\{0,\ldots,l\}$, we deduce
\begin{eqnarray}
\langle\phi_T(0)^{2\cdot2^l} \rangle&{\lesssim} &\langle\phi
_T(0)^{2^l} \rangle^2
+\mu_d(T)\bigl(\langle\phi_T(0)^{2n} \rangle
^{1-{1}/({n(n+1)})}+1\bigr)\nonumber\\
&&\eqntext{(\mbox{estimate }0)}\\
&\vdots& \nonumber\\
\langle\phi_T(0)^{2\cdot2^{l-q}} \rangle&{\lesssim}& \langle\phi
_T(0)^{2^{l-q}} \rangle^2
+\mu_d(T)\bigl(\langle\phi_T(0)^{2n} \rangle
^{
{1}/{2^{q}}-{1}/({n(n+1)})}+1\bigr)\nonumber\\
&&\eqntext{(\mbox{estimate }q) }\\
&\vdots& \nonumber\\
\langle\phi_T(0)^{2\cdot2^0} \rangle&{\lesssim} &\underbrace
{\langle\phi_T(0) \rangle^2}_{\stackrel{\mathrm{Lemma}\mbox{
\fontsize
{8.36pt}{10.36pt}\selectfont{\ref{lem:app-corr}}}}{=}0}{}+{}\mu
_d(T)\bigl(\langle\phi_T(0)^{2n} \rangle^{{1}/{n}-
{1}/({n(n+1)})}+1\bigr)\nonumber\\
&&\eqntext{(\mbox{estimate }l).}
\end{eqnarray}
We then take the power $2^q$ of each (estimate $q$) and obtain using
Young's inequality:
%
%
\begin{eqnarray}\label{eq:referee-2.1}\quad
\langle\phi_T(0)^{2n} \rangle&{\lesssim} &\langle\phi_T(0)^{n}
\rangle^2+\mu
_d(T)\bigl(\langle\phi_T(0)^{2n} \rangle^{1-{1}/({n(n+1)})}+1\bigr), \nonumber\\
&\vdots& \nonumber\\
\langle\phi_T(0)^{2\cdot2^{l-q}} \rangle^{2^q}&{\lesssim}& \langle
\phi_T(0)^{2^{l-q}} \rangle^{2^{q+1}}\nonumber\\
&&{}+\mu_d(T)^{2^q}\bigl(\langle\phi
_T(0)^{2n} \rangle^{1-{2^q}/({n(n+1)})}+1\bigr),\nonumber\\[-8pt]\\[-8pt]
\langle\phi_T(0)^{2^{l-q}} \rangle^{2^{q+1}}&{\lesssim}& \langle
\phi_T(0)^{2^{l-q-1}} \rangle^{2^{q+2}}\nonumber\\
&&{}+\mu_d(T)^{2^{q+1}}\bigl(\langle
\phi_T(0)^{2n} \rangle^{1-{2^{q+1}}/({n(n+1)})}+1\bigr),\nonumber\\
&\vdots& \nonumber\\
\langle\phi_T(0)^{2} \rangle^{n} &{\lesssim} &\mu_d(T)^{n}\bigl(\langle
\phi_T(0)^{2n} \rangle^{1-{1}/({n+1})}+1\bigr).\nonumber
\end{eqnarray}
Since the multiplicative constants in each line of (\ref
{eq:referee-2.1}) only depend on $\alpha,\beta,d,n$ and $q$, a
linear combination of these $l+1$ inequalities with suitable positive
coefficients allows us to cancel the respective terms both on
the left- and right-hand sides, which yields
%
%
\begin{equation}\label{eq:var-step1}
\langle\phi_T(0)^{2n} \rangle\lesssim\sum_{q=0}^{l}\mu
_d(T)^{2^q}\bigl(\langle\phi_T(0)^{2n} \rangle^{1-{2^q}/({n(n+1)})}+1\bigr).
\end{equation}
Using Young's inequality, each term gives the same contribution and
(\ref{eq:var-step1}) turns into
%
%
\begin{equation}
\langle\phi_T(0)^{2n} \rangle\lesssim\mu_d(T)^{n(n+1)}. \label
{eq:prop1-delta}
\end{equation}
Formula (\ref{eq:estim-prop}) is then proved for all $q\leq2n$ using
H\"older's inequality in probability.

\textit{Step} 2. Estimate for the Green's function.

Let $p>2$ be as in Lemma \ref{lem:int-grad}.
We shall prove that for all $q\geq1$ and $R\gg1$ the following holds
%
%
\begin{eqnarray}
\label{eq:int-grad-d=2+}
&&\mbox{for }d=2\qquad \int_{R< |z|\leq2R}|\nabla_zG_T(z,0)
|^q\,dz\nonumber\\[-8pt]\\[-8pt]
&&\phantom{\mbox{for }d=2\qquad}\qquad\lesssim
R^{2\max\{1,{q/ p}\}}R^{-q}\min\bigl\{1,\sqrt{T}R^{-1}\bigr\}^{q},\nonumber
\\
\label{eq:int-grad+}
&&\mbox{for }d>2\qquad \int_{R< |z|\leq2R}|\nabla_zG_T(z,0)
|^q\,dz\nonumber\\[-8pt]\\[-8pt]
&&\phantom{\mbox{for }d>2\qquad}\qquad\lesssim
R^{d\max\{1,{q/ p}\}}(R^{1-d})^q.\nonumber
\end{eqnarray}
We split the argument into two parts to treat $q\geq p$ and $q< p$,
respectively.
For $q\geq p$, we use the discrete $L^p-L^{ q}$ estimate:
\[
\biggl(\int_{R< |z|\leq2R}|\nabla_zG_T(z,0) |^q\,dz\biggr)^{1/q} \leq
\biggl(\int_{R< |z|\leq2R}|\nabla_zG_T(z,0) |^{ p}\,dz\biggr)^{1/{ p}}.
\]
Combined with (\ref{eq:int-grad}) in Lemma \ref{lem:int-grad}, it
proves (\ref{eq:int-grad-d=2+}) and (\ref{eq:int-grad+}).

For $q< p$, we simply use H\"older's inequality with
exponents $(\frac{p}{q},\frac{p}{p-q})$ in the form
\[
\biggl(R^{-d}\int_{R< |z|\leq2R}|\nabla_zG_T(z,0) |^q\,dz\biggr)^{1/q}
\lesssim\biggl(R^{-d}\int_{R< |z|\leq2R}|\nabla_zG_T(z,0) |^{
p}\,dz\biggr)^{1/{ p}},
\]
that we also combine with (\ref{eq:int-grad}).

\textit{Step} 3. General estimate.

Let $\chi\geq0$ be a random variable.
In order to prove (\ref{eq:key-prop-d=2}), we will need to estimate
terms of the form
\[
\int_{\Z^d}\langle\chi|\nabla_z G_T(z,0)|^q \rangle^{1/r} \,dz
\]
for $q,r>1$. Relying on (\ref{eq:int-grad-d=2+}) and (\ref
{eq:int-grad+}), we show that
%
%
\begin{eqnarray}\label{eq:gen-estim}\quad
&&\int_{\Z^d}\langle\chi|\nabla_z G_T(z,0)|^q \rangle^{1/r}
\,dz\nonumber\\[-8pt]\\[-8pt]
&&\qquad\lesssim \langle\chi\rangle^{1/r}\cases{
1, &\quad if $\displaystyle  d\max\biggl\{1,1-\frac{1}{r}+\frac{q}{rp}\biggr\}
+(1-d)\frac
{q}{r}<0$,\cr
&\quad $d\geq2$,\cr
\ln T, &\quad if $\displaystyle 2\max\biggl\{1,1-\frac{1}{r}+\frac{q}{rp}\biggr\}-\frac
{q}{r}=0$,\cr
&\quad $d=2$.}\nonumber
\end{eqnarray}
Note that\vspace*{1pt} there is no overlap in (\ref{eq:gen-estim}). For $d>2$, we
will only make use of the estimate
with $d\max\{1,1-\frac{1}{r}+\frac{q}{rp}\}+(1-d)\frac{q}{r}<0$. For
$d=2$, we will use the estimate both
with $2\max\{1,1-\frac{1}{r}+\frac{q}{rp}\}-\frac{q}{r}<0$, and with
$2\max\{1,1-\frac{1}{r}+\frac{q}{rp}\}-\frac{q}{r}=0$,
which requires a specific argument.

Let $i_{\min}\in\N, i_{\min}\sim1$ be such that
Lemma \ref{lem:int-grad} holds for $R\geq2^{i_{\min}}$.
To prove (\ref{eq:gen-estim}), we use a dyadic decomposition of $\Z^d$
in annuli of radii $R_i=2^i$:
%
%
\begin{eqnarray}\label{eq:ajout-ref}
&&\int_{\Z^d}\langle\chi|\nabla_z G_T(z,0)|^q \rangle^{1/r}
\,dz
\nonumber\\
&&\qquad= \int_{|z|\leq2^{i_{\min}}}\langle\chi|\nabla_z
G_T(z,0)|^q \rangle^{1/r} \,dz\\
&&\qquad\quad{} +
\sum_{i= i_{\min}}^\infty\int_{R_i< |z|\leq R_{i+1}}\langle
\chi|\nabla_z G_T(z,0)|^q \rangle^{1/r} \,dz.\nonumber
\end{eqnarray}
Using Corollary \ref{coro:unif-bound-grad}, we bound the first term of
the right-hand side by
\[
\int_{|z|\leq2^{i_{\min}}}\langle\chi|\nabla_z
G_T(z,0)|^q \rangle^{1/r} \,dz \lesssim\langle\chi\rangle^{1/r}.
\]
For the second term of the right-hand side, we appeal to H\"older's inequality
with $(r,\frac{r}{r-1})$:
\begin{eqnarray*}
&&\sum_{i= i_{\min}}^\infty\int_{R_i< |z|\leq
R_{i+1}}\langle\chi|\nabla_z G_T(z,0)|^q \rangle^{1/r} \,dz\\
&&\qquad \lesssim \sum_{i= i_{\min}}^\infty(R_i^d)^{1-1/r}
\biggl(\int
_{R_i < |z|\leq R_{i+1}}\langle\chi|\nabla_z G_T(z,0)|^q \rangle \,dz\biggr)^{1/r},
\end{eqnarray*}
so that (\ref{eq:ajout-ref}) turns into
\begin{eqnarray*}
&&\int_{\Z^d}\langle\chi|\nabla_z G_T(z,0)|^q \rangle^{1/r}
\,dz
\\
&&\qquad= \langle\chi\rangle^{1/r}+\sum_{i= i_{\min}}^\infty
(R_i^d)^{1-1/r}\biggl\langle\chi\int_{R_i< |z|\leq R_{i+1}}|\nabla_z
G_T(z,0)|^q \,dz \biggr\rangle^{1/r}.
\end{eqnarray*}
Using then (\ref{eq:int-grad-d=2+}) and (\ref{eq:int-grad+}), we get
\begin{eqnarray*}
&&\biggl\langle\chi\int_{R_i< |z|\leq R_{i+1}}|\nabla_z G_T(z,0)|^q \,dz
\biggr\rangle
\\
&&\qquad\lesssim\cases{
\langle\chi\rangle R_i^{2\max\{1,{q/p}\}}R_i^{-q}\min\bigl\{
1,\sqrt
{T}R_i^{-1}\bigr\}^{q}, &\quad $d=2$, \cr
\langle\chi\rangle R_i^{d\max\{1,{q/p}\}}(R_i^{1-d})^q, &\quad
$d>2$.}
\end{eqnarray*}
Hence,
\begin{eqnarray*}
&&\int_{\Z^d}\langle\chi|\nabla_z G_T(z,0)|^q \rangle^{1/r} \,dz\\
&&\qquad\lesssim \cases{
\displaystyle \langle\chi\rangle^{1/r} \sum_{i=0}^\infty R_i^{2\max\{1,1-
{1/r}+{q}/({rp})\}-{q/r}}\min\bigl\{1,\sqrt{T}R_i^{-1}\bigr\}^{q/r},
&\quad
$d= 2$, \cr
\displaystyle \langle\chi\rangle^{1/r} \sum_{i=0}^\infty R_i^{d\max\{1,1-
{1/r}+{q}/({rp})\}+(1-d){q/r}},&\quad $d> 2$.}
\end{eqnarray*}
We distinguish two cases.
If $d\max\{1,1-\frac{1}{r}+\frac{q}{rp}\}+(1-d)\frac{q}{r}<0$, then
\[
\int_{\Z^d}\langle\chi|\nabla_z G_T(z,0)|^q \rangle^{1/r} \,dz
\lesssim
\langle\chi\rangle^{1/r} \sum_{i=0}^\infty R_i^{d\max\{1,1-
{1/r}+{q}/({rp})\}+(1-d){q/r}} \lesssim\langle\chi\rangle^{1/r}.
\]
This proves the first estimate of (\ref{eq:gen-estim}).
For $d=2$, and $2\max\{1,1-\frac{1}{r}+\frac{q}{rp}\}-\frac
{q}{r}=0$, then
we obtain
\begin{eqnarray*}
\int_{\Z^2}\langle\chi|\nabla_z G_T(z,0)|^q \rangle^{1/r} \,dz
&\lesssim&
\langle\chi\rangle^{1/r} \sum_{i=0}^\infty\min\bigl\{1,\sqrt
{T}R_i^{-1}\bigr\}
^{q/r} \\
&\lesssim& \langle\chi\rangle^{1/r}\Biggl(\ln T+\sum_{i=0}^\infty
R_i^{-q/r} \Biggr)
\\
&\lesssim& \langle\chi\rangle^{1/r}(1+\ln T).
\end{eqnarray*}
This proves the second estimate of (\ref{eq:gen-estim}).

\textit{Step} 4. Proof of (\ref{eq:key-prop-d=2}).

Let $n\geq1$ and $n\geq m \geq1$. We first treat the first
term of the left-hand side of (\ref{eq:key-prop-d=2}).
In that case H\"older's inequality in probability with exponents $(n+1,\frac
{n+1}{n})$ and the stationarity of $\nabla\phi_T$ show
%
%
\begin{eqnarray}\label{eq:prop-step2-1}\quad
&&\int_{\Z^d} \bigl\langle\phi_T(0)^{2(m-1)}\bigl(|\nabla\phi
_T(z)|+1\bigr)^2|\nabla_{z}G_T(z,0)|^{2} \bigr\rangle \,dz \nonumber\\
&&\qquad\lesssim\int_{\Z^d}\bigl(\bigl\langle|\nabla\phi_T(z)|^{2(n+1)} \bigr\rangle
^{{1}/({n+1})}+1\bigr)\nonumber\\
&&\qquad\quad\hspace*{13.4pt}{}\times\bigl\langle|\phi_T(0)|^{{2(m-1)(n+1)}/{n}}|\nabla
_{z}G_T(z,0)|^{{2(n+1)}/{n}} \bigr\rangle^{{n}/({n+1})} \,dz\\
&&\qquad= \bigl(\bigl\langle|\nabla\phi_T(0)|^{2(n+1)} \bigr\rangle^{
{1}/({n+1})}+1
\bigr)\nonumber\\
&&\qquad\quad{}\times\int_{\Z^d}\bigl\langle|\phi_T(0)|^{{2(m-1)(n+1)}/{n}}|\nabla
_{z}G_T(z,0)|^{{2(n+1)}/{n}} \bigr\rangle^{{n}/({n+1})}\,dz.\nonumber
\end{eqnarray}
We apply Lemma \ref{lem:Lp-grad} to bound the first ensemble average in
(\ref{eq:prop-step2-1}):
%
%
\begin{eqnarray}\label{eq:prop-step2-3}
&&\bigl\langle|\nabla\phi_T(0)|^{2(n+1)} \bigr\rangle\nonumber\\
&&\qquad\lesssim \Biggl\langle\sum
_{i=1}^d|\nabla\phi_T(0)|^{2}\bigl(\phi_T(0)^{2n}+\phi_T(\ee_i)^{2n}\bigr)
\Biggr\rangle
\nonumber\\[-8pt]\\[-8pt]
&&\hspace*{-14.51pt}\qquad\stackrel{\mathrm{stationarity}}{=} 2 \Biggl\langle\sum
_{i=1}^d|\nabla\phi_T(0)|^{2}\phi_T(0)^{2n} \Biggr\rangle\nonumber\\
&&\qquad\hspace*{-5.91pt} \stackrel{\mbox{\fontsize{8.36pt}{10.36pt}\selectfont{(\ref
{eq:phi-grad-phi})}}}{\lesssim}\langle\phi_T(0)^{2n} \rangle.\nonumber
\end{eqnarray}
We now want to apply Step 3 to the right-hand side integral of (\ref
{eq:prop-step2-1}), that is, setting $q=\frac{2(n+1)}{n}$, $r=\frac
{n+1}{n}$ and $\chi=|\phi_T(0)|^{{2(m-1)(n+1)}/{n}}$.
Estimate (\ref{eq:gen-estim}) involves the number
%
%
\begin{eqnarray}\label{eq:prop-step2-2}
&&d\max\biggl\{1,1-\frac{1}{r}+\frac{q}{rp}\biggr\}+(1-d)\frac{q}{r}\nonumber\\[-8pt]\\[-8pt]
&&\qquad=d \max\biggl\{
1,\frac
{1}{n+1}+\frac{2}{ p}\biggr\}+2(1-d).\nonumber
\end{eqnarray}
We distinguish the cases $d>2$ and $d=2$. For $d>2$, we have that the
number (\ref{eq:prop-step2-2}) is equal to $d+2(1-d)=2-d$ and thus
negative for $n$ sufficiently large since $p>2$.
Hence, (\ref{eq:gen-estim}) yields
\begin{eqnarray*}
&&{\int_{\Z^d}\bigl\langle|\phi_T(0)|^{{2(m-1)(n+1)}/{n}}|\nabla
_{z}G_T(z,0)|^{{2(n+1)}/{n}} \bigr\rangle^{{n}/({n+1})} \,dz}\\
&&\qquad\lesssim \bigl\langle|\phi_T(0)|^{{2(m-1)(n+1)}/{n}} \bigr\rangle
^{{n}/({n+1})}
\leq \langle|\phi_T(0)|^{2n} \rangle^{({m-1})/{n}},
\end{eqnarray*}
where in the last inequality we appealed to Jensen in probability using
\[
\frac{2(m-1)(n+1)}{n}\leq\frac{2(n-1)(n+1)}{n}\leq2n.
\]
The combination of this with (\ref{eq:prop-step2-1}) and (\ref
{eq:prop-step2-3}) yields as desired
\begin{eqnarray*}
&&\int_{\Z^d} \bigl\langle\phi_T(0)^{2(m-1)}\bigl(|\nabla\phi
_T(z)|+1\bigr)^2|\nabla_{z}G_T(z,0)|^{2} \bigr\rangle \,dz\\
&&\qquad\lesssim\langle\phi(0)^{2n} \rangle^{{1}/({n+1})+
({m-1})/{n}}+1
=\langle\phi(0)^{2n} \rangle^{{m/n}-{1}/({n(n+1)})}+1.
\end{eqnarray*}

We turn to the case $d=2$. We note that the number (\ref
{eq:prop-step2-2}) is zero for $n$ large enough since $p>2$. Thus, from
(\ref{eq:gen-estim}), we infer as we did above that
\begin{eqnarray*}
&&\int_{\Z^d} \bigl\langle\phi_T(0)^{2(m-1)}\bigl(|\nabla\phi
_T(z)|+1\bigr)^2|\nabla_{z}G_T(z,0)|^{2} \bigr\rangle \,dz\\
&&\qquad\lesssim
(\ln T)\bigl(\langle\phi(0)^{2n} \rangle^{{m/n}-
{1}/({n(n+1)})}+1\bigr).
\end{eqnarray*}
Let us now treat the second term of the left-hand side of (\ref
{eq:key-prop-d=2}), which differs from the first term only when $m\geq2$.
As for the first term, H\"older's inequality in probability with
$(\frac
{n+1}{m},\frac{n+1}{n-m+1})$, the stationarity of $\nabla\phi_T$ and
Lemma \ref{lem:Lp-grad} imply
%
%
\begin{eqnarray}\label{eq:prop-step2-5}
&&\int_{\Z^d} \bigl\langle\bigl(|\nabla\phi_T(z)|+1\bigr)^{2m}|\nabla
_{z_i}G_T(z,0)|^{2m} \bigr\rangle \,dz\nonumber\\
&&\qquad\lesssim\int_{\Z^d} \bigl(\bigl\langle|\nabla\phi_T(z)|^{2(n+1)} \bigr\rangle
^{
{m}/({n+1})}+1\bigr)\nonumber\\
&&\qquad\quad\hspace*{13.3pt}{}\times\bigl\langle|\nabla_{z}G_T(z,0)|^{{2(n+1)m}/({n-m+1})}
\bigr\rangle^{({n-m+1})/({n+1})} \,dz \\
&&\qquad\lesssim \bigl(\langle\phi_T(0)^{2n} \rangle^{{m}/({n+1})}+1\bigr)\nonumber\\
&&\qquad\quad{}\times\int
_{\Z^d} \bigl\langle|\nabla_{z}G_T(z,0)|^{{2(n+1)m}/({n-m+1})} \bigr\rangle
^{({n-m+1})/({n+1})} \,dz.\nonumber
\end{eqnarray}
We use (\ref{eq:gen-estim}) with $\chi\equiv1$, $q=\frac
{2(n+1)m}{n-m+1}$ and $r=\frac{n+1}{n-m+1}$, in which case we have
%
%
\begin{eqnarray}\label{eq:prop-step2-4}
&&
d\max\biggl\{1,1-\frac{1}{r}+\frac{q}{rp}\biggr\}
+(1-d)\frac{q}{r}\nonumber\\[-8pt]\\[-8pt]
&&\qquad=d \max\biggl\{1,\frac
{m}{n+1}+\frac{2m}{ p}\biggr\}+(1-d)2m.\nonumber
\end{eqnarray}
We claim that this number is negative for $n$ sufficiently large.
Indeed, if $\max\{1$, $\frac{m}{n+1}+\frac{2m}{ p}\}=1$, then
\begin{eqnarray*}
d \max\biggl\{1,\frac{m}{n+1}+\frac{2m}{ p}\biggr\}
+(1-d)2m&=&d+2m(1-d)\\
&=&(2m-1)(1-d)+1 <0
\end{eqnarray*}
since $d\geq2$ and $m\geq2$.
Otherwise, $\max\{1,\frac{m}{n+1}+\frac{2m}{ p}\}=\frac
{m}{n+1}+\frac
{2m}{ p}$, and
\begin{eqnarray*}
d \max\biggl\{1,\frac{m}{n+1}+\frac{2m}{ p}\biggr\}+(1-d)2m&=& 2m\biggl(d
\biggl(\frac
{1}{2(n+1)}+\frac{1}{ p}\biggr) +1-d\biggr) \\&<& 2m\biggl(1-\frac
{d}{2}\biggr) \leq0
\end{eqnarray*}
for $d\geq2$ and $n$ large enough since $\frac{1}{ p}<\frac{1}{2}$.
This shows that (\ref{eq:prop-step2-4}) is negative so that we obtain
by (\ref{eq:gen-estim})
\[
\int_{\Z^d} \bigl\langle|\nabla_{z}G_T(z,0)|^{{2(n+1)m}/({n-m+1})}
\bigr\rangle^{
({n-m+1})/({n+1})} \,dz \lesssim1.
\]
Combining this with (\ref{eq:prop-step2-5}) yields
\begin{eqnarray*}
\int_{\Z^2} \bigl\langle\bigl(|\nabla\phi_T(z)|+1\bigr)^{2m}|\nabla
_{z}G_T(z,0)|^{2m} \bigr\rangle \,dz &\lesssim& \langle\phi_T(0)^{2n}
\rangle^{{m}/({n+1})}+1\\
&=& \langle\phi_T(0)^{2n} \rangle^{{m/n}-{m}/({n(n+1)})}+1
\\
&\leq& \langle\phi_T(0)^{2n} \rangle^{{m/n}-{1}/({n(n+1)})}+1.
\end{eqnarray*}
This concludes the proof of the proposition.

\subsection{\texorpdfstring{Proof of Theorem \protect\ref{th:main}}{Proof of Theorem 2.1.}}

Let us define the spatial average of a function $h\dvtx\Z^d\to
\R
$ with the mask $\eta_L$ by
\[
\llangle h \rrangle_L:=\int_{\Z^d}h(x)\eta_L(x)\,dx,
\]
where $\eta_L$ satisfies
%
%
\begin{eqnarray}\label{eq:eta}
\eta_L\dvtx\Z^d \to[0,1]\qquad \operatorname{supp}(\eta_L)
&\subset&(-L,L)^d,\nonumber\\[-8pt]\\[-8pt]
\int_{\Z^d} \eta_L(x) \,dx&=&1,\qquad |\nabla\eta_L|\lesssim
L^{-d-1}.\nonumber
\end{eqnarray}
The claim of the theorem is that there exists $q$ depending only on
$\alpha,\beta$, and $d$ such that
\[
\operatorname{var}[\llangle T^{-1}\phi_T^2+(\nabla\phi_T+\xi)\cdot
A(\nabla\phi_T+\xi) \rrangle_L]\lesssim L^{-d}\mu_d(T)^q,
\]
where $\mu_d(T)=1$ for $d>2$ and $\mu_d(T)=\ln T$ for $d=2$.
Since we are not interested in the precise value of $q$, we adopt the
convention that $q$ is a nonnegative exponent which only depends on
$\alpha,\beta$, and $d$ but which may vary from line to line in the proof.

Starting point is the estimate provided by Lemmas \ref
{lem:var-estim} and \ref{lem:depend-coeff}
%
%
\begin{eqnarray}\label{eq:main-var-start}\quad
&&\operatorname{var}[\llangle T^{-1}\phi_T^2+(\nabla\phi_T+\xi
)\cdot A(\nabla\phi_T+\xi) \rrangle_L]\nonumber\\[-8pt]\\[-8pt]
&&\qquad\lesssim\biggl\langle\sum_{e}\sup_{a(e)}\biggl|\frac{\partial}{\partial
a(e)}\llangle T^{-1}\phi_T^2+(\nabla\phi_T+\xi)\cdot A(\nabla\phi
_T+\xi) \rrangle_L\biggr|^2 \biggr\rangle.\nonumber
\end{eqnarray}

\textit{Step} 1. In this step, using the notation $e=[z,z+\ee_i]$, we
establish
the formula
%
%
\begin{eqnarray}\label{eq:theo-step1-1}
&&\frac{\partial}{\partial a(e)}\llangle T^{-1}\phi_T^2+(\nabla\phi
_T+\xi)\cdot A(\nabla\phi_T+\xi) \rrangle_L\nonumber\\
&&\qquad=2\int_{\Z^d} \bigl(\nabla_i
\phi_T(z)+\xi_i\bigr)\nabla_{z_i}G_T(z,x)\nonumber\\[-8pt]\\[-8pt]
&&\qquad\quad\hspace*{21pt}{}\times
\Biggl(\sum_{j=1}^da(x-\ee_j,x)\nabla_j^*\eta_L(x)\bigl(\nabla^*_j \phi
_T(x)+\xi
_j\bigr)\Biggr) \,dx \nonumber\\
&&\qquad\quad{} +\eta_L(z)(\nabla_i \phi_T+\xi_i)^2(z).\nonumber
\end{eqnarray}
Indeed, by definition of $\llangle\cdot\rrangle_L$ we have
\begin{eqnarray*}
&&\frac{\partial}{\partial a(e)}\llangle T^{-1}\phi_T^2+(\nabla\phi
_T+\xi)\cdot A(\nabla\phi_T+\xi) \rrangle_L\\
&&\qquad=\int_{\Z^d}\eta_L(x)\,\frac{\partial}{\partial a(e)}\bigl(
T^{-1}\phi
_T^2+(\nabla\phi_T+\xi)\cdot A(\nabla\phi_T+\xi)\bigr)(x)\,dx.
\end{eqnarray*}
We note
\begin{eqnarray*}
&&\frac{\partial}{\partial a(e)}\bigl(T^{-1}\phi_T^2+
(\nabla
\phi_T+\xi)\cdot A(\nabla\phi_T+\xi)\bigr)(x)\\
&&\qquad=\biggl(2T^{-1}\phi_T\,\frac{\partial\phi_T}{\partial a(e)} +2
\nabla\,
\frac{\partial\phi_T}{\partial a(e)} \cdot A(\nabla\phi_T+\xi)\\
&&\qquad\quad\hspace*{49pt}{}+
(\nabla\phi_T+\xi)\cdot\frac{\partial A}{\partial a(e)}(\nabla
\phi
_T+\xi) \biggr)(x) \\
&&\qquad=2T^{-1}\biggl(\phi_T\,\frac{\partial\phi_T}{\partial a(e)}\biggr)(x)+
2\biggl(\nabla\,\frac{\partial\phi_T}{\partial a(e)} \cdot A(\nabla
\phi
_T+\xi)\biggr)(x)\\
&&\qquad\quad{}+(\nabla_i \phi_T+\xi_i)^2(z)\delta(x-z),
\end{eqnarray*}
so that
%
%
\begin{eqnarray}\label{eq:main-diff}\quad
&&\frac{\partial}{\partial a(e)}\llangle T^{-1}\phi_T^2+(\nabla\phi
_T+\xi)\cdot A(\nabla\phi_T+\xi) \rrangle_L\nonumber\\
&&\qquad=2\int_{\Z^d}\biggl(\eta_L\biggl(T^{-1}\phi_T\,\frac{\partial\phi
_T}{\partial a(e)}+ \nabla\,\frac{\partial\phi_T}{\partial a(e)}
\cdot
A(\nabla\phi_T+\xi)\biggr)\biggr)(x)\,dx \\
&&\qquad\quad{} +\eta_L(z)(\nabla_i \phi_T+\xi_i)^2(z).\nonumber
\end{eqnarray}
Using the discrete integration by parts formula of Definition \ref
{def:ibp}, the first term of the right-hand side of (\ref
{eq:main-diff}) turns into
%
%
\begin{eqnarray}\label{eq:theo-step1-2}
&&\int_{\Z^d}\biggl(\eta_L \biggl(T^{-1}\phi_T\,\frac{\partial
\phi
_T}{\partial a(e)}+ \nabla\,\frac{\partial\phi_T}{\partial a(e)}
\cdot
A(\nabla\phi_T+\xi) \biggr)\biggr)(x)\,dx\nonumber\\
&&\qquad= - \int_{\Z^d} \frac{\partial\phi_T}{\partial a(e)}(x)\nabla
^*\cdot
\bigl(\eta_LA(\nabla\phi_T+\xi) \bigr)(x) \,dx\\
&&\qquad\quad{} + \int_{\Z
^d}\biggl(\eta
_L T^{-1}\phi_T\,\frac{\partial\phi_T}{\partial a(e)}\biggr)(x)
\,dx.\nonumber
\end{eqnarray}
We now use the following discrete Leibniz rule:
\begin{eqnarray*}
\nabla^*\cdot\bigl(\eta_LA(\nabla\phi_T+\xi)
\bigr)(x)
&=& \eta_L(x)\bigl(\nabla^*\cdot A(\nabla\phi_T+\xi)\bigr)(x)\\
&&{} +
\sum
_{j=1}^d \nabla^*_j\eta_L(x) [A(\nabla\phi_T+\xi)]_j(x-\ee_j),
\end{eqnarray*}
where $[A(\nabla\phi_T+\xi)]_j$ denotes the $j$th coordinate of the
vector $A(\nabla\phi_T+\xi)$.
For notational convenience, we take advantage of the diagonal structure
of $A$ (although this is not crucial) to rewrite
the latter equality in the form
%
%
\begin{eqnarray}\label{eq:theo-step1-3}
&&\nabla^*\cdot\bigl(\eta_LA(\nabla\phi_T+\xi)
\bigr)(x)\nonumber\\
&&\qquad= \eta_L(x)\bigl(\nabla^*\cdot A(\nabla\phi_T+\xi)\bigr)(x) \\
&&\qquad\quad{} + \sum
_{j=1}^d a(x-\ee_j,x)\nabla^*_j\eta_L(x) \bigl(\nabla^*_j \phi_T(x)+\xi
_j\bigr).\nonumber
\end{eqnarray}
The combination of (\ref{eq:theo-step1-3}) with (\ref{eq:theo-step1-2})
and the use of the equation satisfied by $\phi_T$,
\[
T^{-1}\phi_T-\nabla^*\cdot A(\nabla\phi_T+\xi)=0,
\]
yield
\begin{eqnarray*}
&&\int_{\Z^d}\biggl(\eta_L \biggl(T^{-1}\phi_T\,\frac{\partial
\phi
_T}{\partial a(e)}+ \nabla\,\frac{\partial\phi_T}{\partial a(e)}
\cdot
A(\nabla\phi_T+\xi) \biggr)\biggr)(x)\,dx\\
&&\qquad= - \int_{\Z^d} \frac{\partial\phi_T}{\partial a(e)}(x)\Biggl(
\sum
_{j=1}^d a(x-\ee_j,x)\nabla^*_j\eta_L(x) \bigl(\nabla^*_j \phi_T(x)+\xi_j\bigr)
\Biggr) \,dx.
\end{eqnarray*}
Using now Lemma \ref{lem:diff-phi}, this turns into
%
%
\begin{eqnarray}\label{eq:main-diff-2}
&&\int_{\Z^d}\biggl(\eta_L \biggl(T^{-1}\phi_T\,\frac{\partial
\phi
_T}{\partial a(e)}+ \nabla\,\frac{\partial\phi_T}{\partial a(e)}
\cdot
A(\nabla\phi_T+\xi) \biggr)\biggr)(x)\,dx\nonumber\\
&&\qquad\stackrel{\mbox{\fontsize{8.36pt}{10.36pt}\selectfont{(\ref{eq:diff-phi-1})}}}{=}
\int_{\Z^d} \bigl(\nabla_i \phi_T(z)+\xi _i\bigr)\nabla
_{z_i}G_T(z,x) \\
&&\qquad\quad\hspace*{24.3pt}{}\times \Biggl( \sum_{j=1}^d a(x-\ee_j,x)\nabla^*_j\eta_L(x)
\bigl(\nabla^*_j \phi_T(x)+\xi_j\bigr) \Biggr) \,dx.\nonumber
\end{eqnarray}
Inserting (\ref{eq:main-diff-2}) into (\ref{eq:main-diff}) proves
(\ref
{eq:theo-step1-1}).

\textit{Step} 2. In this step, we provide the estimate
%
%
\begin{eqnarray}\label{eq:theo-step2-1}\qquad
&&\sup_{a(e)}\biggl|\frac{\partial}{\partial a(e)}\llangle T^{-1}\phi
_T^2+(\nabla\phi_T+\xi)\cdot A(\nabla\phi_T+\xi) \rrangle_L
\biggr|\nonumber\\
&&\qquad\lesssim \int_{\Z^d}|\nabla_zG_T(z,x)| |\nabla^* \eta
_L(x)|\bigl(|\nabla
^* \phi_T(x)|^2+|\nabla\phi_T(z)|^2+1\bigr) \,dx\\
&&\qquad\quad{} +\eta_L(z)\bigl(|\nabla\phi_T(z)|^2+1\bigr).
\nonumber
\end{eqnarray}
Indeed, from Step 1, the boundedness of $a$, and $|\xi|=1$, we infer that
%
%
\begin{eqnarray}\label{eq:theo-step2-2}
&&\sup_{a(e)}\biggl|\frac{\partial}{\partial a(e)}\llangle T^{-1}\phi
_T^2+(\nabla\phi_T+\xi)\cdot A(\nabla\phi_T+\xi) \rrangle_L
\biggr|\nonumber\\
&&\qquad\lesssim \int_{\Z^d} \Bigl(\sup_{a(e)}|\nabla_i \phi_T(z)|+1\Bigr)\sup
_{a(e)}|\nabla_{z_i}G_T(z,x)|| \nabla^*\eta_L(x)|\nonumber\\[-8pt]\\[-8pt]
&&\qquad\quad\hspace*{13.4pt}{}\times\Bigl(\sup
_{a(e)}|\nabla^*
\phi_T(x)|+1\Bigr) \,dx\nonumber\\
&&\qquad\quad{} +\eta_L(z)\Bigl(\sup_{a(e)}|\nabla_i
\phi_T(z)|^2+1\Bigr).\nonumber
\end{eqnarray}
Hence, in the remainder of this step, we have to deal with the suprema
over $a(e)$.
Recalling that $e=[z,z+\ee_i]$, the two following inequalities are
consequences of Lemmas \ref{lem:diff-Green} and \ref{lem:diff-phi}:
\begin{eqnarray*}
\sup_{a(e)}|\nabla_{z_i} G_T(z,x)|&\stackrel{\mbox{\fontsize
{8.36pt}{10.36pt}\selectfont{(\ref{eq:bd-G(x,e)})}}}{\lesssim}&
|\nabla
_{z_i} G_T(z,x)|\qquad \mbox{for all }x\in\Z^d,\\
\sup_{a(e)}|\nabla_i \phi_T(z)|&\stackrel{\mbox{\fontsize
{8.36pt}{10.36pt}\selectfont{(\ref{eq:diff-phi-3})}}}{\lesssim}&
|\nabla
_i \phi_T(z)|+1.
\end{eqnarray*}
The last inequality we need is
\[
\sup_{a(e)}|\nabla^* \phi_T(x)|\lesssim|\nabla^* \phi_T(x)|+\sup
_{a(e)}|\nabla_i \phi_T(z)|+1 \stackrel{\mbox{\fontsize
{8.36pt}{10.36pt}\selectfont{(\ref{eq:diff-phi-3})}}}{\lesssim}
|\nabla
^* \phi_T(x)|+|\nabla_i \phi_T(z)|+1.
\]
It is then proved combining the boundedness of $a$ and the following
bound on the derivative of $\nabla^* \phi_T(x)$ with respect to $a(e)$:
\begin{eqnarray*}
\biggl|\frac{\partial}{\partial a(e)}\,\nabla^* \phi_T(x)\biggr|&=&
\biggl|\nabla_x^*\, \frac{\partial}{\partial a(e)}\,\phi_T(x)\biggr|\\
&\stackrel{\mbox{\fontsize{8.36pt}{10.36pt}\selectfont{(\ref
{eq:diff-phi-1})}}}{=}& \bigl|\nabla_x^* \bigl(\bigl(\nabla_i \phi
_T(z)+\xi
_i\bigr)\nabla_{z_i}G_T(z,x)\bigr) \bigr| \\
&=&\bigl|\bigl(\nabla_i \phi_T(z)+\xi_i\bigr)\nabla_{z_i}\nabla_x^* G_T(z,x)
\bigr| \\
&\leq& 2\bigl(|\nabla_i \phi_T(z)|+|\xi_i|\bigr)\sup_{\Z^d\times\Z
^d}|\nabla
G_T|\\
&{\lesssim}& |\nabla_i \phi_T(z)|+1,
\end{eqnarray*}
where we have used the uniform bound on $\nabla G_T$ provided by
Corollary \ref{coro:unif-bound-grad}.

Combining these three inequalities with (\ref
{eq:theo-step2-2}) yields
\begin{eqnarray*}
&&\sup_{a(e)}\biggl|\frac{\partial}{\partial a(e)}\llangle T^{-1}\phi
_T^2+(\nabla\phi_T+\xi)\cdot A(\nabla\phi_T+\xi) \rrangle_L\biggr|\\
&&\qquad\lesssim \int_{\Z^d}\bigl(|\nabla\phi_T(z)|+1\bigr)|\nabla_zG_T(z,x)|
|\nabla
^* \eta_L(x)|\\
&&\qquad\quad\hspace*{12.7pt}{}\times\bigl(|\nabla^* \phi_T(x)|+|\nabla\phi_T(z)|+1\bigr) \,dx\\
&&\qquad\quad{} +\eta_L(z)\bigl(|\nabla\phi_T(z)|^2+1\bigr)
\end{eqnarray*}
from which we deduce (\ref{eq:theo-step2-1}).

\textit{Step} 3. In this step, we argue that
%
%
\begin{eqnarray}
&&\operatorname{var}[\llangle T^{-1}\phi_T^2+(\nabla\phi_T+\xi
)\cdot A(\nabla\phi_T+\xi) \rrangle_L]\nonumber\\
\label{eq:main:term11}
&&\qquad\lesssim \biggl\langle\int_{\Z^d}\biggl(\int_{\Z^d}|\nabla
_zG_T(z,x)||\nabla^* \eta_L(x)||\nabla^* \phi_T(x)|^2 \,dx\biggr)^2\,dz
\biggr\rangle\\
\label{eq:main:term12}
&&\qquad\quad{} +\biggl\langle\int_{\Z^d}\biggl(\int_{\Z^d}|\nabla_zG_T(z,x)||\nabla^*
\eta_L(x)||\nabla\phi_T(z)|^2 \,dx\biggr)^2\,dz \biggr\rangle\\
\label{eq:main:term13}
&&\qquad\quad{} +\biggl\langle\int_{\Z^d}\biggl(\int_{\Z^d}|\nabla_zG_T(z,x)||\nabla^*
\eta_L(x)| \,dx\biggr)^2\,dz \biggr\rangle\\
\label{eq:main:term2}
&&\qquad\quad{} +\biggl\langle\int_{\Z^d}\eta_L(z)^2\bigl(|\nabla\phi_T(z)|^2+1\bigr)^2\,dz
\biggr\rangle.
\end{eqnarray}
Indeed, inserting (\ref{eq:theo-step2-1}) in (\ref
{eq:main-var-start}) yields
\begin{eqnarray*}
&&\operatorname{var}[\llangle T^{-1}\phi_T^2+(\nabla\phi_T+\xi
)\cdot A(\nabla\phi_T+\xi) \rrangle_L]\\
&&\qquad\lesssim \biggl\langle\sum_{e} \biggl(\int_{\Z^d}|\nabla_zG_T(z,x)|
|\nabla^* \eta_L(x)|\bigl(|\nabla^* \phi_T(x)|^2+|\nabla\phi
_T(z)|^2+1\bigr) \,dx \biggr)^2 \biggr\rangle\\
&&\qquad\quad{} +\biggl\langle\sum_{e} \eta_L^2(z)\bigl(|\nabla\phi_T(z)|^2+1\bigr)^2 \biggr\rangle.
\end{eqnarray*}
We then use Young's inequality in the first term of the right-hand side
of this
inequality and we replace the sum $\sum_e$ over edges
$[z,z+\ee_i]$ by $d$ times the sum over $z\in\Z^d$ to establish this step.

It now remains to estimate the terms (\ref{eq:main:term11}),
(\ref{eq:main:term12}), (\ref{eq:main:term13}) and (\ref
{eq:main:term2}) to conclude the proof of the theorem.

\textit{Step} 4. Estimate of (\ref{eq:main:term2}):
%
%
\begin{equation}\label{eq:theo-step4}
\biggl\langle\int_{\Z^d}\eta_L(z)^2\bigl(|\nabla\phi_T(z)|^2+1\bigr)^2\,dz \biggr\rangle
\lesssim
\mu_d(T)^q L^{-d}.
\end{equation}
Indeed, by stationarity we have
\[
\langle|\nabla\phi_T(z)|^4 \rangle\lesssim\sum_{i=1}^d \langle
|\phi_T(z+\ee_i)|^4+|\phi_T(z)|^4 \rangle=2d\langle\phi_T(0)^4
\rangle,
\]
so that
\begin{eqnarray*}
\biggl\langle\int_{\Z^d}\eta_L(z)^2\bigl(|\nabla\phi_T(z)|^2+1\bigr)^2 \,dz \biggr\rangle
&\lesssim& \biggl\langle\int_{\Z^d}\eta_L(z)^2\bigl(|\nabla\phi_T(z)|^4+1\bigr)
\,dz \biggr\rangle\\
&=& \int_{\Z^d}\eta_L(z)^2\bigl(\langle|\nabla\phi_T(z)|^4 \rangle+1\bigr)
\,dz \\
&\lesssim& \bigl(\langle\phi_T(0)^4 \rangle+1\bigr)\int_{\Z^d}\eta_L(z)^2 \,dz.
\end{eqnarray*}
On the one hand, it follows from Proposition \ref{prop:main} that
\[
\langle\phi_T(0)^4 \rangle\lesssim\mu_d(T)^q,
\]
with $q=\gamma(4)$.
On the other hand, it follows from (\ref{eq:eta}) that
\[
\int_{\Z^d}\eta_L(z)^2 \,dz \lesssim L^{-d}.
\]
This establishes Step 4.

\textit{Step} 5. Estimate of (\ref{eq:main:term13}):
%
%
\begin{equation}\label{eq:theo-step5}
\biggl\langle\int_{\Z^d}\biggl(\int_{\Z^d}|\nabla_zG_T(z,x)||\nabla^* \eta
_L(x)| \,dx\biggr)^2\,dz \biggr\rangle\lesssim\mu_d(T)^q L^{-d}.
\end{equation}
We expand the square
\begin{eqnarray*}
&&\biggl\langle\int_{\Z^d}\biggl(\int_{\Z^d}|\nabla_zG_T(z,x)||\nabla^* \eta
_L(x)| \,dx\biggr)^2\,dz \biggr\rangle\\
&&\qquad= \biggl\langle\int_{\Z^d}\int_{\Z^d}\int_{\Z^d}|\nabla^* \eta
_L(x)||\nabla^* \eta_L(x')||\nabla_zG_T(z,x)||\nabla_zG_T(z,x')|
\,dx \,dx' \,dz \biggr\rangle\\
&&\qquad= \int_{\Z^d}\int_{\Z^d} |\nabla^* \eta_L(x)||\nabla^* \eta_L(x')|
\int_{\Z^d}\langle|\nabla_zG_T(z,x)||\nabla_zG_T(z,x')| \rangle
\,dz \,dx \,dx'.
\end{eqnarray*}
We then use Cauchy--Schwarz' inequality in probability and the
stationarity of $G_T$:
\begin{eqnarray*}
&&\langle|\nabla_zG_T(z,x)||\nabla_zG_T(z,x')| \rangle\\
&&\qquad\leq\langle|\nabla_zG_T(z,x)|^2 \rangle^{1/2}\langle|\nabla
_zG_T(z,x')|^2 \rangle^{1/2} \\
&&\qquad=\langle|\nabla_zG_T(z-x,0)|^2 \rangle^{1/2}\langle|\nabla
_zG_T(z-x',0)|^2 \rangle^{1/2} .
\end{eqnarray*}
Hence, with the notation
\[
h(y):=\langle|\nabla_yG_T(y,0)|^2 \rangle^{1/2},
\]
we have by definition of $\eta_L$:
\begin{eqnarray*}
&&\biggl\langle\int_{\Z^d}\biggl(\int_{\Z^d}|\nabla_zG_T(z,x)||\nabla^* \eta
_L(x)| \,dx\biggr)^2\,dz \biggr\rangle\\
&&\qquad\leq \int_{\Z^d}\int_{\Z^d} |\nabla^* \eta_L(x)||\nabla^*
\eta_L(x')|
\int_{\Z^d}h(z-x)h(z-x') \,dz \,dx \,dx' \\
&&\qquad\lesssim L^{-2(d+1)} \int_{|x|\leq L}\int_{|x'|\leq L} \int_{\Z^d}
h(z-x)h(z-x') \,dz \,dx \,dx' \\
&&\qquad= L^{-2(d+1)} \int_{|x|\leq L}\int_{|x'|\leq L} \int_{\Z^d}
h(z')h(z'+x-x') \,dz' \,dx \,dx' \\
&&\qquad\leq L^{-d-2} \int_{|y|\leq2L} \int_{\Z^d} h(z')h(z'-y) \,dz' \,dy.
\end{eqnarray*}
We note that
\[
\int_{R<|y| \leq2R}h^2(y) \,dy = \biggl\langle\int_{R< |y| \leq
2R}|\nabla_yG_T(y,0)|^2 \,dy \biggr\rangle.
\]
On the one hand, for $R\gg1$ we have according to Lemma \ref
{lem:int-grad} (for $q=2$)
\begin{eqnarray*}
\mbox{for }d=2\qquad \int_{R< |y| \leq2R}h^2(y) \,dy
&\lesssim&
R^{2-2}\min\bigl\{1,\sqrt{T}R^{-1}\bigr\}^2\\
&=&\min\bigl\{1,\sqrt{T}R^{-1}\bigr\}^2,\\
\mbox{for }d>2\qquad \int_{R< |y| \leq2R}h^2(y) \,dy
&\lesssim&
R^d(R^{1-d})^2\\
&=&R^{2-d}.
\end{eqnarray*}
On the other hand, for $R\sim1$, Corollary \ref{coro:unif-bound-grad} implies
\[
\mbox{for }d\geq2\qquad \int_{|y| \leq R}h^2(y) \,dy \lesssim1 .
\]
Hence, we are in position to apply Lemma \ref{lem:hh}, which yields as desired
\[
\int_{|y|\leq2L} \int_{\Z^d} h(z')h(z'-y) \,dz' \,dy \lesssim
L^2\mu_d(T).
\]
Note that for $d=2$, we have used the elementary fact that $\max\{
1,\ln\sqrt{T}L^{-1}\}\lesssim\ln T$.

\textit{Step} 6. Estimate of (\ref{eq:main:term12}):
%
%
\begin{equation}\label{eq:theo-step6}\qquad\quad
\biggl\langle\int_{\Z^d}\biggl(\int_{\Z^d}|\nabla_zG_T(z,x)||\nabla^* \eta
_L(x)||\nabla\phi_T(z)|^2 \,dx\biggr)^2\,dz \biggr\rangle\lesssim\mu_d(T)^q L^{-d}.
\end{equation}
As in Step 5,
\begin{eqnarray*}
&&\biggl\langle\int_{\Z^d}\biggl(\int_{\Z^d}|\nabla_zG_T(z,x)||\nabla^* \eta
_L(x)||\nabla\phi_T(z)|^2 \,dx\biggr)^2\,dz \biggr\rangle\\
&&\qquad= \int_{\Z^d}\int_{\Z^d} |\nabla^* \eta_L(x)||\nabla^*
\eta_L(x')|\\
&&\qquad\quad\hspace*{32.6pt}{}\times
\int_{\Z^d}\langle|\nabla_zG_T(z,x)||\nabla_zG_T(z,x')||\nabla
\phi_T(z)|^4 \rangle \,dz \,dx \,dx'.
\end{eqnarray*}
This time, we use H\"older's inequality with $(p,p,\frac{p}{p-2})$ in
probability (where $p>2$ is the exponent in Lemma \ref{lem:int-grad}):
\begin{eqnarray*}
&&\langle|\nabla_zG_T(z,x)||\nabla_zG_T(z,x')||\nabla\phi_T(z)|^4
\rangle\\
&&\qquad\leq\langle|\nabla_zG_T(z,x)|^p \rangle^{{1/p}}\langle
|\nabla_zG_T(z,x')|^p \rangle^{{1/p}}\bigl\langle|\nabla\phi
_T(z)|^{{4p}/({p-2})} \bigr\rangle^{({p-2})/{p}}.
\end{eqnarray*}
By stationarity of $G_T$ and $\phi_T$, we obtain with Proposition \ref
{prop:main}
\begin{eqnarray*}
&&\langle|\nabla_zG_T(z,x)||\nabla_zG_T(z,x')||\nabla\phi_T(z)|^4
\rangle\\
&&\qquad\lesssim\mu_d(T)^q \langle|\nabla_zG_T(z-x,0)|^p \rangle^{
{1/p}}\langle|\nabla_zG_T(z-x',0)|^p \rangle^{{1/p}}.
\end{eqnarray*}
Hence, with the notation
\[
h(y):=\langle|\nabla_yG_T(y,0)|^p \rangle^{1/p},
\]
by definition of $\eta_L$:
\begin{eqnarray*}
&&\biggl\langle\int_{\Z^d}\biggl(\int_{\Z^d}|\nabla_zG_T(z,x)||\nabla^* \eta
_L(x)||\nabla\phi_T(z)|^2 \,dx\biggr)^2\,dz \biggr\rangle\\
&&\qquad\lesssim \mu_d(T)^q \int_{\Z^d}\int_{\Z^d} |\nabla^* \eta
_L(x)||\nabla
^* \eta_L(x')| \int_{\Z^d}h(z-x)h(z-x') \,dz \,dx \,dx' \\
&&\qquad\lesssim\mu_d(T)^q L^{-d-2} \int_{|y|\leq2L} \int_{\Z^d}
h(z')h(z'-y) \,dz' \,dy.
\end{eqnarray*}
As in Step 5, we shall establish that for $R\gg1$
%
%
\begin{eqnarray}\label{eq:theo-step6-1}
&&\mbox{for }d=2\qquad \int_{R< |y| \leq2R}h^2(y) \,dy \lesssim
\min\bigl\{1,\sqrt{T}R^{-1}\bigr\}^2,\nonumber\\[-8pt]\\[-8pt]
&&\mbox{for }d>2\qquad \int_{R< |y| \leq2R}h^2(y) \,dy \lesssim
R^{2-d},\nonumber
\end{eqnarray}
and for $R\sim1$
%
%
\begin{equation}\label{eq:theo-step6-2}
\mbox{for }d\geq2\qquad \int_{|y| \leq R}h^2(y) \,dy \lesssim1.
\end{equation}
Once this is done, Lemma \ref{lem:hh} implies as desired
\[
\int_{|y|\leq2L} \int_{\Z^d} h(z')h(z'-y) \,dz' \,dy \lesssim
L^2\mu_d(T),
\]
using in addition that $\max\{1,\ln\sqrt{T}L^{-1}\}\lesssim\ln T$
for $d=2$.
As above, (\ref{eq:theo-step6-2}) is a direct consequence of
Corollary \ref{coro:unif-bound-grad}.
We now deal with (\ref{eq:theo-step6-1}).
Note that according to Lemma \ref{lem:int-grad}, we have for $R\gg1$
%
%
\begin{eqnarray}\label{eq:theo-step6-3}
&&\mbox{for }d=2\qquad \int_{R< |y| \leq2R}h^p(y) \,dy \lesssim
R^{2-p}\min\bigl\{1,\sqrt{T}R^{-1}\bigr\}^p,\nonumber\\[-8pt]\\[-8pt]
&&\mbox{for }d>2\qquad \int_{R< |y| \leq2R}h^p(y) \,dy \lesssim
R^d(R^{1-d})^p.\nonumber
\end{eqnarray}
We now argue that this yields (\ref{eq:theo-step6-1}).
Indeed, by Jensen's inequality
\begin{eqnarray*}
&&\biggl( R^{-d}\int_{R< |x|\leq2R}h^2(x)\,dx\biggr)^{1/2} \\
&&\qquad \leq
\biggl(
R^{-d}\int_{R< |x|\leq2R}h^p(x)\,dx\biggr)^{1/p} \\
&&\hspace*{-5.4pt}\qquad\stackrel{\mbox{\fontsize{8.36pt}{10.36pt}\selectfont{(\ref
{eq:theo-step6-3})}}}{\lesssim} \cases{
\bigl( R^{-2}R^{2-p}\min\bigl\{1,\sqrt{T}R^{-1}\bigr\}^p\bigr)^{1/p}, &\quad $d=2$, \vspace*{2pt}\cr
( R^{-d}R^d(R^{1-d})^p)^{1/p}, &\quad $d>2$,}
\\
&&\qquad= \cases{
R^{-1} \min\bigl\{1,\sqrt{T}R^{-1}\bigr\}, &\quad $d=2$, \cr
R^{1-d}, &\quad $d>2$,}
\end{eqnarray*}
which implies (\ref{eq:theo-step6-1}).

\textit{Step} 7. Estimate of (\ref{eq:main:term11}):
%
%
\begin{equation}\label{eq:theo-step7}
\biggl\langle\int_{\Z^d}\biggl(\int_{\Z^d}|\nabla_zG_T(z,x)||\nabla^* \eta
_L(x)||\nabla^* \phi_T(x)|^2 \,dx\biggr)^2\,dz \biggr\rangle\lesssim\mu
_d(T)^q L^{-d}.\hspace*{-28pt}
\end{equation}
As in Steps 5 and 6,
\begin{eqnarray*}
&&\biggl\langle\int_{\Z^d}\biggl(\int_{\Z^d}|\nabla_zG_T(z,x)||\nabla^* \eta
_L(x)||\nabla^* \phi_T(x)|^2 \,dx\biggr)^2\,dz \biggr\rangle\\
&&\qquad= \int_{\Z^d}\int_{\Z^d} |\nabla^* \eta_L(x)||\nabla^* \eta
_L(x')| \\
&&\qquad\quad\hspace*{32.3pt}{}\times
\int_{\Z^d}\langle|\nabla_zG_T(z,x)||\nabla_zG_T(z,x')|\\
&&\qquad\quad{}\hspace*{65.4pt}\times|\nabla
^* \phi_T(x)|^2 |\nabla^* \phi_T(x')|^2 \rangle \,dz
\,dx \,dx'.
\end{eqnarray*}
H\"older's inequality with $(p,p,\frac{2p}{p-2},\frac{2p}{p-2})$ in
probability (where $p>2$ is the exponent in Lemma \ref{lem:int-grad})
then yields
\begin{eqnarray*}
&&\langle|\nabla_zG_T(z,x)||\nabla_zG_T(z,x')||\nabla^* \phi
_T(x)|^2 |\nabla^* \phi_T(x')|^2 \rangle\\
&&\qquad\leq\langle|\nabla_zG_T(z,x)|^p \rangle^{{1/p}}\langle
|\nabla_zG_T(z,x')|^p \rangle^{{1/p}}\\
&&\qquad\quad\hspace*{0pt}{}\times\bigl\langle|\nabla^* \phi
_T(x)|^{{4p}/({p-2})} \bigr\rangle^{({p-2})/({2p})}\bigl\langle|\nabla^*
\phi_T(x')|^{{4p}/({p-2})} \bigr\rangle^{({p-2})/({2p})}.
\end{eqnarray*}
The stationarity of $G_T$ and $\phi_T$, and Proposition \ref
{prop:main} show
\begin{eqnarray*}
&&\langle|\nabla_zG_T(z,x)||\nabla_zG_T(z,x')||\nabla^* \phi
_T(x)|^2|\nabla^* \phi_T(x')|^2 \rangle\\
&&\qquad\lesssim \mu_d(T)^q \langle|\nabla_zG_T(z-x,0)|^p \rangle^{
{1/p}}\langle|\nabla_zG_T(z-x',0)|^p \rangle^{{1/p}}.
\end{eqnarray*}
We may now conclude as in Step 6.

The theorem follows from the combination of Step 3 with (\ref
{eq:theo-step4}), (\ref{eq:theo-step5}), (\ref{eq:theo-step6}) and
(\ref{eq:theo-step7}).

\textit{Step} 8. Extension to the energy density of the corrector
field for $d>2$.

Let $A_{L,\infty}$ be defined by
\[
\xi\cdot A_{L,\infty}\xi:=\int_{\Z^d} \bigl(\nabla\phi(x)+\xi\bigr)\cdot
A(x)\bigl(\nabla\phi(x)+\xi\bigr)\mu_L(x) \,dx,
\]
for all $L\gg1$.
The claim is
\[
\operatorname{var}[\xi\cdot A_{L,\infty}\xi] \lesssim L^{-d} ,
\]
for $d>2$.
It is proved using (\ref{eq:estim-var-hom}) provided we show
%
%
\begin{equation} \label{eq:pr-th-step8-0}
\operatorname{var}[\xi\cdot A_{L,\infty}\xi] \leq\liminf_{T\to
\infty}\operatorname{var}[\xi\cdot A_{L,T}\xi].
\end{equation}
As we shall prove, the following two convergences hold:
%
%
\begin{eqnarray} \label{eq:pr-th-step8-1}
&&\biggl\langle\int_{\Z^d}\bigl(\xi+\nabla\phi_T(x)\bigr)\cdot A(x)\bigl(\xi+\nabla
\phi_T(x)\bigr)\mu_L(x)\,dx \biggr\rangle\nonumber\\
&&\qquad \to \biggl\langle\int_{\Z^d}\bigl(\xi+\nabla\phi(x)\bigr)\cdot A(x)\bigl(\xi
+\nabla\phi(x)\bigr)\mu_L(x)\,dx \biggr\rangle\\
&&\qquad= \xi\cdot A_\ho\xi,\nonumber
\end{eqnarray}
which in fact amounts to $\langle(\xi+\nabla\phi_T)\cdot A(\xi
+\nabla\phi_T) \rangle\to\langle(\xi+\nabla\phi)\cdot A(\xi
+\nabla\phi) \rangle$
by stationarity, and
%
%
\begin{eqnarray}\label{eq:pr-th-step8-2}
&&\int_{\Z^d}\bigl(\xi+\nabla\phi_T(x)\bigr)\cdot A(x)\bigl(\xi+\nabla
\phi
_T(x)\bigr)\mu_L(x)\,dx\nonumber\\
&&\qquad \rightharpoonup \int_{\Z^d}\bigl(\xi+\nabla\phi(x)\bigr)\cdot A(x)\bigl(\xi
+\nabla\phi(x)\bigr)\mu_L(x)\,dx\\
&&\eqntext{\mbox{weakly in probability.}}
\end{eqnarray}
We may now conclude the proof of (\ref{eq:pr-th-step8-0}). Expanding
the variance, one has
\begin{eqnarray*}
\operatorname{var}[\xi\cdot A_{L,T}\xi] &=& \biggl\langle\biggl(\int_{\Z
^d}\bigl(\xi+\nabla\phi_T(x)\bigr)\cdot A(x)\bigl(\xi+\nabla\phi_T(x)\bigr)\mu
_L(x)\,dx\biggr)^2 \biggr\rangle\\
&&{} -\biggl\langle\int_{\Z^d}\bigl(\xi+\nabla\phi_T(x)\bigr)\cdot A(x)\bigl(\xi+\nabla
\phi_T(x)\bigr)\mu_L(x)\,dx \biggr\rangle^2.
\end{eqnarray*}
By (\ref{eq:pr-th-step8-1}), the second term of the right-hand side
converges to
$(\xi\cdot A_\ho\xi)^2$ as $T\to\infty$,
whereas by lower-semicontinuity of quadratic functionals, (\ref
{eq:pr-th-step8-2}) implies that
\begin{eqnarray*}
&&\biggl\langle\biggl(\int_{\Z^d}\bigl(\xi+\nabla\phi(x)\bigr)\cdot A(x)\bigl(\xi+\nabla
\phi(x)\bigr)\mu_L(x)\,dx\biggr)^2 \biggr\rangle\\
&&\qquad\leq\liminf_{T\to\infty}
\biggl\langle\biggl(\int_{\Z^d}\bigl(\xi+\nabla\phi_T(x)\bigr)\cdot A(x)\bigl(\xi+\nabla
\phi_T(x)\bigr)\mu_L(x)\,dx\biggr)^2 \biggr\rangle,
\end{eqnarray*}
which shows (\ref{eq:pr-th-step8-0}).

It remains to prove (\ref{eq:pr-th-step8-1}) and (\ref{eq:pr-th-step8-2}).
Note that by stationarity, (\ref{eq:pr-th-step8-1}) is a consequence of
%
%
\begin{equation}\label{eq:pr-th-step8-3}
{\lim_{T\to\infty}}|A_T-A_\ho| = 0,
\end{equation}
for all $d\geq2$, where $\xi\cdot A_T\xi:=\langle(\xi+\nabla\phi
_T)\cdot A(\xi+\nabla\phi_T) \rangle$.
Starting point for (\ref{eq:pr-th-step8-3}) is the definition of $A_T$
and $A_\ho$ from which we deduce
%
%
\begin{eqnarray}\label{eq:AT-Ahom}
\xi\cdot(A_T-A_\ho)\xi&=& \langle(\xi+\nabla\phi_T)\cdot A(\xi
+\nabla\phi_T)-(\xi+\nabla\phi) \cdot A(\xi+\nabla\phi)
\rangle\hspace*{-35pt}\nonumber\\
&=&\langle\xi\cdot A(\nabla\phi_T-\nabla\phi) \rangle+\langle
\nabla\phi_T\cdot A(\xi+\nabla\phi_T) \rangle\\
&&{}-\langle\nabla
\phi\cdot A(\xi+\nabla\phi) \rangle.\nonumber
\end{eqnarray}
Let us treat each term separately.
For the second term, we shall argue that (\ref{eq:app-corr}) yields:
for every stationary field $\zeta\dvtx\Z^d\to\R$
such that $\langle\zeta^2 \rangle<\infty$, one has
%
%
\begin{equation}\label{eq:VF-proba}
T^{-1}\langle\phi_T\zeta\rangle+\langle\nabla\zeta\cdot A(\xi
+\nabla\phi_T) \rangle= 0,
\end{equation}
so that one may replace the second term of the right-hand side of (\ref
{eq:AT-Ahom}) by $-T^{-1}\langle\phi_T^2 \rangle$.
For the first term, we shall use the following weak convergence of
$\nabla\phi_T(x)$ to $\nabla\phi(x)$
in probability: for every random variable $\chi$ taking values in $\R^d$
with $\langle|\chi|^2 \rangle<\infty$, one has for all $x\in\Z^d$,
%
%
\begin{equation} \label{eq:kun-1}
\lim_{T\to\infty}\bigl\langle\chi\cdot\bigl(\nabla\phi_T(x)-\nabla\phi
(x)\bigr) \bigr\rangle= 0,
\end{equation}
so that taking $x=0$ and $\chi\equiv A(0)\xi$ shows that the first
term in the right-hand side of (\ref{eq:AT-Ahom}) vanishes as
$T\uparrow
\infty$.
For the last term, combining (\ref{eq:kun-1}) and (\ref{eq:VF-proba}),
we will prove
%
%
\begin{equation}\label{eq:kun-2}
\langle\nabla\phi\cdot A(\xi+\nabla\phi) \rangle= 0.
\end{equation}
We directly draw the conclusion: the combination of (\ref{eq:AT-Ahom}),
(\ref{eq:kun-1}),
(\ref{eq:VF-proba}) and (\ref{eq:kun-2})
shows that
\[
\limsup_{T\to\infty} |\xi\cdot(A_T-A_\ho)\xi| = \limsup
_{T\to\infty
}T^{-1}\langle\phi_T^2 \rangle,
\]
which implies (\ref{eq:pr-th-step8-3}) by Proposition \ref{prop:main}.

We give the arguments for (\ref{eq:VF-proba}), (\ref{eq:kun-1}) and
(\ref{eq:kun-2}) for the reader's convenience
(we could also directly appeal to \cite{Kunnemann-83}).
Multiplying the defining equation for $\phi_T$ by $\zeta$ yields
%
%
\begin{equation}\label{eq:VF-reel}
T^{-1}(\phi_T\zeta)(z)-\bigl(\nabla^*\cdot A(\xi+\nabla\phi
_T)\bigr)(z)\zeta(z)
= 0.
\end{equation}
We then use the discrete Leibniz rule in the form
%
%
\begin{eqnarray}\label{eq:disc-leib1}\quad
\nabla^*\cdot\bigl(\zeta A(\xi+\nabla\phi_T)\bigr)(z) &=& \bigl(\nabla^*\cdot
A(\xi+\nabla\phi_T)\bigr)(z)\zeta(z)\nonumber\\[-8pt]\\[-8pt]
&&{} + \sum_{j=1}^d
\nabla^*_j\zeta(z)[A(\xi+\nabla\phi_T)(z-\ee_j)]_j.\nonumber
\end{eqnarray}
Since $\nabla\phi_T$, $\zeta$ and $A$ are jointly stationary random
fields, the expectation of the left-hand side of
(\ref{eq:disc-leib1}) vanishes, and
%
%
\begin{eqnarray}\label{eq:expec1}\qquad
\bigl\langle\bigl(\nabla^*\cdot A(\xi+\nabla\phi_T)\bigr)(z)\zeta(z) \bigr\rangle
&=&-\Biggl\langle\sum_{j=1}^d \nabla^*_j\zeta(z)\bigl[A\bigl(\xi+\nabla\phi
_T(z-\ee_j)\bigr)\bigr]_j \Biggr\rangle
\nonumber\\[-8pt]\\[-8pt]
&=& -\langle\nabla\zeta\cdot A (\xi+\nabla\phi_T) \rangle,\nonumber
\end{eqnarray}
noting that $\nabla^*_j\zeta(z)=\nabla_j \zeta(z-\ee_j)$.
We then take the expectation of (\ref{eq:VF-reel}) and use (\ref
{eq:expec1}) to obtain (\ref{eq:VF-proba}).

We recall the standard a priori estimate which one derives
from (\ref{eq:VF-proba}):
\[
\langle T^{-1}\phi_T(x)^2+|\nabla\phi_T(x)|^2 \rangle\lesssim1.
\]
Since the left-hand side does not depend on $x$ by stationarity, there exists
$g\dvtx\Z^d \to\R^d$
such that up to extraction, $\nabla\phi_T(x)$ converges to $g(x)$
weakly in probability for
all $x\in\Z^d$.
By construction, $g$ is a gradient field, and is jointly stationary
with $A$.
By the boundedness of $\langle T^{-1}\phi_T^2 \rangle^{1/2}$, one may
pass to
the limit in (\ref{eq:VF-proba}),
and obtain for every stationary field $\zeta$
%
%
\begin{equation}\label{eq:kun-3}
\langle\nabla\zeta\cdot A(\xi+\nabla\phi) \rangle= 0.
\end{equation}
As noticed by K\"unnemann in \cite{Kunnemann-83}, this characterizes
the gradient of the corrector, so that $g\equiv\nabla\phi$.
This proves (\ref{eq:kun-1}) by definition of weak convergence in probability.

We then use (\ref{eq:kun-3}) for $\zeta=\phi_T$ and pass to the limit
$T\uparrow\infty$ in (\ref{eq:kun-3})
by the weak convergence (\ref{eq:kun-1}). This proves (\ref{eq:kun-2}).

We finally turn to the proof of (\ref{eq:pr-th-step8-2}).
By definition, (\ref{eq:pr-th-step8-2}) is proved if for all bounded
random variables $\chi$,
%
%
\begin{eqnarray}\label{eq:th-step8-wc}
&&\lim_{T\to\infty}\biggl\langle\chi\int_{\Z^d}\bigl(\xi+\nabla\phi
_T(x)\bigr)\cdot A(x)\bigl(\xi+\nabla\phi_T(x)\bigr)\mu_L(x)\,dx \biggr\rangle
\nonumber\\[-8pt]\\[-8pt]
&&\qquad= \biggl\langle\chi\int_{\Z^d}\bigl(\xi+\nabla\phi(x)\bigr)\cdot A(x)\bigl(\xi
+\nabla\phi(x)\bigr)\mu_L(x)\,dx \biggr\rangle.\nonumber
\end{eqnarray}
W.l.o.g. we may assume that $\chi$ takes values in $[0,1]$.
By lower-semicontinuity of quadratic functionals in probability, and
since $\chi\geq0$, the weak
convergence (\ref{eq:kun-1}) of $\nabla\phi_T(x)$ to $\nabla\phi(x)$
in $L^2$ in probability for all $x\in\Z^d$
yields
\begin{eqnarray*}
&&\liminf_{T\to\infty}\biggl\langle\chi\int_{\Z^d}\bigl(\xi+\nabla\phi
_T(x)\bigr)\cdot A(x)\bigl(\xi+\nabla\phi_T(x)\bigr)\mu_L(x)\,dx \biggr\rangle\\
&&\qquad=\int_{\Z^d}\mu_L(x) \Bigl(\liminf_{T\to\infty}\bigl\langle\chi\bigl(\xi
+\nabla\phi_T(x)\bigr)\cdot A(x)\bigl(\xi+\nabla\phi_T(x)\bigr) \bigr\rangle\Bigr)\,dx \\
&&\qquad\geq \int_{\Z^d}\mu_L(x) \bigl\langle\chi\bigl(\xi+\nabla\phi
_T(x)\bigr)\cdot A(x)\bigl(\xi+\nabla\phi_T(x)\bigr) \bigr\rangle \,dx\\
&&\qquad=\biggl\langle\chi\int_{\Z^d}\bigl(\xi+\nabla\phi(x)\bigr)\cdot A(x)\bigl(\xi
+\nabla\phi(x)\bigr)\mu_L(x)\,dx \biggr\rangle.
\end{eqnarray*}
Likewise,
\begin{eqnarray*}
&&\liminf_{T\to\infty}\biggl\langle(1-\chi)\int_{\Z^d}\bigl(\xi+\nabla
\phi_T(x)\bigr)\cdot A(x)\bigl(\xi+\nabla\phi_T(x)\bigr)\mu_L(x)\,dx \biggr\rangle\\
&&\qquad\geq \biggl\langle(1-\chi)\int_{\Z^d}\bigl(\xi+\nabla\phi(x)\bigr)\cdot
A(x)\bigl(\xi+\nabla\phi(x)\bigr)\mu_L(x)\,dx \biggr\rangle
\end{eqnarray*}
since $1-\chi\geq0$.
Combined with the convergence of the expectation (\ref
{eq:pr-th-step8-1}) and the trivial identity
$1=\chi+(1-\chi)$, these two inequalities imply (\ref{eq:th-step8-wc})
for $\chi$ taking values in $[0,1]$,
and therefore (\ref{eq:pr-th-step8-2}) as desired.


\section{Proofs of the estimates on the Green functions}\label
{sec:proofs-Green}

Before addressing the proofs proper, let us make a general comment.
In what follows, we shall replace the classical Leibniz rule by
its discrete counterpart.
Although they are essentially the same, the expressions that appear are
more intricate
in the discrete case.
In order to keep the proofs clear, we first present the arguments using
the classical Leibniz rule (though it does not hold at the discrete level)
and we later give a separate argument to show that the various
results still hold with the true discrete version.

\subsection{\texorpdfstring{Proof of Lemma \protect\ref{L11a}}{Proof of Lemma 2.8.}}

Without loss of generality, we may assume $y=0$ and suppress
the $y$-dependance
of $G_T$ in our notation.
We will first give the proof in the continuum case (i.e., using the
classical Leibniz rule) and then
sketch the modifications arising from the discreteness.

We first argue that for any $d$,
%
%
\begin{equation}\label{L11.1}
T^{-1}\int_{\Z^d} G_{T,M}^2 \,dx+\int_{\Z^d}|\nabla G_{T,M}|^2 \,dx
\lesssim M,
\end{equation}
where for $0<M<\infty$, $G_{T,M}$
denotes the following truncated version of $G_T$
\[
G_{T,M} = \min\{G_T,M\} \geq0.
\]
Indeed, we consider $T^{-1}G_T-\nabla^*\cdot A\nabla G_T=\delta$ in its
weak form, that is,
%
%
\begin{equation}\label{L11.0}
T^{-1}\int_{\Z^d} \zeta G_T \,dx+\int_{\Z^d}\nabla\zeta\cdot
A\nabla G_T \,dx = \zeta(0)
\end{equation}
and select $\zeta=G_{T,M}$. Since $G_{T,M} G_T\ge G_{T,M}^2$ and
provided that
$\nabla G_{T,M}\cdot A\nabla G_T\geq\nabla G_{T,M}\cdot A\nabla
G_{T,M}$, we obtain
(\ref{L11.1}) by uniform ellipticity. Indeed, since $A$ is diagonal,
\begin{eqnarray*}
&&{\nabla G_{T,M}\cdot A\nabla G_T (x)} \\
&&\qquad = \sum_{i=1}^d a(x+\ee_i,x)
\bigl(G_{T,M}(x+\ee_i)-G_{T,M}(x)\bigr)\bigl(G_{T}(x+\ee_i)-G_{T}(x)\bigr) \\
&&\qquad\geq \sum_{i=1}^d a(x+\ee_i,x) \bigl(G_{T,M}(x+\ee_i)-G_{T,M}(x)\bigr)^2 \\
&&\qquad \geq \alpha|\nabla G_{T,M}(x)|^2.
\end{eqnarray*}

\textit{Step} 1. Proof of (i) for $d>2$.

Following \cite{Grueter-Widman-82}, Theorem 1.1,
we argue that (\ref{L11.1}) implies a weak-$L^{{d}/({d-2})}$
estimate, that is,
%
%
\begin{equation}\label{L11.4}
{\mathcal L}_d(\{G_T\ge M\}) \lesssim M^{-{d}/({d-2})}.
\end{equation}
For this purpose, we appeal to Sobolev's inequality in $\Z^d$,
that is,
\[
\biggl(\int_{\Z^d} G_{T,M}^{2d/(d-2)} \,dx\biggr)^{({d-2})/({2d})}
\lesssim
\biggl(\int_{\Z^d}|\nabla G_{T,M}|^2 \,dx\biggr)^{1/2},
\]
which is a consequence of \cite{Zhou-93}, Lemma 2.1 (when ``$n\to
\infty
$''), or
\cite{Delmotte-97}, Theorem~4.4 (when ``$r\to\infty$'').
Via Chebyshev's inequality and (\ref{L11.1}), this yields
\[
M {\mathcal L}_d(\{G_T\ge M\})^{({d-2})/({2d})} \lesssim M^{1/2},
\]
which is (\ref{L11.4}).

We now argue that the weak-$L^{{d}/({d-2})}$ estimate
(\ref{L11.4})
in $\Z^d$ yields a strong $L^q$-estimate on balls $\{|x|\le R\}$
for $1\leq q<\frac{d}{d-2}$. More precisely, we have
%
%
\begin{equation}\label{L11.7}
\int_{|x|\le R}G_T^q \,dx \lesssim R^d (R^{2-d})^q.
\end{equation}
Indeed, we have on the one hand
%
%
\begin{eqnarray}\label{L11.5}
\int_{G_T>M} G_T^q \,dx
&=& q\int_M^\infty{\mathcal L}_d(\{G_T>M'\}) {M'}^{q-1} \,dM'\nonumber\\
&&{} +M^q
{\mathcal L}_d(\{|G_T|> M\})\\
&\stackrel{\mbox{\fontsize{8.36pt}{10.36pt}\selectfont{(\ref
{L11.4})}}}{\lesssim}& M^{q-{d}/({d-2})},\nonumber
\end{eqnarray}
where we have used $q<\frac{d}{d-2}$.
On the other hand, we have trivially
%
%
\begin{equation}\label{L11.6}
\int_{\{G_T\le M\}\cap\{|x|\le R\}}G_T^q \,dx \lesssim R^d M^q.
\end{equation}
With the choice of $M=R^{2-d}$, the combination of (\ref{L11.5})
and (\ref{L11.6}) yields (\ref{L11.7}).

In order to increase the exponent $q$ in (\ref{L11.7}), one combines
a Cacciopoli estimate\footnote{This is the only place where we use the
Leibniz rule.} for monotone functions of $G_T$ with a
Poincar\'e--Sobolev estimate to obtain a ``reverse H\"older'' inequality
(as in the proof of Harnack's inequality, see \cite{Han-Lin-97},
Chapter 4, Method II). We start with the
Cacciopoli estimate, that is,
%
%
\begin{equation}\label{L11.8}
\int_{2R\le|x|\le4R}|\nabla G_T^{q/2}|^2 \,dx
\lesssim
R^{-2}\int_{R\le|x|\le8R} G_T^q \,dx
\end{equation}
for all $1<q<\infty$.
For that purpose, we test (\ref{L11.0}) with $\zeta=\eta^2 G_T^{q-1}$,
where the spatial cut-off function $\eta$ has the properties
%
%
\begin{eqnarray}\label{L11.-1}\qquad
\eta&\equiv&1 \qquad\mbox{in } \{2R\le|x|\le4R\},\nonumber\\[-8pt]\\[-8pt]
\eta&\equiv&0 \qquad\mbox{outside } \{R\le|x|\le8R\},\qquad
|\nabla\eta| \lesssim R^{-1},\qquad 0 \le\eta\le1.\nonumber
\end{eqnarray}
This yields
%
%
\begin{equation}\label{L11.-4}
T^{-1}\int_{\Z^d}\eta^2 G_T^q \,dx+\int_{\Z^d}\nabla(\eta^2
G_T^{q-1})\cdot A\nabla G_T \,dx
= 0.
\end{equation}
Since by the uniform ellipticity of $A$, there exists a generic
constant $C<\infty$
(which only depends on $q$, $\alpha$, $\beta$) such that
\begin{eqnarray*}
&&\nabla(\eta^2 G_T^{q-1})\cdot A\nabla G_T\\
&&\qquad=
(q-1) \eta^2 G_T^{q-2} \nabla G_T\cdot A\nabla G_T
+2 \eta G_T^{q-1} \nabla\eta\cdot A \nabla G_T\\
&&\hspace*{-7pt}\qquad\stackrel{\mathrm{Young}}{\ge}
C^{-1} \eta^2 G_T^{q-2} |\nabla G_T|^2-
C G_T^q |\nabla\eta|^2\\
&&\qquad\gtrsim
C^{-1} \eta^2 |\nabla G_T^{q/2}|^2-
C G_T^q |\nabla\eta|^2,
\end{eqnarray*}
we obtain
\[
\int_{\Z^d}\eta^2 |\nabla G_T^{q/2}|^2 \,dx \lesssim
\int_{\Z^d} G_T^q |\nabla\eta|^2 \,dx.
\]
In view of the properties (\ref{L11.-1}) of $\eta$,
this yields (\ref{L11.8}) for $d>2$.

We now derive the ``reverse H\"older'' inequality
%
%
\begin{eqnarray}\label{L11.10}
&&\biggl(R^{-d}\int_{2R\le|x|\le4R} G_T^{{qd}/({d-2})} \,dx
\biggr)^{
({d-2})/({qd})}\nonumber\\[-8pt]\\[-8pt]
&&\qquad\lesssim
\biggl(R^{-d}\int_{R\le|x|\le8R} G_T^q \,dx\biggr)^{{1/q}}.\nonumber
\end{eqnarray}
For that purpose, we appeal to the Poincar\'e--Sobolev estimate (see
\cite{Zhou-93}, Lem\-ma~2.1, or \cite{Delmotte-97}, Theorem 4.4)
on the annulus $\{2R\le|x|\le4R\}$:
\begin{eqnarray*}
\biggl(R^{-d}\int_{2R\le|x|\le4R}|u|^{{2d}/({d-2})}
\,dx\biggr)^{({d-2})/({2d})}
&\lesssim&
\biggl(R^{2-d}\int_{2R\le|x|\le4R}|\nabla u|^2 \,dx\biggr)^{1/2}\\
&&{}+\biggl(R^{-d}\int_{2R\le|x|\le4R}|u|^{2}
\,dx\biggr)^{1/2}.
\end{eqnarray*}
We apply the latter to $u=G_T^{q/2}$:
\begin{eqnarray*}
\biggl(R^{-d}\int_{2R\le|x|\le4R}G_T^{{qd}/({d-2})}
\,dx\biggr)^{({d-2})/({qd})}
&\lesssim&
\biggl(R^{2-d}\int_{2R\le|x|\le4R}|\nabla G_T^{q/2}|^2 \,dx\biggr)^{1/q}\\
&&{} + \biggl(R^{-d}\int_{2R\le|x|\le4R}G_T^{q} \,dx\biggr)^{1/q}.
\end{eqnarray*}
The combination of this with (\ref{L11.8}) yields (\ref{L11.10}).

We now may conclude in the case of $d>2$: indeed, (\ref{L11.10})
allows us to iteratively increase the integrability $q$ in
multiplicative increments of $\frac{d}{d-2}$ in the estimate (\ref{L11.7}).
Since any $p<\infty$ can be reached in finite multiplicative
increments starting from a $1<q<\frac{d}{d-2}$,
the side effect that the annuli get dyadically larger at
every step does not matter qualitatively (in this sense, the
above argument is much less subtle than the proof of the Harnack
inequality). This proves (\ref{11.16b}).

\textit{Step} 2. Proof of (i) for $d=2$.

We now tackle the case of $d=2$, which in fact amounts
to the $L^1$-BMO estimate
%
%
\begin{equation}\label{11.5}
\biggl( R^{-2}\int_{|x|\le R}|u-\bar u_{\{|x|\le R\}}|^q \,dx\biggr)^{1/q}
\lesssim\int_{\Z^2} |f| \,dx
\end{equation}
for
%
%
\begin{equation}\label{L11.-3}
T^{-1} u-\nabla^*\cdot A \nabla u = f,
\end{equation}
where $\bar u_{\{|x|\le R\}}$ denotes the average of $u$ on the ball of
radius $R$.
We fix an exponent $q<\infty$ and
a radius $1\ll R<\infty$ and assume w.l.o.g.
%
%
\begin{equation}\label{11.13}
\bar u_{\{|x|\le R\}} = 0.
\end{equation}
As in (\ref{L11.1}), we have
%
%
\begin{equation}\label{11.7}
\int_{|x|\le R}|\nabla u_M|^2 \,dx \lesssim M \int_{\Z^2}|f| \,dx.
\end{equation}
As opposed to the case of $d>2$, this is the only time we use
the equation (\ref{L11.-3}).

Estimate (\ref{11.7}) is used in connection with the
Poincar\'e--Sobolev inequality with mean value zero, that is,
\[
\biggl(R^{-2}\int_{|x|\le R}\bigl|u_M-(\overline{u}_M)_{\{|x|\le R\}}\bigr|^s
\,dx\biggr)^{1/s}
\lesssim
\biggl(\int_{|x|\le R}|\nabla u_M|^2 \,dx\biggr)^{1/2},
\]
for any $s<\infty$, which we use once for $s=q$, that is,
%
%
\begin{eqnarray}\label{11.14}\hspace*{32pt}
\biggl(R^{-2}\int_{|x|\le R}\bigl|u_M-(\overline{u}_M)_{\{|x|\le R\}}\bigr|^q
\,dx\biggr)^{1/q}
&\lesssim&
\biggl(\int_{|x|\le R}|\nabla u_M|^2
\,dx\biggr)^{1/2}\nonumber\\[-8pt]\\[-8pt]
&\stackrel{\mbox{\fontsize{8.36pt}{10.36pt}\selectfont{(\ref
{11.7})}}}{\lesssim}&
\biggl(M \int_{\Z^2}|f| \,dx\biggr)^{1/2},\nonumber
\end{eqnarray}
and once for arbitrary $s$ (which we think of being larger than $q$) in
the form
%
%
\begin{eqnarray}\label{11.9}
&&\biggl(R^{-2}\int_{|x|\le R}|u_M|^s \,dx\biggr)^{1/s}\nonumber\\
&&\qquad\lesssim
\biggl(\int_{|x|\le R}|\nabla u_M|^2 \,dx\biggr)^{1/2}
+\bigl|(\overline{u}_M)_{\{|x|\le R\}}\bigr|\\
&&\hspace*{-5.9pt}\qquad\stackrel{\mbox{\fontsize{8.36pt}{10.36pt}\selectfont{(\ref
{11.7})}}}{\lesssim}
\biggl(M \int_{\Z^2}|f| \,dx\biggr)^{1/2}
+\biggl(R^{-2}\int_{|x|\le R}|u|^q \,dx\biggr)^{1/q}.\nonumber
\end{eqnarray}

We use (\ref{11.9}) to estimate the peaks of $u$. More precisely,
we claim that for $s>2q$,
%
%
\begin{eqnarray}\label{11.10}
&&\biggl(R^{-2}\int_{\{|x|\le R\}\cap\{|u|>M\}}|u|^q
\,dx\biggr)^{1/q}\nonumber\\[-8pt]\\[-8pt]
&&\qquad\lesssim
M^{1-s/(2q)}\biggl(\int_{\Z^2}|f| \,dx\biggr)^{s/(2q)}
+M^{1-s/q}
\biggl(R^{-2}\int_{|x|\le R}|u|^q \,dx\biggr)^{s/q^2}.\hspace*{-30pt}\nonumber
\end{eqnarray}
The argument for (\ref{11.10}) is similar to the case of $d>2$:
estimate (\ref{11.9}) yields the weak estimate
\begin{eqnarray*}
&&M \bigl(R^{-2} {\mathcal L}_2(\{|x|\le R\}\cap\{
|u|>M\}
)\bigr)^{1/s}\\
&&\qquad\lesssim
\biggl(M \int_{\Z^2}|f| \,dx\biggr)^{1/2}
+\biggl(R^{-2}\int_{|x|\le R}|u|^q \,dx\biggr)^{1/q},
\end{eqnarray*}
which we rewrite as
%
%
\begin{eqnarray}\label{11.11}
&&R^{-2} {\mathcal L}_2\bigl(\{|x|\le
R\}\cap\{|u|>M\}\bigr)\nonumber\\[-8pt]\\[-8pt]
&&\qquad\lesssim
M^{-s/2} \biggl(\int_{\Z^2}|f| \,dx\biggr)^{s/2}
+M^{-s} \biggl(R^{-2}\int_{|x|\le R}|u|^q \,dx\biggr)^{s/q}.\nonumber
\end{eqnarray}
On the other hand, we have
%
%
\begin{eqnarray}\label{11.12}
\int_{\{|x|\le R\}\cap\{|u|>M\}}|u|^q \,dx
&=&q\int_{M}^\infty{\mathcal L}_2(\{|x|\le R\}\cap\{|u|>M'\})
{M'}^{q-1} \,dM'\hspace*{-32pt} \nonumber\\[-8pt]\\[-8pt]
&&{} +M^q {\mathcal L}_2(\{|x|\le1\}\cap\{|u|>
M\}).\nonumber
\end{eqnarray}
Since $s>2q$, the combination of (\ref{11.11}) and (\ref{11.12}) yields
\begin{eqnarray*}
&&R^{-2}\int_{\{|x|\le R\}\cap\{|u|>M\}}|u|^q \,dx\\
&&\qquad\lesssim
M^{q-s/2} \biggl(\int_{\Z^2}|f| \,dx\biggr)^{s/2}
+M^{q-s} \biggl(R^{-2}\int_{|x|\le R}|u|^q \,dx
\biggr)^{s/q},
\end{eqnarray*}
which is (\ref{11.10}).

We now combine (\ref{11.14}) and (\ref{11.10}) as follows
\begin{eqnarray*}
&&\biggl(R^{-2}\int_{|x|\le R}|u|^q \,dx
\biggr)^{1/q}\\
&&\hspace*{-6.08pt}\qquad\stackrel{\mbox{\fontsize{8.36pt}{10.36pt}\selectfont{(\ref
{11.13})}}}{\le}
\biggl(R^{-2}\int_{|x|\le R}\bigl|u-(\overline{u}_M)_{\{|x|\le R\}}\bigr|^q
\,dx\biggr)^{1/q}\\
&&\qquad\le
\biggl(R^{-2}\int_{|x|\le R}
\bigl|u_M-(\overline{u}_M)_{\{|x|\le R\}}\bigr|^q \,dx\biggr)^{1/q}\\
&&\qquad\quad{}+\biggl(R^{-2}\int_{\{|x|\le R\}\cap\{|u|>M\}}|u|^q \,dx
\biggr)^{1/q}\\
&&\hspace*{-23.87pt}\qquad\stackrel{\mbox{\fontsize{8.36pt}{10.36pt}\selectfont{(\ref{11.14})
and (\ref{11.10})}}}{\lesssim}
M^{1/2}\biggl(\int_{\Z^2}|f| \,dx\biggr)^{1/2}+M^{1-s/(2q)}\biggl(\int_{\Z^2}|f| \,dx\biggr)^{s/(2q)}\\
&&\qquad\quad{}
+M^{1-s/q}\biggl(R^{-2}\int_{|x|\le R}|u|^q \,dx
\biggr)^{s/q^2}.
\end{eqnarray*}
We claim that this estimate contains the desired estimate.
Indeed, using the abbreviations
\[
U := \biggl(R^{-2}\int_{|x|\le R}|u|^q \,dx\biggr)^{1/q}
\quad\mbox{and}\quad
F := \int_{\Z^2}|f| \,dx,
\]
we rewrite the above as
%
%
\begin{equation}\label{11.15}
U \lesssim M^{1/2} F^{1/2}+M^{1-s/(2q)} F^{s/(2q)}
+M^{1-s/q } U^{s/q}.
\end{equation}
Since $s>q$, choosing $M\sim U$ sufficiently large, we may absorb the last
term of (\ref{11.15}) into the left-hand side. This yields
\[
U \lesssim U^{1/2} F^{1/2}+U^{1-s/(2q)} F^{s/(2q)}.
\]
Using Young's inequality twice in the right-hand side since $s>2q$,
we obtain
as desired $U\lesssim F$, which shows
\[
\biggl(R^{-2}\int_{|x|\le R}|G_T-{\overline{G}_T}_{|x|\leq R}|^q
\,dx
\biggr)^{1/q} \lesssim1.
\]

\textit{Step} 3. Proof of (ii).

We first derive a weak $L^4$-estimate on $\{|x|\leq R\}$:
%
%
\begin{equation}\label{L11:20}
R^{-2}\mathcal{L}_2(\{G_T>M\}\cap\{|x|\leq R\}) \lesssim M^{-4}.
\end{equation}
For that purpose, we combine (\ref{L11.1}), which for $R\sim\sqrt{T}$
turns into
%
%
\begin{equation}\label{L11:21}
R^{-2}\int_{\Z^2} G_{T,M}^2 \,dx+\int_{\Z^2} |\nabla G_{T,M}|^2 \,dx
\lesssim M,
\end{equation}
with the Poincar\'e--Sobolev estimate
\[
\biggl( R^{-2}\int_{|x|\leq R}|G_{T,M}-{\overline{G}_{T,M}}_{\{|x|\leq
R\}
}|^8 \,dx\biggr)^{1/8} \lesssim\biggl(\int_{|x|\leq R}|\nabla
G_{T,M}|^2 \,dx \biggr)^{1/2}
\]
in form of
\begin{eqnarray*}
\biggl( R^{-2}\int_{|x|\leq R}G_{T,M}^8 \,dx\biggr)^{1/8} &\lesssim&
\biggl(\int_{|x|\leq R}|\nabla G_{T,M}|^2 \,dx \biggr)^{1/2}\\
&&{} +
\biggl(R^{-2}\int_{|x|\leq R}G_{T,M}^2 \,dx\biggr)^{1/2}.
\end{eqnarray*}
This yields (\ref{L11:20}):
\[
\bigl(R^{-2}M^8\mathcal{L}_2(\{G_T>M\}\cap\{|x|\leq R\})\bigr)^{1/8}\leq\biggl(
R^{-2}\int_{|x|\leq R}G_{T,M}^8 \,dx\biggr)^{1/8} \lesssim M^{1/2}.
\]
We now argue that (\ref{L11:20}) yields (\ref{L11.19b}). Indeed, combining
\begin{eqnarray*}
R^{-2}\int_{\{G_T>M\}\cap\{|x|\leq R\}}G_{T}^2 \,dx&=&qR^{-2}\int
_M^\infty\mathcal{L}_2(\{G_T>M'\}\cap\{|x|\leq R\}){M'}\,dM'\\
&&{} +R^{-2}M^2\mathcal{L}_2(\{G_T>M\}\cap\{|x|\leq R\})\\
&\stackrel{\mbox{\fontsize{8.36pt}{10.36pt}\selectfont{(\ref
{L11:20})}}}{\lesssim} & M^{-2}
\end{eqnarray*}
with the trivial inequality
\[
R^{-2}\int_{|x|\leq R} G_{T,M}^2 \,dx \lesssim M^2
\]
for $M=1$ yields property (ii) of the lemma.

\textit{Step} 4. Proof of (iii).

We establish for all $q>1$ and $R\gg1$
%
%
\begin{equation}\label{L11.18}
(2R)^{-d}\int_{|x|\geq2R}G_T^q \,dx \lesssim\frac
{T}{R^2}R^{-d}\int
_{R\leq|x|\leq2R}G_T^q\,dx.
\end{equation}
Indeed, we test (\ref{L11.0}) with $\eta^2 G_T^{q-1}$ where the cut-off
function $\eta$ is chosen as follows
\begin{eqnarray*}
\eta&\equiv&1 \qquad\mbox{in } \{|x|\ge2R\},\\
\eta&\equiv&0 \qquad\mbox{in } \{|x|\le R\},\qquad
|\nabla\eta| \lesssim R^{-1},\qquad
0 \le\eta\le1,
\end{eqnarray*}
yielding
\[
T^{-1}\int_{\Z^d} \eta^2G_T^q\,dx+\int_{\Z^d} \nabla(\eta
^2G_T^{q-1})\cdot A\nabla G_T \,dx=0.
\]
Arguing as for (\ref{L11.-4}), this yields
\[
T^{-1}\int_{|x|\geq2R} G_T^q \,dx+\int_{|x|\geq2R} |\nabla
G_T^{q/2}|^2 \,dx \lesssim R^{-2}\int_{R\leq|x|\leq2R}G_T^q \,dx,
\]
so in particular (\ref{L11.18}).

We now turn to (\ref{11.17}). We introduce the abbreviations
\begin{eqnarray*}
R_k&:=&2^k\sqrt{T}, \\
\Lambda_k&:=&R_k^{-d}\int_{R_{k}\leq|x|\leq R_{k+1}}G_T^q\,dx,
\end{eqnarray*}
so that (\ref{L11.18}) turns into
\[
\Lambda_{k+1}\leq C\frac{T}{R^2_k}\Lambda_k=C4^{-k}\Lambda_k,
\]
where $C$ denotes a constant depending only on $\alpha,\beta$, and $d$.
This yields by iteration
\[
\Lambda_k \leq\Lambda_0 C^{k}\prod_{i=0}^{k-1} 4^{-i} =\Lambda_0
C^{k}4^{-{(k-1)k}/{2}}=\Lambda_0C^{k}2^{-(k-1)k}.
\]
Thus, for all $r>0$,
\begin{eqnarray*}
\ln\biggl(\biggl(\frac{R_k}{\sqrt{T}}\biggr)^{r}\frac{\Lambda
_k}{\Lambda
_0} \biggr) &\leq&kr\ln2+(k+1)\ln C-k^2\ln2\\
&\lesssim& 1
\end{eqnarray*}
for $k$ large enough.
Hence,
\[
\int_{R_{k}\leq|x|\leq R_{k+1}}G_T^q\,dx \lesssim\Lambda_0 R_k^d
\biggl(\frac{R_k}{\sqrt{T}}\biggr)^{-r}.
\]

To conclude the proof of (iii), it remains to argue
that
%
%
\begin{equation}\label{L11.19}
\Lambda_0 \lesssim
\cases{
d=2, & 1 ,\cr
d>2, & $\bigl(\sqrt{T}{}^{2-d}\bigr)^q$.}\vadjust{\goodbreak}
\end{equation}
For $d>2$, this a consequence of (\ref{11.16b}), whereas for $d=2$ we
combine (\ref{11.16}) with (\ref{L11.19b}) as follows:
\begin{eqnarray*}
\Lambda_0 &\leq&T^{-1}\int_{|x|\leq2\sqrt{T}}G_T(x)^q\,dx \\
&\leq& T^{-1}\biggl( \biggl( \underbrace{\int_{|x|\leq2\sqrt
{T}}\bigl|G_T(x)-{\overline{G}_T}_{|x|\leq2\sqrt{T}}\bigr|^q\,dx}_{
\stackrel{\mbox{\fontsize{8.36pt}{10.36pt}\selectfont{(\ref
{11.16})}}}{\lesssim}\sqrt{T}{}^2}\biggr)^{1/q}
\\
&&\hspace*{34.9pt}{}+T^{-1}\biggl( \int_{|x|\leq2\sqrt{T}}\bigl(\underbrace{{\overline
{G}_T}_{|x|\leq2\sqrt{T}}}_{
\stackrel{\mbox{\fontsize{8.36pt}{10.36pt}\selectfont{(\ref
{L11.19b})}}}{\lesssim}\sqrt{T}{}^2}\bigr)^q\,dx\biggr)^{1/q}\biggr)^q \\
&\lesssim& 1.
\end{eqnarray*}

\textit{Step} 5. Modifications due to the discreteness.

The only place where we have used the Leibniz rule is the
proof of the Cacciopoli inequality (\ref{L11.8}). At the discrete level,
we have for $i\in\{1,\ldots,d\}$
%
%
\begin{eqnarray}\label{eq:disc-LR}
&&\nabla_i(\eta^2 G_T^{q-1})(x)\nonumber\\
&&\qquad=\eta^2(x+\ee_i)G_T^{q-1}(x+\ee_i)-\eta^2(x)G_T^{q-1}(x)
\nonumber\\[-8pt]\\[-8pt]
&&\qquad= \frac{\eta^2(x+\ee_i)+\eta^2(x)}{2}\bigl(G_T^{q-1}(x+\ee
_i)-G_T^{q-1}(x)\bigr) \nonumber\\
&&\qquad\quad{} + \frac{\eta^2(x+\ee_i)-\eta^2(x)}{2}\bigl(G_T^{q-1}(x+\ee
_i)+G_T^{q-1}(x)\bigr).\nonumber
\end{eqnarray}
Taking advantage of the diagonal structure of $A$ (although this is not
essential), we obtain
\begin{eqnarray*}
\nabla(\eta^2 G_T^{q-1})\cdot A\nabla G_T(x)
&=&\sum_{i=1}^d\nabla_i(\eta^2 G_T^{q-1})(x)a(x,x+\ee_i)\nabla_i
G_T(x) \\
&\stackrel{\mbox{\fontsize{8.36pt}{10.36pt}\selectfont{(\ref
{eq:disc-LR})}}}{=}& \sum_{i=1}^d a(x,x+\ee_i) \frac{\eta^2(x+\ee
_i)+\eta
^2(x)}{2}\\
&&\hspace*{12.13pt}{}\times\underbrace{\bigl(G_T^{q-1}(x+\ee_i)-G_T^{q-1}(x)\bigr)\nabla_i
G_T(x)}_{\geq0} \\
&&{} + \sum_{i=1}^d a(x,x+\ee_i)\frac{\eta^2(x+\ee_i)- \eta
^2(x)}{2}\\
&&\hspace*{25pt}{}\times\bigl(G_T^{q-1}(x+\ee_i)+G_T^{q-1}(x)\bigr)\nabla_i G_T(x).
\end{eqnarray*}
Since the underbraced term is nonnegative, the lower and upper bounds
on $a$ yield
\begin{eqnarray*}
&&\nabla(\eta^2 G_T^{q-1})\cdot A\nabla G_T(x)\\
&&\qquad\geq\alpha\sum_{i=1}^d\frac{\eta^2(x+\ee_i)+\eta
^2(x)}{2}\bigl(G_T^{q-1}(x+\ee_i)-G_T^{q-1}(x)\bigr)\nabla_i G_T(x) \\
&&\qquad\quad{} - \beta\sum_{i=1}^d |\nabla_i \eta(x)|\frac{\eta(x+\ee
_i)+\eta
(x)}{2}\bigl(G_T^{q-1}(x+\ee_i)+G_T^{q-1}(x)\bigr)|\nabla_i G_T(x)| \\
&&\qquad\hspace*{-7.03pt} \stackrel{\mathrm{Young}}{\geq}  \alpha\sum_{i=1}^d\frac{\eta
^2(x+\ee_i)+\eta^2(x)}{2}\bigl(G_T^{q-1}(x+\ee_i)-G_T^{q-1}(x)\bigr)\nabla_i
G_T(x) \\
&&\qquad\quad{} - \beta C\sum_{i=1}^d\bigl(G_T(x+\ee_i)^q+G_T(x)^q\bigr)|\nabla_i \eta
(x)|^2 \\
&&\qquad\quad{} -\beta C^{-1}\sum_{i=1}^d\underbrace{\biggl(\frac{\eta(x+\ee
_i)+\eta
(x)}{2}\biggr)^2}_{\leq({\eta^2(x+\ee_i)+\eta^2(x)})/{2}}
(\nabla_i G_T(x))^2\\
&&\qquad\quad\hspace*{9.6pt}\hspace*{45.9pt}{}\times\bigl(G_T^{q-2}(x+\ee_i)+G_T^{q-2}(x)\bigr).
\end{eqnarray*}
Using the inequality (proved at the end of the step)
%
%
\begin{eqnarray}\label{eq:bcq-1}
&&2(b^{q-1}-c^{q-1})(b-c) \geq(b-c)^2(b^{q-2}+c^{q-2}) \nonumber\\[-8pt]\\[-8pt]
&&\eqntext{\mbox{for }b,c\geq0, q\geq2,}
\end{eqnarray}
we may absorb the last term of the right-hand side of the latter
inequality into
the first term for $C$ large enough, so that it turns into
%
%
\begin{eqnarray}\label{eq:L11a-step4-1}
&&\nabla(\eta^2 G_T^{q-1})\cdot A\nabla G_T(x)\nonumber\\
&&\qquad \geq (\alpha-2\beta C^{-1})\sum_{i=1}^d\frac{\eta^2(x+\ee
_i)+\eta
^2(x)}{2}\nonumber\\[-8pt]\\[-8pt]
&&\qquad\quad\hspace*{74.6pt}{}\times\bigl(G_T^{q-1}(x+\ee_i)-G_T^{q-1}(x)\bigr)\nabla_i G_T(x)\nonumber\\
&&\qquad\quad{} - \beta C\sum_{i=1}^d\bigl(G_T(x+\ee_i)^q+G_T(x)^q\bigr)|\nabla_i \eta
(x)|^2.\nonumber
\end{eqnarray}
Using now the following inequality
%
%
\begin{equation}\label{eq:bcq-2}
(b^{q-1}-c^{q-1})(b-c) \gtrsim(b^{q/2}-c^{q/2})^2 \qquad\mbox{for
}b,c\geq0, q>1,
\end{equation}
(\ref{eq:L11a-step4-1}) finally turns into
\begin{eqnarray*}
&&\nabla(\eta^2 G_T^{q-1})\cdot A\nabla G_T(x)\\
&&\qquad \gtrsim \sum_{i=1}^d\frac{\eta^2(x+\ee_i)+\eta
^2(x)}{2}\bigl(G_T^{q/2}(x+\ee_i)-G_T^{q/2}(x)\bigr)^2 \\
&&\qquad\quad{} - C\sum_{i=1}^d\bigl(G_T(x+\ee_i)^q+G_T(x)^q\bigr)|\nabla_i \eta(x)|^2.
\end{eqnarray*}
Combining this with (\ref{L11.-4}) yields
%
%
\begin{eqnarray}\label{eq:L11-disc-Caccio}
&&\int_{\Z^d}\eta^2(x)|\nabla G_T^{q/2}(x)|^2 \,dx \nonumber\\[-8pt]\\[-8pt]
&&\qquad\lesssim\int
_{\Z^d}
\bigl(G_T(x+\ee_i)^q+G_T(x)^q\bigr)|\nabla_i \eta(x)|^2
\,dx,\nonumber
\end{eqnarray}
which implies as desired
\[
\int_{2R\leq|x|<4R}|\nabla G_T^{q/2}(x)|^2 \,dx \lesssim R^{-2}\int
_{R\leq|x|< 8R} G_T(x)^q \,dx,
\]
provided that $\eta$ satisfies in addition
\[
\eta(x)=0 \qquad\mbox{for }x \notin\{y\dvtx R+1 \leq|y|\leq8R-1\},
\]
which is no restriction since $R\gg1$.

We quickly sketch the proofs of (\ref{eq:bcq-1}) and (\ref
{eq:bcq-2}) to conclude.
Inequality (\ref{eq:bcq-1}) follows by symmetry from
\begin{eqnarray*}
&&(b^{q-1}-c^{q-1})(b-c) - (b-c)^2c^{q-2}\\
&&\qquad=(b-c)(b^{q-1}-bc^{q-2})
\\
&&\qquad=b(b-c)(b^{q-2}-c^{q-2})\\
&&\qquad=b|b-c||b^{q-2}-c^{q-2}| \geq0.
\end{eqnarray*}
To prove (\ref{eq:bcq-2}) we first note that by homogeneity and
nonnegativity of $b$ and $c$, it is enough
to consider $c=1$ and $b\geq0$.
We introduce the function $h=\R^+\to\R^+$ defined by
\[
h(b) = \cases{
\dfrac{(b^{q/2}-1)^2}{(b^{q-1}-1)(b-1)}, &\quad $b\neq1$, \cr
\dfrac{q^2}{4(q-1)}, &\quad $b=1$.}
\]
Since $h\geq0$, the claim is proved if $h$ is bounded on $\R^+$.
As $h(0)=1$ and $\lim_{b\to\infty}h(b)=1$, it is enough to prove that
$h$ is continuous on $\R^+$.
A Taylor expansion
around $b=1$ yields
\begin{eqnarray*}
(b^{q/2}-1)^2&=& \frac{q^2}{4}(b-1)^2+o\bigl((b-1)^2\bigr), \\
(b^{q-1}-1)(b-1)&=& (q-1)(b-1)^2+o\bigl((b-1)^2\bigr).
\end{eqnarray*}
Hence, $\lim_{b\to1}h(b)=h(1)$, $h$ is continuous and therefore
bounded on $\R^+$, as desired.

\subsection{\texorpdfstring{Proof of Lemma \protect\ref{lem:int-grad}}{Proof of Lemma 2.9.}}

The proof relies on three ingredients: a Meyers' estimate
based on the $L^q$ theory
for the constant-coefficients Helmholtz projection, a Cacciopoli
estimate and the estimates of Lemma \ref{L11a}.

We begin with Meyers' estimates.
Let $u\dvtx\Z^d\to\R$, $f\dvtx\Z^d\to\R$, and $g\dvtx\Z^d\to\R^d$ have
support in $\{|x|<R\}$, and let $u$ satisfy the equation
%
%
\begin{equation}\label{eq:intgrad-step1-1}
-\nabla^*\cdot A(x)\nabla u(x)=\nabla^*\cdot g(x)+f(x) \qquad\mbox{in
}\Z^d.
\end{equation}
We claim that there exists $p>2$ depending only on $\alpha,\beta$, and
$d$ such that for all $R\gg1$,
the following $L^p$-estimate holds
%
%
\begin{eqnarray} \label{eq:Meyers}
\biggl(\int_{\Z^d}|\nabla u(x)|^p \,dx \biggr)^{1/p}&\lesssim&
\biggl(\int
_{\Z^d}|g(x)|^p \,dx \biggr)^{1/p} \nonumber\\
&&{} + R^{1-d({1}/{2}-{1}/{p})}\biggl(\int_{\Z^d}|f(x)|^2 \,dx
\biggr)^{1/2}.
\end{eqnarray}
As in the original paper \cite{Meyers-63} by Meyers, the proof of
(\ref
{eq:Meyers}) relies on a perturbation argument and on the $L^q$
regularity theory for the
Helmholtz projection.

\textit{Step} 1. $L^q$ regularity for the Helmholtz projection.

Let $\mathcal{H}\dvtx L^2(\Z^d,\R^d)\to L^2(\Z^d,\R^d)$ denote the Helmholtz
projection, that is,
the orthogonal projection onto gradient fields for the inner product of
$L^2(\Z^d,\R^d)$.
By definition, $\mathcal H$ is continuous on $L^2(\Z^d,\R^d)$ and satisfies
%
%
\begin{equation}\label{eq:Helm-1}
\|\mathcal H g\|_{L^2(\Z^d,\R^d)} \leq\| g\|_{L^2(\Z^d,\R^d)}.
\end{equation}
Let us show that $\mathcal H$ can be extended to a continuous operator from
$L^q(\Z^d,\R^d)$ to $L^q(\Z^d,\R^d)$ for all $1<q<\infty$.
The proof is standard, appealing to Calder\'on--Zygmund singular
integral theory and
to Marcinkiewicz interpolation theorem (such theorems apply to the
discrete case under investigation
since the associated measure has the so-called ``doubling'' property).
Since $\mathcal H$ commutes with translations,
it is a convolution operator: there exists a matrix-valued kernel $K$
such that
%
%
\begin{equation}\label{eq:Helm-2}
\mathcal H g(x) = \int_{\Z^d}K(x-y)g(y)\,dy.
\end{equation}
From an elementary Fourier series analysis (see \cite{Martinsson-02}
for related arguments),
we infer that the symbol of $K$ coincides with the symbol of the second
derivative
of the Green's function of the Laplace equation studied in \cite
{Martinsson-02}.
In particular, from the analysis of \cite{Martinsson-02}, we learn that
%
%
\begin{equation}\label{eq:Helm-3}
|\nabla K(x)| \lesssim\frac{1}{1+|x|^{d+1}}.
\end{equation}
We are therefore in position to apply Calder\'on--Zygmund's theory (see
\cite{Stein-93}, Theorem 2, page 17), which shows that
$\mathcal H$ is of weak type $(1,1)$ (see the proof of \cite{Stein-93},
Theorem 3,
page 19).
Appealing to Marcinkiewicz' interpolation theorem
(see \cite{Bergh-Lofstrom-76}, Theorem 1.3.1, page 9)
then shows that $\mathcal H$ can be extended to a continuous operator from
$L^q(\Z^d,\R^d)$ to $L^q(\Z^d,\R^d)$ for all $1<q< 2$. A standard
duality argument
(see, \cite{Stein-70}, 2.5(c), page 33, e.g.)
implies that $\mathcal H$ can also be extended to a continuous operator from
$L^q(\Z^d,\R^d)$ to $L^q(\Z^d,\R^d)$ for all $2<q<\infty$.
Let $r>2$ be fixed, and for all $q>1$ let denote by $C_q$
the norm of $\mathcal H$ in $\mathcal L (L^q(\Z^d,\R^d),L^q(\Z^d,\R^d))$.
Then Riesz--Thorin interpolation theorem (see \cite{Bergh-Lofstrom-76},
Theorem 1.1.1,
page 2) shows
that for all $\theta\in(0,1)$, $C_{2\theta+r(1-\theta)}\leq
C_2^\theta C_r^{1-\theta}$, so that
%
%
\begin{equation}\label{eq:Helm-4}
\limsup_{q\to2}C_q \leq1
\end{equation}
since $C_2\leq1$ by (\ref{eq:Helm-1}).

We now turn to the proof of (\ref{eq:Meyers}) proper and
proceed with the perturbation argument.

\textit{Step} 2. Proof of (\ref{eq:Meyers}) for $f\equiv0$.

We first assume that $f\equiv0$, and rewrite the left-hand side of
(\ref{eq:intgrad-step1-1})
as a perturbation of the operator $-\frac{\alpha+\beta
}{2}\bigtriangleup$:
\[
-\frac{\alpha+\beta}{2}\bigtriangleup u = \nabla^*\cdot
\biggl(g+
\biggl(A-\frac{\alpha+\beta}{2}\Id\biggr)\nabla u \biggr)
\]
or equivalently in the form
%
%
\begin{equation}\label{eq:Helm-7}
-\bigtriangleup u = \nabla^*\cdot\biggl(\frac{2}{\alpha+\beta}
\biggl(g+\biggl(A-\frac{\alpha+\beta}{2}\Id\biggr)\nabla u \biggr)\biggr).
\end{equation}
In order to apply the $L^q$ theory for the Helmholtz projection, we
need to show
that
%
%
\begin{equation}\label{eq:Helm-5}
\nabla u \equiv\mathcal H \biggl( \frac{2}{\alpha+\beta}
\biggl(g+\biggl(A-\frac{\alpha+\beta}{2}\Id\biggr)\nabla u \biggr)\biggr).
\end{equation}
Since $\nabla u$ is obviously a gradient, it remains to show that for
all $\zeta\dvtx\Z^d\to\R$ such that $\nabla\zeta\in L^2(\Z^d,\R^d)$
one has
%
%
\begin{eqnarray}\label{eq:Helm-6}
&&\int_{\Z^d}\nabla u (x)\cdot\nabla\zeta(x)\,dx \nonumber\\[-8pt]\\[-8pt]
&&\qquad= \int_{\Z^d}
\biggl( \frac{2}{\alpha+\beta}\biggl(g+\biggl(A-\frac{\alpha+\beta}{2}\Id
\biggr)\nabla u
\biggr)\biggr)(x) \cdot\nabla\zeta(x)\,dx.\nonumber
\end{eqnarray}
To this aim, we multiply (\ref{eq:Helm-7}) by $\zeta$ and integrate by
parts using that $u$, $\nabla u$ and $g$ have compact supports.
This yields (\ref{eq:Helm-6})
and proves therefore (\ref{eq:Helm-5}).
The continuity of $\mathcal H$ from $L^q(\Z^d,\R^d)$ to $L^q(\Z^d,\R
^d)$ proved in Step 1 then implies that
%
%
\begin{eqnarray}\label{eq:intgrad-step1-4}\qquad
&&\biggl(\int_{|x|\leq R}|\nabla u(x)|^q \,dx
\biggr)^{1/q}\nonumber\\[-8pt]\\[-8pt]
&&\qquad\leq C_q \frac{2}{\alpha+\beta}\biggl( \int_{|x|\leq R}
\biggl|g(x)+\biggl(A(x)-\frac{\alpha+\beta}{2}\Id\biggr)\nabla
u(x)\biggr|^q
\,dx \biggr)^{1/q}.\nonumber
\end{eqnarray}
Using the triangle inequality, (\ref{eq:intgrad-step1-4}) turns into
%
%
\begin{eqnarray}\label{eq:intgrad-step1-5}
&&\biggl(\int_{|x|\leq R}|\nabla u(x)|^q \,dx \biggr)^{1/q}
\nonumber\\
&&\qquad\leq C_q\frac{2}{\alpha+\beta}
\biggl(\int_{|x|\leq R}|g(x)|^q \,dx \biggr)^{1/q}\\
&&\qquad\quad{} +C_q\frac{2}{\alpha+\beta}\biggl(\int_{|x|\leq R}\biggl|
\biggl(A(x)-\frac{\alpha+\beta}{2}\Id\biggr)\nabla u(x)\biggr|^q
\,dx
\biggr)^{1/q} .\nonumber
\end{eqnarray}
Since $a\in\mathcal{A}_{\alpha\beta}$, $|(A(x)-\frac
{\alpha
+\beta}{2}\Id)\nabla u(x)|\leq\frac{\beta-\alpha
}{2}|\nabla u(x)|$ and we may absorb
the term
\begin{eqnarray*}
&&C_q\frac{2}{\alpha+\beta}\biggl(\int_{|x|\leq R}
\biggl|
\biggl(A(x)-\frac{\alpha+\beta}{2}\Id\biggr)\nabla u(x)\biggr|^q
\,dx
\biggr)^{1/q}\\
&&\qquad \leq C_q\frac{\beta-\alpha}{\alpha+\beta}\biggl(\int_{|x|\leq
R}|\nabla u(x)|^q \,dx\biggr)^{1/q}
\end{eqnarray*}
into the left-hand side of (\ref{eq:intgrad-step1-5}) provided that
%
%
\begin{equation}\label{eq:intgrad-step1-6}
C_q\underbrace{\frac{\beta-\alpha}{\alpha+\beta}}_{<1}<1.
\end{equation}
The interpolation property (\ref{eq:Helm-4}) ensures there exists $p>2$
such that (\ref{eq:intgrad-step1-6}) holds for all $p\geq q\geq2$.
For such a $q$, we then have
%
%
\begin{equation}\label{eq:meyers12}
\biggl(\int_{|x|\leq R}|\nabla u(x)|^q \,dx \biggr)^{1/q}
\lesssim\biggl(\int_{|x|\leq R}|g(x)|^q \,dx \biggr)^{1/q} ,
\end{equation}
as desired.

\textit{Step} 3. Proof of (\ref{eq:Meyers}) for general $f$.

Note that since $u$ and $g$ have compact supports, equation (\ref
{eq:intgrad-step1-1}) implies
that $\int_{\Z^d}f(x)\,dx=0$.
We first show that there exists $\nabla w\in L^2(\Z^d,\R^d)$
such that for all $\zeta\dvtx\Z^d\to\R$ with $\nabla\zeta\in L^2(\Z
^d,\R
^d)$, one has
%
%
\begin{equation}\label{eq:Helm-8}
\int_{\Z^d}\nabla w(x)\cdot\nabla\zeta(x)\,dx = \int_{\Z
^d}f(x)\zeta(x)\,dx,
\end{equation}
so that (\ref{eq:Helm-5}) turns into
\[
\nabla u \equiv\mathcal H \biggl( \frac{2}{\alpha+\beta}
\biggl(g+\biggl(A-\frac{\alpha+\beta}{2}\Id\biggr)\nabla u \biggr)+\nabla w\biggr).
\]
Provided
%
%
\begin{equation}\label{eq:Meyers-f-0}
\biggl(\int_{\Z^d}|\nabla w(x)|^q\,dx\biggr)^{1/q} \lesssim
R^{1-d({1}/{2}-{1}/{q})}\biggl(\int_{\Z^d}f(x)^2\,dx\biggr)^{1/2},
\end{equation}
for all $2\leq q \leq\tilde q$ for some $\tilde q>2$, we then conclude
as in the case $f\equiv0$ (with potentially a smaller $p$).
To prove the existence of such a $\nabla w$, we proceed by minimization
and consider the problem
%
%
\begin{eqnarray}\label{eq:Helm-9}
&&\inf\biggl\{ \int_{\Z^d}|\nabla\zeta(x)|^2\,dx-\int_{\Z
^d}f(x)\zeta
(x)\,dx ;\nonumber\\[-8pt]\\[-8pt]
&&\qquad\hspace*{23.6pt} \zeta\dvtx\Z^d\to\R, \nabla\zeta\in L^2(\Z^d,\R^d)
\biggr\}.\nonumber
\end{eqnarray}
The same argument as in the proof of Riesz' theorem yields the
existence of a minimizer once
one shows that the functional is coercive.
Let $R$ be large enough so that $f$ has support in $\{|x|<R\}$, and
denote by $\bar\zeta_{\{|x|<R\}}$ the average of $\zeta$
on $\{|x|<R\}$. Since $f$ has zero average, one may subtract the
average of $\zeta$ and
obtain by Cauchy--Schwarz and Poincar\'e's inequalities
\begin{eqnarray*}
\biggl|\int_{\Z^d}f(x)\zeta(x)\,dx\biggr|&=&\biggl|\int
_{|x|<R}f(x)\bigl(\zeta
(x)- \bar\zeta_{\{|x|<R\}}\bigr)\,dx\biggr| \\
&\lesssim& R \biggl(\int_{\Z^d}f(x)^2\,dx\biggr)^{1/2}\biggl(\int
_{|x|<R}|\nabla\zeta(x)|^2\,dx\biggr)^{1/2} \\
&\lesssim& -2R^2 \int_{\Z^d}f(x)^2\,dx+ \frac{1}{2}\int_{\Z
^d}|\nabla
\zeta(x)|^2\,dx.
\end{eqnarray*}
This shows that for all test functions $\zeta$
%
%
\begin{eqnarray}\label{eq:bdded-below}
&&
\int_{\Z^d}|\nabla\zeta(x)|^2\,dx-\int_{\Z^d}f(x)\zeta(x)\,dx \nonumber\\[-8pt]\\[-8pt]
&&\qquad\geq -2R^2
\int_{\Z^d}f(x)^2\,dx+\frac{1}{2}\int_{\Z^d}|\nabla\zeta(x)|^2\,dx,\nonumber
\end{eqnarray}
as desired. This proves the existence of a minimizer $w\dvtx\Z^d\to\R$
such that $\nabla w\in L^2(\R^d,\Z^d)$.
In addition, it satisfies the estimate
%
%
\begin{equation}\label{eq:Helm-10}
\int_{\Z^d}|\nabla w(x)|^2\,dx \leq 4R^2 \int_{\Z^d}f(x)^2\,dx.
\end{equation}
Since $w$ is a minimizer of (\ref{eq:Helm-9}), the first variation of
the energy at $w$ vanishes,
and $w$ satisfies (\ref{eq:Helm-8}).

It remains to estimate the $L^q$ norm of $\nabla w$ for some $q>2$.
To this aim, we argue that
%
%
\begin{eqnarray}\label{eq:Helm-11}
\int_{\Z^d} |\nabla\nabla w(x)|^2\,dx &=& \int_{\Z^d}(\triangle
w(x))^2\,dx \nonumber\\[-8pt]\\[-8pt]
&=& \int_{\Z^d}f(x)^2\,dx.\nonumber
\end{eqnarray}
As in the continuum case, the first identity in (\ref{eq:Helm-11})
follows directly from two integrations
by parts for $w$ with compact support. For general $w$, the boundary
term involves products of
first and second derivatives of $w$ on spheres of large radius~$R$. In
our discrete setting,
these boundary terms can be estimated by the integral of $|\nabla w|^2$
outside the ball of radius $R$, which is
finite since $|\nabla w|^2$ is integrable by construction.
Hence, the boundary terms can be made to vanish in the limit $R\to
\infty$.
The second identity in (\ref{eq:Helm-11}) follows from the fact that
$w$ solves the equation
\[
-\triangle w(x) = f(x) \qquad\mbox{in }\Z^d,
\]
which is a consequence of (\ref{eq:Helm-8}).
We are in position to conclude.
For $d>2$, we appeal to Poincar\'e--Sobolev inequality on $\nabla w$
to turn (\ref{eq:Helm-11}) into
%
%
\begin{equation}\label{eq:Meyers-f-3}
\int_{\Z^d}|\nabla w(x)|^{{2d}/({d-2})}\,dx \lesssim\int_{\Z
^d}f(x)^2\,dx.
\end{equation}
Combined with (\ref{eq:Helm-10}), (\ref{eq:Meyers-f-3})
implies (\ref{eq:Meyers-f-0}) for all $2\leq q \leq\frac{2d}{d-2}$ by
H\"older's inequality.
For $d=2$, we appeal to Poincar\'e--Sobolev inequality on $(\nabla_i
w)^2$ for $i\in\{1,\ldots,d\}$ to
turn (\ref{eq:Helm-10}) and (\ref{eq:Helm-11}) into
%
%
\begin{eqnarray} \label{eq:Meyers-f-4}
\int_{\Z^2}(\nabla_i w(x))^{4}\,dx &\lesssim& \biggl( \int_{\Z
^2}|\nabla_i
(\nabla_i w(x))^2|\,dx\biggr)^2 \nonumber\\
&\lesssim& \int_{\Z^2}|\nabla\nabla w(x)|^2\,dx\int_{\Z^2} |\nabla
w(x)|^2\,dx \\
&\lesssim& R^2 \biggl(\int_{\Z^2}f(x)^2\,dx\biggr)^2.\nonumber
\end{eqnarray}
Combined with (\ref{eq:Helm-10}), (\ref{eq:Meyers-f-4})
implies (\ref{eq:Meyers-f-0}) for all $2\leq q\leq4$ by H\"older's inequality.

\textit{Step} 4. Cacciopoli estimate.

We need the following finer version of (\ref{L11.8}):
For all $\kappa\in\R$,
%
%
\begin{eqnarray}\label{eq:int-grad-step2-1}
&&\int_{2R\leq|x|\leq16R}|\nabla G_T(x)|^2 \,dx\nonumber\\
&&\qquad\lesssim R^{-2}\int
_{R\leq|x|\leq32R} \bigl(G_T(x)-\kappa\bigr)^2 \,dx \\
&&\qquad\quad{} +T^{-1}|\kappa|\int_{R\leq|x|\leq32R} |G_T(x)-\kappa|
\,dx.\nonumber
\end{eqnarray}
This variant of Cacciopoli's estimate can be proved
along the lines of (\ref{L11.8}), multiplying the equation by $\eta
^2(G_T(x)-\kappa)$ instead of $\eta^2G_T(x)$.
The zero order term then brings the new term in the
right-hand side of (\ref{eq:int-grad-step2-1}).
By Young and Cauchy--Schwarz' inequalities, the second term of the
right-hand side is controlled by
\begin{eqnarray*}
&&
T^{-1}|\kappa|\int_{R\leq|x|\leq32R} |G_T(x)-\kappa| \,dx \\
&&\qquad =
T^{-1}|\kappa|R \int_{R\leq|x|\leq32R} R^{-1}|G_T(x)-\kappa| \,dx \\
&&\qquad\lesssim T^{-2}|\kappa|^2R^2+R^{-2}\int_{R\leq|x|\leq32R}
|G_T(x)-\kappa|^2 \,dx .
\end{eqnarray*}
Hence, it only remains to estimate the first term of the right-hand
side of
(\ref{eq:int-grad-step2-1}).
To this aim, we appeal to H\"older's inequality with exponents
$(p/2,p/(p-2))$ for $p\geq2$:
\begin{eqnarray*}
&&\int_{R\leq|x|\leq32R} |G_T(x)-\kappa|^2 \,dx\\
&&\qquad \leq \biggl(\int
_{R\leq|x|\leq32R} |G_T(x)-\kappa|^p \,dx
\biggr)^{2/p}(R^d)^{(p-2)/p}\\
&&\qquad=R^{d-2d/p}\biggl(\int_{R\leq|x|\leq32R} |G_T(x)-\kappa|^p
\,dx\biggr)^{2/p}.
\end{eqnarray*}
Using these last two estimates and the elementary inequality
$(a^2+b^2)^{1/2}\lesssim a+b$,
(\ref{eq:int-grad-step2-1}) turns into
%
%
\begin{eqnarray}\label{eq:int-grad-step2-1-d=2}
&&\biggl(\int_{2R\leq|x|\leq16R}|\nabla G_T(x)|^2 \,dx
\biggr)^{1/2}\nonumber\\
&&\qquad\lesssim R^{-1}R^{d/2-d/p}\biggl(\int_{R\leq|x|\leq32R}
|G_T(x)-\kappa|^p \,dx\biggr)^{1/p}\\
&&\qquad\quad{}+T^{-1}R^{d/2+1}|\kappa|,\nonumber
\end{eqnarray}
that we will use with $\kappa={\overline{G}_T}_{\{R\leq|x|\leq32R\}}$.

\textit{Step} 5. In this step, we use Steps 1 and 4 to argue that
%
%
\begin{eqnarray}\label{eq:int-grad-step3-4b}
&&\biggl(\int_{4R<|x|\leq8R}|\nabla G_T(x)|^p \,dx
\biggr)^{1/p}\nonumber\\
&&\qquad\lesssim (R^{-1}+RT^{-1})\biggl(\int_{2R\leq|x|\leq
16R}\bigl|G_T(x)-{\overline{G}_T}_{\{R\leq|x|\leq32R\}}\bigr|^p \,dx \biggr)^{1/p}
\hspace*{-25pt}\\
&&\qquad\quad
{}+T^{-1}R^{d/p+1}{\overline{G}_T}_{\{R\leq|x|\leq32R\}}.\nonumber
\end{eqnarray}
We apply Meyers' estimate (\ref{eq:Meyers}) to the function
$u=\eta(G_T-{\overline{G}_T}_{\{R\leq|x|\leq32R\}})$, where the
cut-off function $\eta\dvtx\Z^d\to[0,1]$ is such that
%
%
\begin{eqnarray}\label{eq:int-grad-step3-1}
\eta(x)&=&1 \qquad\mbox{for }4R\leq|x|\leq8R,\nonumber\\[-8.5pt]\\[-8.5pt]
\eta(x)&=&0 \qquad\mbox{for
}\cases{
|x|\leq2R+1, \cr
|x|\geq16R-1,}\qquad
|\nabla\eta|\lesssim R^{-1}.\nonumber
\end{eqnarray}
For all $i\in\{1,\ldots,d\}$, the discrete Leibniz rule yields
\[
\nabla_i u(x)=\eta(x)\nabla_i G_T(x) +\bigl(G_T(x+\ee_i)-{\overline
{G}_T}_{\{
R\leq|x|\leq32R\}}\bigr)\nabla_i \eta(x).
\]
Based on this, a direct calculation shows
\begin{eqnarray*}
&&-\nabla^* \cdot A\nabla u(x)\\[-0.5pt]
&&\qquad= -\underbrace{\eta(x) \nabla^*\cdot A\nabla G_T(x)}_{
\stackrel
{\mbox{\fontsize{8.36pt}{10.36pt}\selectfont{(\ref
{eq:int-grad-step3-1}) and (\ref{eq:disc-Green})}}}{=}\eta(x)
T^{-1}G_T(x)}{}-{}\sum_{i=1}^d\nabla_i^* \eta(x)a(x-\ee_i,x)\nabla_i^*
G_T(x)\\[-0.5pt]
&&\qquad\quad{} -\sum_{i=1}^d\nabla_i^*\bigl(\bigl(G_T(x+\ee_i)-{\overline{G}_T}_{\{
R\leq
|x|\leq32R\}}\bigr)a(x,x+\ee_i)\nabla_i \eta(x)\bigr) \\[-0.5pt]
&&\qquad= \nabla^*\cdot\Biggl(-\sum_{i=1}^d\bigl(G_T(x+\ee_i)-{\overline
{G}_T}_{\{
R\leq|x|\leq32R\}}\bigr)a(x,x+\ee_i)\nabla_i \eta(x)\ee_i\Biggr)\\[-0.5pt]
&&\qquad\quad{} -\sum_{i=1}^d\nabla^*_i \eta(x) a(x-\ee_i,x) \nabla_i^*
G_T(x)-\eta
(x) T^{-1}G_T(x).
\end{eqnarray*}
This identity has the form of (\ref{eq:intgrad-step1-1}) provided we
define the functions $f$ and $g$ by
\begin{eqnarray*}
f(x)&=& -\sum_{i=1}^d\nabla^*_i \eta(x) a(x-\ee_i,x) \nabla_i^*
G_T(x)-\eta(x) T^{-1}G_T(x),\\[-0.5pt]
g(x)&=& -\sum_{i=1}^d\bigl(G_T(x+\ee_i)-{\overline{G}_T}_{\{R\leq|x|\leq
32R\}}\bigr)a(x,x+\ee_i)\nabla_i \eta(x)\ee_i.
\end{eqnarray*}
Since $u$, $f$ and $g$ have support in $\{|x|\leq16R\}$, we may apply
estimate (\ref{eq:Meyers}) which yields
\begin{eqnarray*}
&&\biggl(\int_{|x|\leq16R}|\nabla u(x)|^p \,dx \biggr)^{1/p}
\\[-0.5pt]
&&\qquad\lesssim\Biggl(\sum_{i=1}^d\int_{\Z^d}|\nabla_i\eta
(x)|^p\bigl|G_T(x+\ee
_i)-{\overline{G}_T}_{\{R\leq|x|\leq32R\}}\bigr|^p \,dx \Biggr)^{1/p}
\\[-0.5pt]
&&\qquad\quad{} +R^{1-d({1}/{2}-{1}/{p})}\\[-0.5pt]
&&\qquad\quad\hspace*{10.8pt}{}\times\biggl(\int_{\Z^d}\bigl(|\nabla^* \eta
(x)|^2|\nabla^* G_T(x)|^2+T^{-2}\eta(x)^2 G_T(x)^2\bigr) \,dx
\biggr)^{1/2}.
\end{eqnarray*}
Using the property (\ref{eq:int-grad-step3-1}) of $\eta$, and the
triangle inequality, we are left with
%
%
\begin{eqnarray}\label{eq:int-grad-step3-4bb}
&&\biggl(\int_{|x|\leq16R}|\nabla u(x)|^p \,dx \biggr)^{1/p}\nonumber\\
&&\qquad\lesssim R^{-1}\biggl(\int_{2R\leq|x|\leq16R}\bigl|G_T(x)-{\overline
{G}_T}_{\{R\leq|x|\leq32R\}}\bigr|^p \,dx \biggr)^{1/p} \nonumber\\[-8pt]\\[-8pt]
&&\qquad\quad{} +R^{{d}/{p}-d/2}\biggl(\int_{2R\leq|x|\leq16R}|\nabla G_T(x)|^2 \,dx
\biggr)^{1/2} \nonumber\\
&&\qquad\quad{} +R^{{d}/{p}-d/2+1}\biggl(\int_{2R\leq|x|\leq16R}T^{-2}G_T(x)^2 \,dx
\biggr)^{1/2}.\nonumber
\end{eqnarray}
Let us rearrange the terms.
For the third term, the triangle inequality and H\"older's inequality
with exponents $(p/2,p/(p-2))$ show that
\begin{eqnarray*}
&&\biggl( \int_{2R\leq|x|\leq16R}G_T(x)^2 \,dx \biggr)^{1/2}\\
&&\qquad\lesssim
R^{d/2}{\overline{G}_T}_{\{R\leq|x|\leq32R\}}\\
&&\qquad\quad{}
+R^{d/2-d/p}\biggl( \int_{2R\leq|x|\leq16R}\bigl|G_T(x)-{\overline
{G}_T}_{\{
R\leq|x|\leq32R\}}\bigr|^p \,dx \biggr)^{1/p} ,
\end{eqnarray*}
whereas for the second term we appeal to the Cacciopoli estimate (\ref
{eq:int-grad-step2-1-d=2}) with $\kappa={\overline{G}_T}_{\{R\leq
|x|\leq32R\}}$.
Hence, (\ref{eq:int-grad-step3-4bb}) finally turns into
%
%
\begin{eqnarray}\label{eq:int-grad-step3-4bbb}\qquad
&&\biggl(\int_{|x|\leq16R}|\nabla u(x)|^p \,dx
\biggr)^{1/p}\nonumber\\
&&\qquad\lesssim (R^{-1}+RT^{-1})\biggl(\int_{2R\leq|x|\leq
16R}\bigl|G_T(x)-{\overline{G}_T}_{\{R\leq|x|\leq32R\}}\bigr|^p \,dx \biggr)^{1/p}
\\
&&\qquad\quad{}
+T^{-1}R^{d/p+1}{\overline{G}_T}_{\{R\leq|x|\leq32R\}}.\nonumber
\end{eqnarray}
We are in position to conclude the proof of this step.
For all $i\in\{1,\ldots,d\}$, the discrete Leibniz rule yields $\nabla_i
u(x)=\eta(x)\nabla_i G_T(x)+G_T(x+\ee_i)\nabla_i \eta(x)$.
Hence, (\ref{eq:int-grad-step3-1}) implies
that $\nabla u(x)=\nabla G_T(x)$ for $4R\leq|x|\leq8R$, so that
(\ref{eq:int-grad-step3-4bbb}) yields (\ref{eq:int-grad-step3-4b}).

\textit{Step} 6. Proof of (\ref{eq:int-grad}).

We claim that (\ref{eq:int-grad}) follows from (\ref
{eq:int-grad-step3-4b})
and the estimates of Lemma \ref{L11a}.

We distinguish two regimes: $R\leq\sqrt{T}$ and $R\geq
\sqrt{T}$.
We begin with $R\leq\sqrt{T}$.
For the first term of the right-hand side of (\ref
{eq:int-grad-step3-4b}), we
appeal to the BMO estimate (\ref{11.16}) of Lemma \ref{L11a} for $d=2$
and to\vadjust{\goodbreak} the decay estimate (\ref{11.16b}) with ``$q=p$'' for $d>2$, so that
%
%
\begin{eqnarray}\label{eq:int-grad-step4-1.1v}
&&(R^{-1}+RT^{-1})\biggl(\int_{2R\leq|x|\leq
16R}\bigl|G_T(x)-{\overline{G}_T}_{\{R\leq|x|\leq32R\}}\bigr|^p \,dx
\biggr)^{1/p} \nonumber\\[-8pt]\\[-8pt]
&&\qquad\lesssim R^{-1}\bigl(R^{d}R^{(2-d)p}\bigr)^{1/p}=R^{d/p-d+1}.\nonumber
\end{eqnarray}
For the second term, we estimate the average using (\ref{L11.19b}) for $d=2$
\[
{\overline{G}_T}_{\{R\leq|x|\leq32R\}} \lesssim R^{-2}\sqrt{T}{}^2
{\overline{G}_T}_{\{|x|\leq32\sqrt{T}\}} \stackrel{\mbox
{\fontsize
{8.36pt}{10.36pt}\selectfont{(\ref{L11.19b})}}}{\lesssim} R^{-2}T ,
\]
and using (\ref{11.16b}) with ``$q=1$'' for $d>2$
\[
{\overline{G}_T}_{\{R\leq|x|\leq32R\}} \stackrel{\mbox{\fontsize
{8.36pt}{10.36pt}\selectfont{(\ref{11.16b})}}}{\lesssim} R^{2-d}
\lesssim R^{-d}T,
\]
since $R \leq\sqrt{T}$.
Hence, in both cases,
%
%
\begin{equation}\label{eq:int-grad-step4-1.1vv}
{T^{-1}R^{d/p+1}{\overline{G}_T}_{\{R\leq|x|\leq32R\}}} \lesssim
R^{d/p-d+1} .
\end{equation}
From (\ref{eq:int-grad-step4-1.1v}) and (\ref
{eq:int-grad-step4-1.1vv}), we then conclude that
(\ref{eq:int-grad})
holds for $R\leq\sqrt{T}$.

We now deal with the case $R\geq\sqrt{T}$.
For the first term of the right-hand side of (\ref
{eq:int-grad-step3-4b}), we
use the decay estimate (\ref{11.17}) with exponents ``$q=p,r=k+2p$,''
which yields
\begin{eqnarray*}
{\int_{2R\leq|x|\leq16R}\bigl|G_T(x)-{\overline{G}_T}_{\{R\leq|x|\leq
32R\}}\bigr|^p \,dx }
&\lesssim&\int_{R\leq|x|\leq32R}G_T(x)^p \,dx \\
&\stackrel{\mbox{\fontsize{8.36pt}{10.36pt}\selectfont{(\ref
{11.17})}}}{\lesssim} & R^dR^{(2-d)p}\bigl(\sqrt{T}R^{-1}\bigr)^{k+2p},
\end{eqnarray*}
and therefore
%
%
\begin{eqnarray}\label{eq:int-grad-step3-41}
&&(R^{-1}+RT^{-1})\biggl(\int_{2R\leq|x|\leq16R}\bigl|G_T(x)-{\overline
{G}_T}_{\{R\leq|x|\leq32R\}}\bigr|^p \,dx \biggr)^{1/p} \nonumber\\[-8pt]\\[-8pt]
&&\qquad\lesssim
R^{d/p-d+1} \bigl(\sqrt{T}R^{-1}\bigr)^{k/p}.\nonumber
\end{eqnarray}
For the second term, we proceed the same way, and appeal to (\ref
{11.17}) with exponents ``$q=1,r=k/p+2$,'' which yields
\[
{\overline{G}_T}_{\{R\leq|x|\leq32R\}} \lesssim R^{2-d}\bigl(\sqrt
{T}R^{-1}\bigr)^{k/p+2} = TR^{-d}\bigl(\sqrt{T}R^{-1}\bigr)^{k/p},
\]
and therefore
%
%
\begin{equation}\label{eq:int-grad-step3-42}
T^{-1}R^{d/p+1}{\overline{G}_T}_{\{R\leq|x|\leq32R\}}
\lesssim
R^{{d}/{p}-d+1}\bigl(\sqrt{T}R^{-1}\bigr)^{k/p}.
\end{equation}
From (\ref{eq:int-grad-step3-41}) and (\ref{eq:int-grad-step3-42}), we
then deduce that
(\ref{eq:int-grad})
holds for $R\geq\sqrt{T}$ as well.

\subsection{\texorpdfstring{Proof of Corollaries \protect\ref{coro:unif-bound}
and \protect\ref{coro:unif-bound-grad}}{Proof of Corollaries 2.2 and 2.3.}}

These results are easy consequences of Lemmas \ref{L11a}
and \ref{lem:int-grad}. We include their proofs for
convenience.

\subsubsection{\texorpdfstring{Proof of Corollary \protect\ref{coro:unif-bound}}{Proof of Corollary 2.2.}}

W.l.o.g. we assume $y=0$ and skip the dependence on $y$ in
the notation.
We distinguish two regimes: $|x|\leq\sqrt{T}$ and $|x|\geq\sqrt{T}$.

In the first case, we use (\ref{L11.19b}) and the
intermediate results (\ref{L11.7}) in the proof of Lemma \ref{L11a},
which yield
\begin{eqnarray*}
&&\mbox{for }d=2\qquad  \int_{|x|\leq\sqrt{T}}G_T^2(x) \,dx \lesssim T,
\\
&&\mbox{for }d>2\qquad  \int_{|x|\leq\sqrt{T}}G_T^q(x) \,dx \lesssim
\sqrt{T}{}^d\bigl(\sqrt{T}{}^{2-d}\bigr)^q ,
\end{eqnarray*}
and imply for $q=\frac{d-1}{d-2}\in(1,\frac{d}{d-2})$ by the
$L^2-L^\infty$ estimate
%
%
\begin{equation}\label{eq:proof-unif-bd-1}
G_T(x) \lesssim
\sqrt{T} \qquad\mbox{for }|x|\leq\sqrt{T}.
\end{equation}
For $|x|\geq\sqrt{T}$, we use the decay estimate (\ref{11.17}) of
Lemma \ref{L11a}
with ``$q=d,r=d(d+1)+1$''
\[
\int_{R\leq|x|\leq2R}G_T^d(x) \,dx \lesssim R^d\bigl(\sqrt
{T}R^{-1}\bigr)^{d(d+1)+1}=\sqrt{T}{}^d\bigl(\sqrt{T}R^{-1}\bigr)^{d^2+1},
\]
so that we may deduce
%
%
\begin{equation}\label{eq:proof-unif-bd-2}
G_T(x) \lesssim\bigl(\sqrt{T}R^{-1}\bigr)^{d+{1/d}}\sqrt{T}
\qquad\mbox
{for }R\leq|x|\leq2R.
\end{equation}
We then define $h_T\in L^1(\R^d)$ by
\[
h_T(x) \sim\cases{
\sqrt{T} 2^{-k(d+{1/d})}, &\quad $\sqrt{T}2^k \leq|x| \leq
\sqrt{T}2^{k+1},\qquad k\in\N$, \vspace*{2pt}\cr
\sqrt{T}, &\quad $|x| \leq\sqrt{T}$,}
\]
so that $G_T(x)\leq h_T(x)$ for all $x\in\Z^d$.
This concludes the proof since the factors in (\ref
{eq:proof-unif-bd-1}) and (\ref{eq:proof-unif-bd-2}) only
depend on $\alpha,\beta$ and $d$.

\subsubsection{\texorpdfstring{Proof of Corollary \protect\ref{coro:unif-bound-grad}}{Proof of Corollary 2.3.}}

We divide the proof in three steps. We first prove that the
Green function
$G_T(x,y)$ is symmetric so that $\nabla_x G_T(x,y)=\nabla_x G_T(y,x)$.
In the second step,
we show the uniform bound for $|x-y| \geq R$ sufficiently large, and in
the third
step for $|x-y|\leq R$.

\textit{Step} 1. Symmetry of $G_T$.

Let $y,\tilde y\in\Z^d$.
Testing the defining equation (\ref{eq:disc-Green}) with $x\mapsto
G_T(x,\tilde y)$ yields
\[
\int_{\Z^d}T^{-1}G_T(x,y)G_T(x,\tilde y)\,dx +\int_{\Z^d}\nabla
G_T(x,\tilde y)\cdot A(x)\nabla G_T(x,y)\,dx = G_T(y,\tilde y).
\]
Since $A$ is symmetric, the left-hand side of this identity is
symmetric in $y$
and $\tilde y$. Hence, the right-hand side is also symmetric, that is,
$G_T(y,\tilde y)=G_T(\tilde y,y)$.

Let $R\sim1$ be sufficiently large so that Lemma \ref
{lem:int-grad} applies.

\textit{Step} 2. Estimate for $|x-y|\geq R$.

For $q=2$, formula (\ref{eq:int-grad}) yields for all $k\in
\N$
\[
\int_{2^kR\leq|x-y|\leq2^{k+1}R}|\nabla_x G_T(x,y)|^2 \,dx \lesssim
(2^kR)^d((2^kR)^{1-d})^2=(2^kR)^{2-d} \stackrel{d\geq2}{\lesssim}
1 .
\]
Hence, by the discrete $L^2-L^\infty$ estimate, this shows
%
%
\begin{equation}\label{eq:pr:unifbd-grad-1}
|\nabla_x G_T(x,y)| \lesssim1 \qquad\mbox{for }|x-y|\geq R.
\end{equation}

\textit{Step} 3. Estimate for $|x-y|\leq R$.

We now use an {a priori} estimate.
Let
$i\in\{1,\ldots,d\}$ be fixed. We set $u(x):=G_T(x,y+\ee
_i)-G_T(x,y)=\nabla_{y_i}G_T(x,y)$.
This function solves the equation
%
%
\begin{equation}\label{eq:pr:unifbd-grad-2}
T^{-1}u-\nabla^*\cdot A\nabla u = f \qquad\mbox{in }\Z^d,
\end{equation}
where $f(x)= \delta(y+\ee_i-x)-\delta(y-x)$.
Since $f$ satisfies $\int_{\Z^d}f(x)=0$, one has by integration by parts,
ellipticity of $A$ and Poincar\'e's inequality
\begin{eqnarray*}
\int_{\Z^d}T^{-1}u(x)^2\,dx+\alpha\int_{\Z^d}|\nabla u(x)|^2\,dx
&\leq& \int_{\Z^d}f(x)u(x)\,dx \\
&=&\int_{|x-y|\leq R}f(x)\bigl(u(x)-\bar u_{\{|x-y|\leq R\}}\bigr)\,dx \\
&\lesssim& R \biggl(\int_{\{|x-y|\leq R\}}|\nabla u(x)|^2\,dx\biggr)^{1/2}.
\end{eqnarray*}
Hence,
\[
\int_{\Z^d}|\nabla u(x)|^2\,dx \lesssim R^2 \sim1.
\]
This shows that $\sup|\nabla u| \lesssim1$.
Therefore, for all $x$ such that $|x-y|\leq R$, we have using Step 2
and the fact that $R$ is of order 1
\[
|u(x)| \leq {R \sup}|\nabla u|+ {\sup_{|z-y|\geq R}}|u(z)| \lesssim
1.
\]
Recalling that $u(x)=\nabla_{y_i}G_T(x,y)$, we conclude by Step 1 that
this implies
$|\nabla_{y_i}G_T(y,x)|\lesssim1$, as desired.


\section{Proofs of the other auxiliary lemmas}\label{sec:proofs-aux}

\subsection{\texorpdfstring{Proof of Lemma \protect\ref{lem:var-estim}}{Proof of Lemma 2.3.}}

W.l.o.g. we may assume
%
%
\begin{equation}\label{5.1}
\sum_{i=1}^\infty\biggl\langle\sup_{a_i}
\biggl|\frac{\partial X}{\partial a_i}\biggr|^2\biggr\rangle<
\infty.
\end{equation}
Let $X_n$ denote the expected value of $X$ conditioned on
$a_1,\ldots,a_n$, that is,
\[
X_n(a_1,\ldots,a_n) = \langle X | a_1,\ldots,a_n \rangle.
\]
We will establish the following two inequalities for $n<\tilde n\in
\mathbb{N}$:
%
%
\begin{eqnarray}
\label{5.2}
\langle X_n^2\rangle-\langle X_n\rangle^2&\le&
\sum_{i=1}^n\biggl\langle\sup_{a_i} \biggl|\frac{\partial
X}{\partial
a_i}\biggr|^2\biggr\rangle\operatorname{var}[a_1],
\\
\label{5.5}
\langle(X_{\tilde n}-X_n)^2\rangle&\le&
\sum_{i=n+1}^{\tilde n}
\biggl\langle\sup_{a_i} \biggl|\frac{\partial X}{\partial a_i}
\biggr|^2\biggr\rangle\operatorname{var}[a_1].
\end{eqnarray}
Before proving (\ref{5.2}) and (\ref{5.5}), we draw the conclusion.
There is a slight technical difficulty due to the fact that there
are infinitely many random variables.

From (\ref{5.5}) and (\ref{5.1}), we learn that $\{X_n\}
_{n\uparrow\infty}$
is a Cauchy sequence in $L^2$ w.r.t. probability. Hence, there exists
a square integrable function $\tilde X$ of $a$ such that
%
%
\begin{equation}\label{5.3}
\lim_{n\uparrow\infty}\langle(\tilde X-X_n)^2\rangle= 0.
\end{equation}
By construction of $X_n$, (\ref{5.3}) implies
\[
\langle\tilde X | a_1,\ldots,a_n \rangle
=
\langle X | a_1,\ldots,a_n \rangle
\qquad\mbox{for a. e. } (a_1,\ldots,a_n) \mbox{ and all } n\in
\mathbb{N}.
\]
This means that the
random variables $X$ and
$\tilde X$ agree on all measurable finite rectangular cylindrical
sets, that is, measurable sets of the form $A_1\times\cdots\times A_n
\times\mathbb{R}\times\cdots,$ where $n$ is finite.
Since these sets are stable under intersection and generate the entire
$\sigma$-algebra of measurable sets, the
random variables $X$ and $\tilde X$ are uniquely determined by their
value on these sets \cite{Klenke-08}, Satz 14.12. Hence, the two
random variables coincide, yielding
%
%
\begin{equation}\label{5.4}
\tilde X = X \qquad\mbox{almost surely}.
\end{equation}
From (\ref{5.2}), (\ref{5.3}) and (\ref{5.4}), we obtain in the limit
$n\uparrow\infty$
as desired
\[
\operatorname{var}[X] = \langle X^2\rangle-\langle X\rangle^2
\le
\sum_{i=1}^\infty\biggl\langle\sup_{a_i}
\biggl|\frac{\partial X}{\partial a_i}\biggr|^2\biggr\rangle\operatorname{var}[a_1].
\]
We now turn to (\ref{5.2}) and (\ref{5.5}). Notice that we have
the decomposition
\[
\langle X_n^2\rangle-\langle X_n\rangle^2 =
\sum_{i=1}^n(\langle X_i^2\rangle-\langle X_{i-1}^2\rangle
),
\]
where we have set $X_0\dvtx\equiv\langle X\rangle$ so that
$\langle X_n\rangle^2=\langle X_0^2\rangle$. Hence, (\ref{5.2})
reduces to
%
%
\begin{equation}\label{5.6}
\langle X_i^2\rangle-\langle X_{i-1}^2\rangle
\le
\biggl\langle\sup_{a_i} \biggl|\frac{\partial X}{\partial a_i}
\biggr|^2\biggr\rangle
\operatorname{var}[a_1].
\end{equation}
Likewise,
\[
\langle(X_{\tilde n}-X_n)^2\rangle
= \langle X_{\tilde n}^2\rangle-\langle X_n^2\rangle
= \sum_{i=n+1}^{\tilde n}(\langle X_i^2\rangle-\langle
X_{i-1}^2\rangle),
\]
so that also (\ref{5.5}) reduces to (\ref{5.6}).

We finally turn to (\ref{5.6}). We note that by our
assumption that $\{a_i\}_{i\in\mathbb{N}}$
are i.i.d., we have
\begin{eqnarray*}
\langle X_i^2(a_1,\ldots,a_i)\rangle
&=&
\biggl\langle\int X_i^2(a_1,\ldots,a_{i-1},a_i') \beta(da_i')
\biggr\rangle,\\
X_{i-1}(a_1,\ldots,a_{i-1})
&=&
\int X_i(a_1,\ldots,a_{i-1},a_i'') \beta(da_i''),
\end{eqnarray*}
where $\beta$ denotes the distribution of $a_1$.
Hence, we obtain
\begin{eqnarray*}
&&\langle X_i^2\rangle-\langle X_{i-1}^2\rangle\\
&&\qquad=
\biggl\langle
\int X_i^2(a_1,\ldots,a_{i-1},a_i') \beta(da_i')
-\biggl(\int X_i(a_1,\ldots,a_{i-1},a_i'') \beta(da_i'')
\biggr)^2
\biggr\rangle\\
&&\qquad=
\biggl\langle{\int\int\frac{1}{2}}
\bigl(X_i(a_1,\ldots,a_{i-1},a_i')-X_i(a_1,\ldots
,a_{i-1},a_i'')
\bigr)^2 \beta(da_i')
\beta(da_i'')\biggr\rangle\\
&&\qquad\le
\biggl\langle{\int\int
\sup_{a_i'''}\biggl|\frac{\partial X_i}{\partial a_i}(a_1,\ldots
,a_{i-1},a_{i}''')\biggr|^2
\frac{1}{2}} (a_i'-a_i'')^2 \beta(da_i')
\beta(da_i'')\biggr\rangle\\
&&\qquad=
\biggl\langle
\sup_{a_i'''}\biggl|\frac{\partial X_i}{\partial a_i}(a_1,\ldots
,a_{i-1},a_{i}''')\biggr|^2
\biggr\rangle
\biggl(\int(a_i')^2 \beta(da_i')-\biggl(\int a_i'' \beta
(da_i'')
\biggr)^2\biggr)\\
&&\qquad=
\biggl\langle\sup_{a_i'''}\biggl|\frac{\partial X_i}{\partial
a_i}(a_1,\ldots,a_{i-1},a_{i}''')\biggr|^2 \biggr\rangle
\operatorname{var}[a_1].
\end{eqnarray*}
We conclude by noting that by the definition of $X_i$ and Jensen's inequality
\[
\biggl| \frac{\partial X_i}{\partial a_i}(a_1,\ldots,a_{i}) \biggr|^2
=
\biggl|\biggl\langle\frac{\partial X}{\partial a_i} \Big| a_1,\ldots,a_{i} \biggr\rangle
\biggr|^2
\le
\biggl\langle
\biggl|\frac{\partial X}{\partial a_i}\biggr|^2 \Big| a_1,\ldots,a_{i}
\biggr\rangle,
\]
so that
\begin{eqnarray*}
&&\biggl\langle
\sup_{a_i'}\biggl|\frac{\partial X_i}{\partial a_i}(a_1,\ldots
,a_{i-1},a_{i}')\biggr|^2
\biggr\rangle\\
&&\qquad\le\biggl\langle\!\biggl\langle\sup_{a_i'}
\biggl|\frac{\partial X}{\partial a_i}(a_1,\ldots
,a_{i-1},a_i',a_{i+1},\ldots)\biggr|^2
\Big| a_1,\ldots,a_{i}\biggr\rangle\!\biggr\rangle\\
&&\qquad=
\biggl\langle\sup_{a_i'}
\biggl|\frac{\partial X}{\partial a_i}(a_1,\ldots
,a_{i-1},a_i',a_{i+1},\ldots)\biggr|^2
\biggr\rangle.
\end{eqnarray*}


\subsection{\texorpdfstring{Proof of Lemma \protect\ref{lem:diff-Green}}{Proof of Lemma 2.5.}}

Let us divide the proof in four steps.

\textit{Step} 1. Proof of (\ref{eq:diff-Green}).

We recall the definition of the operator
\[
(Lu)(x)=\sum_{x',|x'-x|=1}a(x,x')\bigl(u(x)-u(x')\bigr).
\]
For convenience, we set $e=[z,z']$, $z'=z+\ee_i$.
We recall that $G_T(\cdot,y)$, $y\in\Z^d$, is defined via
%
%
\begin{equation}\label{eq:Green-fun}
(T^{-1}+L)G_T(\cdot,y)(x)=\delta(x-y),\qquad x\in\Z^d.
\end{equation}
Hence, we obtain by differentiating (\ref{eq:Green-fun})
\begin{eqnarray*}
&&\biggl((T^{-1}+L)\,\frac{\partial}{\partial a(e)}G_T(\cdot,y)
\biggr)(x)+\bigl(G_T(z,y)-G_T(z',y)\bigr)\delta(x-z) \\
&&\qquad{} +\bigl(G_T(z',y)-G_T(z,y)\bigr)\delta(x-z')=0,
\end{eqnarray*}
which, in view of (\ref{eq:Green-fun}), can be rewritten as
%
%
\begin{eqnarray}\label{7.2}\quad
&&(T^{-1}+L) \biggl(\frac{\partial}{\partial a(e)}G_T(\cdot
,y)+\bigl(G_T(z,y)-G_T(z',y)\bigr)G_T(\cdot,z) \nonumber\\[-8pt]\\[-8pt]
&&\qquad\hspace*{96pt}{}
+\bigl(G_T(z',y)-G_T(z,y)\bigr)G_T(\cdot,z')\biggr)\equiv0.\nonumber
\end{eqnarray}
From this, we would like to conclude
%
%
\begin{eqnarray}\label{7.1}
&&\frac{\partial}{\partial a(e)}G_T(\cdot
,y)+\bigl(G_T(z,y)-G_T(z',y)\bigr)G_T(\cdot
,z)\nonumber\\[-8pt]\\[-8pt]
&&\qquad{}
+\bigl(G_T(z',y)-G_T(z,y)\bigr)G_T(\cdot,z')\equiv0,\nonumber
\end{eqnarray}
which is nothing but (\ref{eq:diff-Green}).

In order to draw this conclusion, we will appeal to the
following uniqueness result in $L^2(\Z^d)$: any $u\in L^2(\Z^d)$
which satisfies
$
((T^{-1}+L)u)(x) = 0$ for all $x\in\mathbb{Z}^d
$
vanishes identically. However, we cannot apply this directly to
$u$ given by the left-hand side of (\ref{7.1}), since we do not know a
priori that
$\frac{\partial}{\partial a(e)} G_T(\cdot,y)$ is in $L^2(\Z^d)$.

For that\vspace*{1pt} purpose, we replace the derivative $\frac{\partial
}{\partial a(e)}$ by
the difference quotient. We thus fix a step size $h\not=0$ and
introduce the
abbreviations
\[
G_T(x,y) := G_T(x,y;a) \quad\mbox{and}\quad G_T'(x,y) := G_T(x,y,a'),
\]
where the coefficients $a'$ are defined by modifying $a$ only at edge
$e$ by
the increment $h$, that is,
\[
a'(e) = a(e)+h \quad\mbox{and}\quad a'(e') = a(e') \qquad\mbox{for
all } e'\not= e.
\]
We further denote by $L_T:=T+L_{a}$ and $L_T':=T+L_{a'}$ the operators
with coefficients $a$ and $a'$,
respectively. We mimic the derivation of (\ref{7.2}) on the discrete level:
from (\ref{eq:Green-fun}), we obtain
\begin{eqnarray*}
0&=&
\frac{1}{h}\bigl(L_T G_T(\cdot,y)-L_T' G_T'(\cdot,y)\bigr)\\
&=&
L_T\frac{1}{h}\bigl(G_T(\cdot,y)-G_T'(\cdot,y)\bigr)
+\frac{1}{h}(L_T-L_T')G_T'(\cdot,y)\\
&=&
L_T\frac{1}{h}\bigl(G_T(\cdot,y)-G_T'(\cdot,y)\bigr)
+ \bigl(G_T'(z,y)-G_T'(z',y)\bigr) \delta(\cdot-z )\\
&&{}
+\bigl(G_T'(z',y)-G_T'(z,y)\bigr) \delta(\cdot-z')\\
&=&
L_T\biggl(\frac{1}{h}\bigl(G_T(\cdot,y)-G_T'(\cdot,y)\bigr)
+ \bigl(G_T'(z,y)-G_T'(z',y)\bigr) G_T(\cdot,z )\\
&&\hspace*{124.3pt}{}
+\bigl(G_T'(z',y)-G_T'(z,y)\bigr) G_T(\cdot,z')\biggr).
\end{eqnarray*}
Since for fixed $h\not=0$,
\begin{eqnarray*}
u_h :\!&\equiv&
\frac{1}{h}\bigl(G_T(\cdot,y)-G_T'(\cdot,y)\bigr)+\bigl(G_T'(z,y)-G_T'(z',y)\bigr) G_T(\cdot,z )\\
&&{}
+\bigl(G_T'(z',y)-G_T'(z,y)\bigr) G_T(\cdot,z')
\end{eqnarray*}
does inherit the integrability properties of $G_T(\cdot,y)$ and
$G_T'(\cdot,y)$ from Corollary~\ref{coro:unif-bound}, we now may
conclude that $u_h\in L^2(\Z^d)$, and therefore $u_h\equiv0$, that is,
\begin{eqnarray*}
&&\frac{1}{h}\bigl(G_T(x,y)-G_T'(x,y)\bigr)
+\bigl(G_T'(z,y)-G_T'(z',y)\bigr) G_T(x,z )\\
&&\qquad{}+\bigl(G_T'(z',y)-G_T'(z,y)\bigr) G_T(x,z') = 0
\end{eqnarray*}
for every $x\in\mathbb{Z}^d$.
Since by Lemma 6, $G_T(x,y;\cdot)$ is continuous in $a(e)$,
we learn that $G_T(x,y;\cdot)$ is continuously differentiable w.r.t. $a(e)$
and that (\ref{eq:diff-Green}) holds.

We set for abbreviation
%
%
\begin{eqnarray}\label{eq:GT(e,e)}
G_T(x,e) &:=& G_T(x,z)-G_T(x,z'), \nonumber\\
G_T(e,y) &:=& G_T(z,y)-G_T(z',y), \\
G_T(e,e) &:=& G_T(z,z)+G_T(z',z')-G_T(z,z')-G_T(z',z).\nonumber
\end{eqnarray}

\textit{Step} 2. Proof of
%
%
\begin{eqnarray}\label{eq:diff2}
\frac{\partial}{\partial a(e)}G_T(x,e) &=&
-G_T(e,e)G_T(x,e),\nonumber\\[-8pt]\\[-8pt]
\frac{\partial}{\partial
a(e)}G_T(e,y) &=& -G_T(e,e)G_T(e,y).\nonumber
\end{eqnarray}
This is a consequence of (\ref{eq:diff-Green}) for $y=z,z'$:
\begin{eqnarray*}
&&{\frac{\partial}{\partial a(e)}
\bigl(G_T(x,z)-G_T(x,z')\bigr) }\\
&&\qquad= {\frac{\partial}{\partial a(e)}G_T(x,z)-\frac{\partial
}{\partial a(e)}G_T(x,z')}\\
&&\hspace*{-6.2pt}\qquad\stackrel{\mbox{\fontsize{8.36pt}{10.36pt}\selectfont{(\ref
{eq:diff-Green})}}}{=} { -\bigl(G_T(x,z)-G_T(x,z') \bigr)
\bigl(G_T(z,z)-G_T(z',z) \bigr)} \\
&&\qquad\quad{} +\bigl(G_T(x,z)-G_T(x,z') \bigr)\bigl(G_T(z,z')-G_T(z',z')\bigr)
\\
&&\qquad=-\bigl(G_T(z,z)+G_T(z',z')-G_T(z,z')-G_T(z,z') \bigr)
\\
&&\qquad\quad\hspace*{6.2pt}{}\times\bigl(G_T(x,z)-G_T(x,z') \bigr)
\end{eqnarray*}
and for $x=z,z'$, respectively.

\textit{Step} 3. Conclusion.

Note that Corollary \ref{coro:unif-bound-grad} implies
%
%
\begin{equation}\label{eq:G+-}
|G_T(e,e)|\lesssim1.
\end{equation}
The combination of (\ref{eq:diff2}) with (\ref{eq:G+-}) yields
\[
\biggl| \frac{\partial}{\partial a(e)}G_T(x,e)\biggr|\lesssim
|G_T(x,e)|,\qquad \biggl| \frac{\partial}{\partial
a(e)}G_T(e,y)
\biggr|\lesssim|G_T(e,y)|.
\]
Since $a(e)$ is bounded, this also yields
\[
{\sup_{a(e)}} |G_T(x,e)|\sim|G_T(x,e)|,\qquad {\sup_{a(e)}}
|G_T(e,y)|\sim
|G_T(e,y)|,
\]
which is nothing but (\ref{eq:bd-G(x,e)}).


\subsection{\texorpdfstring{Proof of Lemma \protect\ref{lem:diff-phi}}{Proof of Lemma 2.4.}}

We recall that $e=[z,z']$, $z'=z+\ee_i$.

\textit{Step} 1. Proof of (\ref{eq:diff-phi-1}).

We first give a heuristic argument for (\ref{eq:diff-phi-1})
based on the defining equation
%
%
\begin{equation} \label{eq:pr-diff-phi-4}
T^{-1}\phi_T(x)-\bigl(\nabla^*\cdot A\bigl(\nabla\phi_T(x)+\xi\bigr)\bigr)(x) = 0.
\end{equation}
Differentiating (\ref{eq:pr-diff-phi-4}) w.r.t. $a(e)$ yields as in
Step 1 of the proof of Lemma \ref{lem:diff-Green}
%
%
\begin{eqnarray} \label{eq:pr-diff-phi-5}
&&T^{-1}\,\frac{\partial\phi_T}{\partial a(e)}(x)-\biggl(\nabla^*\cdot
A\nabla\,\frac{\partial\phi_T}{\partial a(e)}\biggr)(x) \nonumber\\[-8pt]\\[-8pt]
&&\qquad\hspace*{0pt}{}-\bigl(\nabla
_i\phi
_T(z)+\xi_i\bigr)\bigl(\delta(x-z)-\delta(x-z')\bigr) = 0.\nonumber
\end{eqnarray}
Provided we have $\frac{\partial\phi_T}{\partial a(e)}\in L^2(\Z^d)$,
this yields by definition of $G_T$
\[
\frac{\partial\phi_T}{\partial a(e)}(x) =
-\bigl(\nabla_i\phi_T(z)+\xi_i\bigr)\bigl(G_T(x,z')-G_T(x,z)\bigr),
\]
which is (\ref{eq:diff-phi-1}).

In order to turn the above into a rigorous argument, we need
to argue that $\phi_T(x)$ is differentiable w.r.t. $a(e)$ and that
$\frac{\partial\phi_T}{\partial a(e)}\in L^2(\Z^d)$.
Starting\vspace*{1pt} point is the representation formula from Step 2 of the proof
of Lemma \ref{lem:depend-coeff}, that is,
%
%
\begin{equation}\label{eq:pr-diff-phi-7}
\phi_T(x) = \int_{\Z^d}G_T(x,y)\nabla^*\cdot(A(y)\xi) \,dy.
\end{equation}
Combined with Corollary \ref{coro:unif-bound}, (\ref{eq:pr-diff-phi-7})
and (\ref{eq:diff-Green}) in Lemma \ref{lem:diff-Green} show that
$\phi
_T(x)$ is differentiable w.r.t. $a(e)$.
We may now switch the order of the differentiation and the sum as follows:
%
%
\begin{eqnarray}\label{eq:pr-diff-phi-9}\qquad
\frac{\partial\phi_T}{\partial a(e)}(x)&=&-\nabla_{z_i} G_T(x,z)\xi
_i\nonumber\\
&&{}-\int_{\Z^d}\nabla_{z_i}G_T(x,z)\nabla_{z_i}G_T(z,y)\nabla
^*\cdot
(A(y)\xi) \,dy \\
&=&-\underbrace{\nabla_{z_i} G_T(x,z)}_{\in L^2_x(\Z^d)}
\biggl(\xi
_i+\int_{\Z^d}\underbrace{\nabla_{z_i}G_T(z,y)}_{\in L^1_y(\Z
^d)}\underbrace{\nabla^*\cdot(A(y)\xi)}_{\in L^\infty(\Z
^d)}\,dy\biggr),\nonumber
\end{eqnarray}
since $G_T(\cdot,z)\in L^2(\Z^d)$ by definition of the Green's
function, $G_T(z,\cdot)\in L^1(\Z^d)$ by Corollary \ref
{coro:unif-bound} and $A$ is bounded. This proves that $\frac{\partial
\phi_T}{\partial a(e)}\in L^2(\Z^d)$.

\textit{Step} 2. Proof of
%
%
\begin{eqnarray}
\label{eq:diff-phi-bounded}
{\sup_{a(e)}}|\phi_T(x)|&\lesssim&|\phi_T(x)|+\bigl(|\nabla_i \phi
_T(z)|+1\bigr) |\nabla_{z_i} G_T(z,x)|,\\
\label{eq:diff-phi-induc=0}
{\sup_{a(e)}}\biggl|\frac{\partial\phi_T(x)}{\partial a(e)}
\biggr|&\lesssim&\bigl(|\nabla_i \phi_T(z)|+1\bigr)| \nabla_{z_i}
G_T(z,x)|.
\end{eqnarray}
We argue that it is enough to prove (\ref{eq:diff-phi-3}). Indeed, the
combination of (\ref{eq:diff-phi-1}), (\ref{eq:bd-G(x,e)}) and (\ref
{eq:diff-phi-3}) with the boundedness of $a$ implies (\ref
{eq:diff-phi-bounded}) and (\ref{eq:diff-phi-induc=0}).
In order to prove (\ref{eq:diff-phi-3}), we proceed as follows
%
%
\begin{eqnarray}\label{eq:pr-diff-phi-10}\qquad
-\biggl(\nabla_i \frac{\partial\phi_T}{\partial a(e)}\biggr)(z)&=&
\frac
{\partial\phi_T}{\partial a(e)}(z)- \frac{\partial\phi_T}{\partial
a(e)}(z') \nonumber\\
&\stackrel{\mbox{\fontsize{8.36pt}{10.36pt}\selectfont{(\ref
{eq:diff-phi-1})}}}{=}& \bigl(\nabla_{i}\phi_T(z)+\xi_i\bigr)\bigl(G_T(z,z)-G_T(z,z')\bigr)
\nonumber\\
&&{} -\bigl(\nabla_{i}\phi_T(z)+\xi_i\bigr)\bigl(G_T(z',z)-G_T(z',z')\bigr)
\nonumber\\[-8pt]\\[-8pt]
&=&\bigl(\nabla_{i}\phi_T(z)+\xi_i\bigr)\nonumber\\
&&{}\times\bigl(G_T(z,z)-G_T(z,z')-G_T(z',z)+G_T(z',z')\bigr)
\nonumber\\
&=& \bigl(\nabla_{i}\phi_T(z)+\xi_i\bigr)G_T(e,e),\nonumber
\end{eqnarray}
where we used the abbreviation
\[
G_T(e,e) = G_T(z,z)-G_T(z,z')-G_T(z',z)+G_T(z',z').
\]
Recalling that Corollary \ref{coro:unif-bound-grad} implies
\[
G_T(e,e) \lesssim1,
\]
inequality (\ref{eq:diff-phi-3}) follows now from (\ref
{eq:pr-diff-phi-10}) and the boundedness of $a$.

\textit{Step} 3. Proof of (\ref{eq:diff-phi-2}).

For $n\geq0$, the chain rule yields
\[
\frac{\partial[\phi_T(x)^{n+1}]}{\partial a(e)}
=(n+1)\phi_T(x)^{n} \,\frac{\partial\phi_T(x)}{\partial a(e)}.
\]
Using (\ref{eq:diff-phi-bounded}) and (\ref{eq:diff-phi-induc=0}),
this implies
\begin{eqnarray*}
\sup_{a(e)}\biggl|\frac{\partial[\phi_T(x)^{n+1}]}{\partial
a(e)}\biggr|
&\lesssim&\bigl(|\phi_T(x)|+\bigl(|\nabla_i \phi_T(z)|+1\bigr)
|\nabla_{z_i} G_T(z,x)|\bigr)^n\\
&&{} \times\bigl(\bigl(|\nabla_i \phi_T(z)
|+1\bigr)|
\nabla_{z_i} G_T(z,x)|\bigr),
\end{eqnarray*}
which turns into (\ref{eq:diff-phi-2}) using Young's inequality.

\subsection{\texorpdfstring{Proof of Lemma \protect\ref{lem:depend-coeff}}{Proof of Lemma 2.6.}}

We first prove the claim for $G_T$ and deduce the result for $\phi_T$
appealing to an integral representation
using the Green's function.

\textit{Step} 1. Properties of $G_T$.

The product topology is the topology of componentwise convergence.
Hence, we consider an arbitrary sequence $\{a_\nu\}_{\nu\uparrow
\infty}
\subset{\mathcal A}_{\alpha\beta}$
of coefficients such that
%
%
\begin{equation}\label{6.1}
\lim_{\nu\uparrow\infty}a_\nu(e) = a(e) \qquad\mbox{for all
edges } e.
\end{equation}
Fix $y\in\mathbb{Z}^d$; by the uniform bounds on $G_T(\cdot,y;a_\nu)$
from Corollary \ref{coro:unif-bound},
we can select a subsequence $\nu'$ such that
%
%
\begin{equation}\label{6.2}
u_T(x) := \lim_{\nu'\uparrow\infty}G_T(x,y;a_{\nu'})
\qquad\mbox{exists for all } x\in\mathbb{Z}^d.
\end{equation}
It remains to argue that $u_T(x)=G_T(x,y;a)$.
Because of (\ref{6.1}) and (\ref{6.2}), we can pass to the limit
in $(T^{-1}G_T(\cdot,y;a_{\nu'})+L_{a_{\nu'}} G_T(\cdot,y;a_{\nu
'}))(x)=\delta(x-y)$ to obtain
%
%
\begin{equation}\label{6.3}
(T^{-1}u_T+L_{a} u_T)(x) = \delta(x-y) \qquad\mbox{for all } x\in
\mathbb{Z}^d.
\end{equation}
Moreover, the uniform decay of $G_T(\cdot,y;a_\nu)$ from Corollary
\ref{coro:unif-bound} is preserved in the limit, so that $u_T\in
L^1(\Z^d)\subset L^2(\Z^d)$. Note that Riesz's representation theorem
on $L^2(\Z^d)$ yields uniqueness for the solution of (\ref{6.3}) in
$L^2(\Z^d)$. Hence, we conclude as desired that
$u_T(\cdot)=G_T(\cdot,y;a)$. Borel measurability of $G_T(x,y;\cdot)$ in
the sense of Lemma \ref {lem:var-estim} follows from continuity w.r.t.
the product topology, cf. \cite{Klenke-08}, Satz~14.8.

\textit{Step} 2. Properties of $\phi_T$.

Corollary \ref{coro:unif-bound} ensures that $G_T(x,\cdot
)\in
L^1(\Z^d)$ for all $x\in\Z^d$ and one may then define a function
$\tilde\phi_T$ by
%
%
\begin{equation}\label{eq:measu-phi}
\tilde\phi_T(x)=\int_{\Z^d}G_T(x,y)\nabla^*\cdot(A(y) \xi_i) \,dy.
\end{equation}
Since $G_T(\cdot+z,\cdot+z)$ has the same law as $G_T(\cdot,\cdot)$ by
uniqueness of the
Green's function and joint stationarity of the coefficient $A$, $\tilde
\phi_T(\cdot+z)$ has the
same law as $\tilde\phi_T$. This shows that $\tilde\phi_T$ is stationary.
In addition, $\tilde\phi_T$ is a solution of (\ref{eq:app-corr}) by
construction.
Hence, by the uniqueness of stationary solutions of (\ref
{eq:app-corr}), $\tilde\phi_T = \phi_T$ almost surely,
so that by the measurability properties we may assume $\tilde\phi_T
\equiv\phi_T$.

Introducing for $R\geq1$
\[
\phi_{T,R}(x):=\int_{|y|\leq R}G_T(x,y)\nabla^*\cdot(A(y) \xi_i) \,dy,
\]
one may rewrite (\ref{eq:measu-phi}) as
%
%
\begin{equation}\label{eq:measu-phi2}
\phi_T(x)=\lim_{R\to\infty} \phi_{T,R}(x).
\end{equation}
From Step 1, $\phi_{T,R}(x)$ is a continuous function of $a$ since
$G_T(x,y)$ is and the formula for $\phi_{T,R}(x)$ involves only a \textit{finite} number of operations.
Note that Corollary~\ref{coro:unif-bound} implies that
\[
\lim_{R\uparrow\infty} \sup_{a\in\mathcal{A}_{\alpha\beta}}\int
_{|y|> R}G_T(x,y;a) \,dy = 0.
\]
Hence, the convergence in (\ref{eq:measu-phi2}) is uniform in $a$ and
the continuity of $\phi_{T,R}$ in $a$ is preserved
at the limit.
Therefore, $\phi_T$ (and continuous functions thereof) are continuous
with respect to the product topology, and hence Borel measurable.


\subsection{\texorpdfstring{Proof of Lemma \protect\ref{lem:Lp-grad}}{Proof of Lemma 2.7.}}

We first sketch the proof in the continuous case, that is, with
$\mathbb{Z}^d$ replaced by $\mathbb{R}^d$.

\textit{Step} 1. Continuous version.

Starting point is
the defining equation (\ref{eq:app-corr}) of the corrector $\phi_T$
in its
continuous version, that is,
%
%
\begin{equation}\label{8.1}
T^{-1}\phi_T-\nabla\cdot A(\nabla\phi_T+\xi) = 0 \qquad\mbox
{in }
\mathbb{R}^d.
\end{equation}
We multiply (\ref{8.1}) with $\phi_T^{n+1}$ and obtain
by Leibniz' rule:
%
%
\begin{eqnarray}\label{8.2}\quad
0&=&T^{-1}\phi_T^{n+2}+\bigl(-\nabla\cdot A(\nabla\phi_T+\xi
)\bigr)
\phi_T^{n+1}\nonumber\\
&=&T^{-1}\phi_T^{n+2}
-\nabla\cdot\bigl(\phi_T^{n+1} A(\nabla\phi_T+\xi)\bigr)
+ \nabla\phi_T^{n+1}\cdot A(\nabla\phi_T+\xi)\nonumber\\[-8pt]\\[-8pt]
&=&T^{-1}\phi_T^{n+2}
-\nabla\cdot\bigl(\phi_T^{n+1} A(\nabla\phi_T+\xi)\bigr)\nonumber\\
&&{}
+ (n+1) \phi_T^n \nabla\phi_T\cdot A(\nabla\phi_T+\xi).\nonumber
\end{eqnarray}
We then take the expected value. Since the random fields $A$ and
$\phi_T$ are jointly stationary, and thus also $\phi_T^{n+1}
A(\nabla
\phi_T+\xi)$,
we obtain
\[
\langle T^{-1}\phi_T^{n+2} \rangle+(n+1)\langle\phi_T^n \nabla
\phi_T\cdot A(\nabla\phi_T+\xi) \rangle= 0,
\]
and therefore
\[
\langle\phi_T^n \nabla\phi_T\cdot A(\nabla\phi_T+\xi) \rangle
\leq0
\]
since $n+2$ is even.
By the uniform ellipticity of $A$ and since $\phi_T^n\ge0$ ($n$ is
even) and $|\xi|=1$, this yields
the estimate
\[
\langle\phi_T^n|\nabla\phi_T|^2\rangle\lesssim
\langle\phi_T^n|\nabla\phi_T|\rangle.
\]
Applying Cauchy--Schwarz' inequality in probability
on the right-hand side
of this inequality yields the continuum version of (\ref
{eq:phi-grad-phi}), that is,
\[
\langle\phi_T^n|\nabla\phi_T|^2\rangle\lesssim
\langle\phi_T^n\rangle.
\]
We now turn to our discrete case.

\textit{Step} 2. Discrete version.

We need a discrete version
of the Leibniz rule
$\nabla\cdot(fg) = f \nabla\cdot g + \nabla f\cdot g$
used in (\ref{8.2}).
Let $f\in L^2_{\mathrm{loc}}(\Z^d)$ and $g\in L^2_{\mathrm{loc}}(\Z
^d,\R
^d)$, then
this formula is replaced by
%
%
\begin{eqnarray}\label{eq:Leib-disc-0}
\nabla^*\cdot(fg)(z)&=&\sum_{j=1}^d \bigl(f(z)[g(z)]_j-f(z-\ee
_j)[g(z-\ee_j)]_j \bigr) \nonumber\\[-8pt]\\[-8pt]
&=& f(z)\nabla^*\cdot g(z)+\sum_{j=1}^d
\nabla^*_jf(z)[g(z-\ee_j)]_j.\nonumber
\end{eqnarray}
We also need a substitute for the identity
$\nabla\phi_T^{n+1}=(n+1)\phi_T^n \nabla\phi_T$ used in (\ref
{8.2}). This
substitute is provided by the two calculus estimates
%
%
\begin{eqnarray}
\label{8.5}
(\tilde\phi^{n+1}-\phi^{n+1}) (\tilde\phi-\phi)&\gtrsim&
(\tilde\phi^n+\phi^n) (\tilde\phi-\phi)^2,\\
\label{8.8}
|\tilde\phi^{n+1}-\phi^{n+1}|&\lesssim&
(\tilde\phi^n+\phi^n) |\tilde\phi-\phi|.
\end{eqnarray}
For the convenience of the reader, we sketch their proof:
by the well-known formula for $\tilde\phi^{n+1}-\phi^{n+1}$,
they are equivalent to
\[
\sum_{m=0}^n\phi^m {\tilde\phi}^{n-m}\sim
\tilde\phi^n+\phi^n.
\]
By homogeneity, we may assume $\tilde\phi=1$, so that the above turns into
\[
\sum_{m=0}^n\phi^m \sim1+\phi^n.
\]
The upper estimate is obvious by H\"older's inequality
since $n$ is even. Also for the lower bound, we use
the evenness of $n$ to rearrange the sum as follows:
\begin{eqnarray*}
\sum_{m=0}^n\phi^m
&=&{ \frac{1}{2}} 1 +
{ \frac{1}{2}} (1+2\phi+\phi^2) +
{ \frac{1}{2}} \phi^2 (1+2\phi+\phi^2) + \cdots\\
&&{} +
{ \frac{1}{2}} \phi^{n-2} (1+2\phi+\phi^2) +
{ \frac{1}{2}} \phi^n\\
&\ge&
{ \frac{1}{2}} (1+\phi^n).
\end{eqnarray*}
After these motivations and preparations, we turn to the
proof of Lemma \ref{lem:Lp-grad} proper.
With $f(z):=\phi_T^{n+1}(z)$ and $g(z):=A(\nabla\phi_T+\xi)(z)$,
(\ref
{eq:Leib-disc-0}) turns into
\begin{eqnarray*}
&&\nabla^*\cdot\bigl(\phi_T^{n+1}(z) A(\nabla\phi_T+\xi
)(z)\bigr)
\\
&&\qquad=\phi_T^{n+1}(z)\nabla^*\cdot A(\nabla\phi_T+\xi)(z)+\sum
_{j=1}^d\nabla^*_j \phi_T^{n+1}(z) \mathop
{\underbrace{[A(\nabla
\phi_T+\xi)(z-\ee_j)]_j}_{
=a(z-\ee_j,z)(\nabla_j \phi_T(z-\ee_j)+\xi_j)}}_
{=a(z-\ee_j,z)(\nabla_j^* \phi_T(z)+\xi_j)}.
\end{eqnarray*}
Hence,
%
%
\begin{eqnarray}\label{eq:Leib-disc-1}
&&-\phi_T^{n+1}(z)\nabla^*\cdot A(\nabla\phi_T+\xi)(z)\nonumber\\
&&\qquad = \sum
_{j=1}^d\nabla^*_j \phi_T^{n+1}(z)a(z-\ee_j,z)\bigl(\nabla_j^* \phi
_T(z)+\xi
_j\bigr)\\
&&\qquad\quad{} -\nabla^*\cdot\bigl(\phi_T^{n+1}(z)
A(\nabla\phi_T+\xi)(z)\bigr).\nonumber
\end{eqnarray}
Multiplying (\ref{eq:app-corr}) with $\phi_T^{n+1}(z)$ and using
(\ref
{eq:Leib-disc-1}) emulate
(\ref{8.2}) and yield
%
%
\begin{eqnarray}\label{8.2-disc}
0 &=& T^{-1}\phi_T^{n+2}(z)-\nabla^*\cdot\bigl(\phi
_T^{n+1}(z) A(\nabla\phi_T+\xi)(z)\bigr)
\nonumber\\[-8pt]\\[-8pt]
&&{} +\sum_{j=1}^d\nabla^*_j \phi_T^{n+1}(z)a(z-\ee_j,z)\bigl(\nabla_j^*
\phi
_T(z)+\xi_j\bigr).\nonumber
\end{eqnarray}
Taking the expectation of (\ref{8.2-disc}) and noting that $\phi
_T^{n+2}\geq0$, we obtain as for
the continuous case
%
%
\begin{eqnarray}\label{8.9}
&&\Biggl\langle\sum_{j=1}^da(z-\ee_j,z)\nabla^*_j \phi_T^{n+1}(z)\nabla
_j^* \phi_T(z) \Biggr\rangle\nonumber\\[-8pt]\\[-8pt]
&&\qquad\lesssim\Biggl\langle\sum_{j=1}^d a(z-\ee
_j,z)|\nabla_j^* \phi_T^{n+1}(z)| \Biggr\rangle.\nonumber
\end{eqnarray}
On the one hand, we have
%
%
\begin{eqnarray}\label{8.10}\hspace*{32pt}
&&\sum_{j=1}^da(z-\ee_j,z)\nabla^*_j \phi_T^{n+1}(z)\nabla_j^*
\phi_T(z)
\nonumber\\
&&\qquad= \sum_{j=1}^da(z-\ee_j,z)\bigl(\phi_T^{n+1}(z)-\phi_T^{n+1}(z-\ee_j)\bigr)\bigl(
\phi_T(z)- \phi_T(z-\ee_j)\bigr)
\\
&&\qquad\stackrel{\mbox{\fontsize{8.36pt}{10.36pt}\selectfont{(\ref
{8.5})}}}{\gtrsim}
\sum_{j=1}^d \bigl(\phi_T^{n}(z)+\phi_T^{n}(z-\ee_j)\bigr)
\bigl(\phi_T(z)-\phi_T(z-\ee_j)\bigr)^2.\nonumber
\end{eqnarray}
On the other hand, we observe
%
%
\begin{eqnarray}\label{8.11}
&&\sum_{j=1}^da(z-\ee_j,z)|\nabla_j^* \phi
_T^{n+1}(z)|\nonumber\\
&&\qquad=\sum_{j=1}^da(z-\ee_j,z)|\phi_T^{n+1}(z)-\phi_T^{n+1}(z-\ee_j)|
\\
&&\qquad\hspace*{-6.08pt}\stackrel{\mbox{\fontsize{8.36pt}{10.36pt}\selectfont{(\ref
{8.8})}}}{\lesssim}
\sum_{j=1}^d\bigl(\phi_T^{n}(z)+\phi_T^{n}(z-\ee_j)\bigr)
|\phi_T(z)-\phi_T(z-\ee_j)|.\nonumber
\end{eqnarray}
Now (\ref{8.9}), (\ref{8.10}) and (\ref{8.11}) combine to
\begin{eqnarray*}
&&\Biggl\langle\sum_{j=1}^d \bigl(\phi_T^{n}(z)+\phi_T^{n}(z-\ee_j)\bigr) \bigl(\phi
_T(z)-\phi_T(z-\ee_j)\bigr)^2 \Biggr\rangle\\
&&\qquad\lesssim
\sum_{j=1}^d\bigl\langle\bigl(\phi_T^{n}(z)+\phi_T^{n}(z-\ee_j)\bigr) |\phi
_T(z)-\phi_T(z-\ee_j)| \bigr\rangle.
\end{eqnarray*}
By stochastic homogeneity, this reduces to
\begin{eqnarray*}
&&
\Biggl\langle\sum_{j=1}^d\bigl(\phi_T^{n}(\ee_j)+\phi_T^{n}(0)\bigr)
\bigl(\phi_T(\ee_j)-\phi_T(0)\bigr)^2\Biggr\rangle
\\
&&\qquad\lesssim
\Biggl\langle\sum_{j=1}^d\bigl(\phi_T^{n}(\ee_j)+\phi_T^{n}(0)\bigr)
|\phi_T(\ee_j)-\phi_T(0)|\Biggr\rangle.
\end{eqnarray*}
An application of Cauchy--Schwarz inequality yields
\[
\Biggl\langle\sum_{j=1}^d\bigl(\phi_T^{n}(\ee_j)+\phi_T^{n}(0)\bigr)
\bigl(\phi_T(\ee_j)-\phi_T(0)\bigr)^2\Biggr\rangle
\lesssim
\Biggl\langle\sum_{j=1}^d\bigl(\phi_T^{n}(\ee_j)+\phi_T^{n}(0)\bigr)
\Biggr\rangle.
\]
A last application of stochastic homogeneity gives as desired
\[
\Biggl\langle\phi_T^n(0) \sum_{j=1}^d\bigl(
\bigl(\phi_T(\ee_j)-\phi_T(0)\bigr)^2+\bigl(\phi_T(0)-\phi_T(-\ee_j)\bigr)^2
\bigr)
\Biggr\rangle
\lesssim
\langle\phi_T^{n}(0)\rangle.
\]


\subsection{\texorpdfstring{Proof of Lemma \protect\ref{lem:hh}}{Proof of Lemma 2.10.}}

The proof relies on a doubly dyadic decomposition of space.
First note that by symmetry,
\begin{eqnarray*}
\int_{|z|\leq|z-x|} h_T(z)h_T(z-x)\,dz&=&\int_{|z|\geq|z-x|}
h_T(z)h_T(z-x)\,dz \\
&\geq&\frac{1}{2}\int_{\Z^d} h_T(z)h_T(z-x)\,dz.
\end{eqnarray*}
Hence, it is enough to consider
\[
\int_{|x|\leq R} \int_{|z|\leq|z-x|}h_T(z)h_T(z-x)\,dz\,dx.
\]
In the three first steps, we treat the case $d>2$. We then sketch the
modification for $d=2$ in the last step.
Let $\tilde{R}\sim1$ be such that (\ref{eq:assump-h}) holds with a
constant independent of $R$ for all $R\geq\tilde R/2$.

\textit{Step} 1. Proof of
%
%
\begin{equation}\label{eq:hh-step1-1}\quad
\int_{R< |x|\leq2R} \int_{|z|\leq|z-x|}h_T(z)h_T(z-x)\,dz\,dx \lesssim
R^2\qquad \mbox{for }R \geq2\tilde R.
\end{equation}
Let $N\in\N$ be such that $\tilde R \leq2^{-N}R \leq2\tilde R$. We
then decompose the sum over ${|z|\leq|z-x|}$ into three contributions:
$R/2 < |z|$,
a dyadic decomposition\vspace*{1pt} for $\tilde R < |z|\leq R/2$ and a remainder on
$|z|\leq\tilde R$. More precisely,
\begin{eqnarray*}
&&\int_{R< |x|\leq2R} \int_{|z|\leq|z-x|}h_T(z)h_T(z-x)\,dz\,dx
\\
&&\qquad= \int_{R< |x|\leq2R}\int_{R/2< |z|\leq|z-x|}h_T(z)h_T(z-x)\,dz\,dx \\
&&\qquad\quad{}+\sum_{n=1}^N\int_{R< |x|\leq2R}\int_{\{2^{-(n+1)}R< |z|\leq
2^{-n}R\}\cap\{|z|\leq|z-x|\}}h_T(z)h_T(z-x)\,dz\,dx\\
&&\qquad\quad{} +\int_{R< |x|\leq2R}\int_{\{|z|\leq2^{-(N+1)}R\}\cap\{|z|\leq
|z-x|\}}h_T(z)h_T(z-x)\,dz\,dx \\
&&\qquad\leq \int_{|x|\leq2R}\underbrace{\int_{R/2< |z|\leq|z-x|}
h_T(z)h_T(z-x)\,dz}_{=I_1}\,dx \\
&&\qquad\quad{} + \sum_{n=1}^N \underbrace{\int_{R< |x|\leq2R}\int_{2^{-(n+1)}R<
|z|\leq2^{-n}R}h_T(z)h_T(z-x)\,dz\,dx}_{=I_2(n)}\\
&&\qquad\quad{}+{}\underbrace{\int_{R< |x|\leq2R}\int_{|z|\leq\tilde
R}h_T(z)h_T(z-x)\,dz\,dx}_{=I_3(N)}.
\end{eqnarray*}
We use Young's inequality, a dyadic decomposition of $\{
|z|>R/2\}$, and the assumption (\ref{eq:assump-h}) to bound $I_1$:
\begin{eqnarray*}
I_1 &\leq& \frac{1}{2}\biggl( \int_{R/2< |z|} h_T(z)^2 \,dz +\int_{R/2<
|z-x|} h_T(z-x)^2 \,dz \biggr)\\
& = & \sum_{k=-1}^\infty\int_{2^kR< |z|\leq2^{k+1}R}h_T^2(z)\,dz
\stackrel{\mbox{\fontsize{8.36pt}{10.36pt}\selectfont{(\ref
{eq:assump-h})}}}{\lesssim} \sum_{k=-1}^\infty\biggl(\frac
{1}{2^{d-2}}\biggr)^k R^{2-d}
\lesssim R^{2-d}.
\end{eqnarray*}
In order to bound $I_2(n)$, we will use the following fact
%
%
\begin{equation}\label{eq:summation}
\bigl( |x|> R \mbox{ and } |z|\leq\tfrac{1}{2}R \bigr)
\Rightarrow\bigl( |z-x|> \tfrac{1}{2}R \bigr).
\end{equation}
We have by Cauchy--Schwarz inequality
\begin{eqnarray*}
I_2(n)
&\leq& \biggl(\int_{|x|\leq2R} \int_{2^{-(n+1)}R< |z|\leq
2^{-n}R}h_T(z)^2\,dz\,dx
\\
&&\hspace*{4.1pt}{}\times\int_{R< |x|\leq2R} \int_{|z|\leq2^{-n}R}h_T(z-x)^2\,dz\,dx \biggr)^{1/2}
\\
&\lesssim& \biggl(R^d \underbrace{\int_{2^{-(n+1)}R< |z|\leq
2^{-n}R}h_T(z)^2\,dz}_{\stackrel{\mbox{\fontsize
{8.36pt}{10.36pt}\selectfont{(\ref{eq:assump-h})}}}{\lesssim}
(2^{-n}R)^{2-d}}
\\
&&\hspace*{3.4pt}{}\times
\underbrace{\int_{R< |x|\leq2R}
\int_{|z|\leq2^{-n}R}h_T(z-x)^2\,dz\,dx}_{
\stackrel{\mbox{\fontsize{8.36pt}{10.36pt}\selectfont{(\ref
{eq:summation})}}}{\lesssim} \underbrace{\int_{|z|\leq2^{-n}R} \int
_{R/2 < |z-x| \leq5R/2}h_T(z-x)^2\,dx\,dz}_{
\stackrel{\mbox{\fontsize{8.36pt}{10.36pt}\selectfont{(\ref
{eq:assump-h})}}}{\lesssim}
\displaystyle \int_{|z|\leq2^{-n}R}R^{2-d} \,dz=(2^{-n}R)^dR^{2-d}}
} \biggr)^{1/2} \\
&{\lesssim} & 2^{-n}R^{2}.
\end{eqnarray*}
We proceed the same way to bound $I_3(N)$. Recalling that $R\geq2
\tilde R\sim1$, it holds that $|z|\leq\tilde R \Rightarrow
|z|\leq R/2$. Hence, we are in position to use (\ref{eq:summation}) and
we obtain
\begin{eqnarray*}
I_3(N) &\leq& \biggl(\int_{|x|\leq2R} \int_{|z|\leq\tilde
R}h_T(z)^2\,dz\,dx
\int_{R< |x|\leq2R} \int_{|z|\leq\tilde R}h_T(z-x)^2\,dz\,dx \biggr)^{1/2}
\\
&\lesssim& \biggl(R^d \underbrace{\int_{|z|\leq\tilde
R}h_T(z)^2\,dz}_{
\stackrel{\mbox{\fontsize{8.36pt}{10.36pt}\selectfont{(\ref
{eq:assump-h-R=1})}}}{\lesssim}1} \
\underbrace{\int_{R< |x|\leq2R}
\int_{|z|\leq\tilde R}h_T(z-x)^2\,dz\,dx}_{
\stackrel{\mbox{\fontsize{8.36pt}{10.36pt}\selectfont{(\ref
{eq:summation})}}}{\lesssim} \underbrace{\int_{|z|\leq\tilde R}
\int
_{R/2 < |z-x| \leq5R/2}h_T(z-x)^2\,dx\,dz}_{
\stackrel{\mbox{\fontsize{8.36pt}{10.36pt}\selectfont{(\ref
{eq:assump-h})}}}{\lesssim}
\displaystyle \int_{|z|\leq\tilde{R}}R^{2-d} \,dz\sim R^{2-d}}
} \biggr)^{1/2} \\
&{\lesssim} &R.
\end{eqnarray*}
Since $\sum_{n=1}^\infty2^{-n}R^2 \sim R^2$ and $|\{|x|\leq2R\}
|R^{2-d}\sim R^{2}$, the bounds on $I_1$, $I_2(n)$ and $I_3(N)$
imply the claim (\ref{eq:hh-step1-1}).

\textit{Step} 2. Proof of
%
%
\begin{equation}\label{eq:hh-step2-1}
\int_{|x|\leq4 \tilde R} \int_{|z|\leq|z-x|}h_T(z)h_T(z-x)\,dz\,dx
\lesssim1.
\end{equation}
This time, we decompose the sum over ${|z|\leq|z-x|}$ in two
contributions only: $|z|\leq\tilde R$ and $\tilde R < |z|$.
We then obtain
\begin{eqnarray*}
&&\int_{|x|\leq4 \tilde R} \int_{|z|\leq
|z-x|}h_T(z)h_T(z-x)\,dz\,dx\\
&&\qquad= \int_{|x|\leq4 \tilde R}\underbrace{\int_{\tilde R < |z|\leq
|z-x|}h_T(z)h_T(z-x)\,dz}_{=I_1'}\,dx\\
&&\qquad\quad{} +{} \underbrace{\int_{|x|\leq4\tilde R} \int_{\{
|z|\leq\tilde R\}\cap\{|z|\leq|z-x|\}}h_T(z)h_T(z-x)\,dz\,dx}_{=I_2'}.
\end{eqnarray*}
Proceeding as for $I_1$ in Step 1 using (\ref{eq:assump-h}) yields
\[
I_1' \lesssim1.
\]
For $I_2'$, we use Cauchy--Schwarz inequality, (\ref{eq:assump-h-R=1}),
and $\tilde R\sim1$:
\[
I_2'\leq\biggl(\int_{|x|\leq4 \tilde R} \int_{|z|\leq\tilde
R}h_T^2(z)\,dz\,dx\biggr)^{1/2}\biggl(\int_{|x|\leq4 \tilde R} \int
_{|z'|\leq5\tilde R}h_T^2(z')\,dz'\,dx\biggr)^{1/2} \lesssim1.
\]
This proves (\ref{eq:hh-step2-1}).

\textit{Step} 3. Proof of (\ref{eq:res-h}).

It only remains to use a dyadic decomposition of the ball of
radius $R$ into the ball of radius $\tilde R$
and annuli of the form $2^{-k}R< |z|\leq2^{-k+1}R$, as follows. Taking
$M$ such that $2\tilde R \leq2^{-M}R \leq4 \tilde R$, it holds that
\begin{eqnarray*}
&&\int_{|x|\leq R} \int_{|z|\leq|z-x|}h_T(z)h_T(z-x)\,dz\,dx\\
&&\qquad= \underbrace{\int_{|x|\leq2^{-M}R}\int_{|z|\leq
|z-x|}h_T(z)h_T(z-x)\,dz\,dx}_{
\stackrel{\mbox{\fontsize{8.36pt}{10.36pt}\selectfont{(\ref
{eq:hh-step2-1})}}}{\lesssim}1}\\
&&\qquad\quad{} +\sum_{n=1}^M\underbrace{\int_{2^{n-M-1}R< |x|\leq2^{n-M}R}\int
_{|z|\leq|z-x|}h_T(z)h_T(z-x)\,dz\,dx}_{
\stackrel{\mbox{\fontsize{8.36pt}{10.36pt}\selectfont{(\ref
{eq:hh-step1-1})}}}{\lesssim}(2^{n-M}R)^2} \\
&&\qquad\lesssim 1+R^2\sum_{n=1}^M4^{-n} \sim R^2,
\end{eqnarray*}
which proves (\ref{eq:res-h}).

\textit{Step} 4. Proof of (\ref{eq:res-h-d=2}).

For the case $d=2$, we use the same strategy as for $d>2$.
The bounds on $I_2(n)$ and $I_3(N)$ are the same as for $d>2$.
However the estimate for $I_1$ is slightly worse.
Indeed, we split the dyadic sums $2^{k}R< |z|\leq2^{k+1}R$ into two
categories in order to take advantage of the
fast decay in (\ref{eq:assump-h-d=2}): the first class is for $k$ such that
$2^kR\leq\sqrt{T}$ and the other class for $k$ such that $2^kR> \sqrt{T}$.
More precisely, setting $\mathcal{I}(R,T):=\{k\in\N\dvtx 2^{k-1}R\leq
\sqrt
{T}\}$, we have
\begin{eqnarray*}
I_1&=&\int_{R/2<|z|\leq|z-x|}h_T(z)h_T(z-x)\,dz \\
&\stackrel{\mathrm{Young}}{\leq}& \int_{R/2<|z|}h_T(z)^2\,dz \\
&=&\sum_{k=-1}^\infty\int_{2^kR< |z|\leq2^{k+1}R}h_T(z)^2\,dz
\\
&= & \underbrace{\sum_{k\in\mathcal{I}(R,T)}\int_{2^{k-1}R<
|z|\leq
2^{k}R}h_T^2(z)\,dz}_{\stackrel{\mbox{\fontsize
{8.36pt}{10.36pt}\selectfont{(\ref{eq:assump-h-d=2})}}}{\lesssim}
\max\{
0,\ln(\sqrt{T}R^{-1})\}} + \underbrace{\sum_{k\in\N\setminus
\mathcal
{I}(R,T)}\int_{2^{k-1}R< |z|\leq2^{k}R}h_T^2(z)\,dz}_{
\stackrel{\mbox{\fontsize{8.36pt}{10.36pt}\selectfont{(\ref
{eq:assump-h-d=2})}}}{\lesssim} \sum\limits_{k\in\N}2^{-2k}\lesssim
1}\\
&\lesssim& \max\{1,\ln(\sqrt{T}R^{-1}) \},
\end{eqnarray*}
which gives the extra factor in (\ref{eq:res-h-d=2}).

\begin{appendix}\label{app}
\section*{Appendix: Heuristics for (\protect\lowercase{\ref{pv10bis}}) and
(\protect
\ref{PV11})}

Let $\bar\phi_i$ and $\bar\phi_{T,i}$ denote for $i\in\{1,\ldots
,d\}$
the solutions of
(\ref{PV9}) and (\ref{PV10}), respectively, with $\xi$ replaced by
the ith unit vector $\ee_i$ of $\mathbb{R}^d$.
We claim that
%
%
\setcounter{equation}{0}
\begin{eqnarray}\qquad
\label{eq:append-1}
\sum_{i=1}^d\sum_{j=1}^d\operatorname{var}\Bigl[\sum\bigl(\ee_j\cdot
(A-\langle A\rangle)\ee_i+2\ee_j\cdot\nabla\bar\phi_{i}\bigr)\eta_L\Bigr]
&=&
d \operatorname{var}[a]\sum\eta_L^2 ,\\
\label{eq:append-2}
\sum_{i=1}^d\langle|\nabla\bar\phi_{T,i}-\nabla\bar\phi
_i|^2\rangle
&=&
\operatorname{var}[a] T^{-2} \sum{\bar G_T}^2,
\end{eqnarray}
where $\bar G_T$ denotes the fundamental solution of
the constant coefficient operator $T^{-1}-\triangle$.
We also denote by $\bar G$ the fundamental solution of the Laplacian.
Since
\[
\sum{\bar G_T}^2 \sim
\cases{
T^{2-d/2}, &\quad for $d<4$,\cr
\ln T, &\quad for $d=4$,\cr
1, &\quad for $d>4$,}
\]
and
\[
\sum\eta_L^2 \sim L^{-d},
\]
(\ref{pv10bis}) and (\ref{PV11}) follow from (\ref{eq:append-1}) and
(\ref
{eq:append-2}), that we prove now.

\textit{Step} 1. Argument for (\ref{eq:append-2}).

Since
%
%
\begin{equation}\label{eq:append-3}
-\triangle(\bar\phi_T-\bar\phi)=-T^{-1}\bar\phi_T,
\end{equation}
one has
%
%
\begin{equation}\label{eq:append-4}
\langle|\nabla(\bar\phi_T-\bar\phi)|^2 \rangle=-T^{-1}\langle
\bar\phi_T(\bar\phi_T-\bar\phi) \rangle.
\end{equation}
Rewriting (\ref{eq:append-3}) in the form
\[
T^{-1}(\bar\phi_T-\bar\phi)-\triangle(\bar\phi_T-\bar\phi
)=-T^{-1}\bar\phi
\]
yields the formula
%
%
\begin{equation}\label{eq:append-5}
(\bar\phi_T-\bar\phi)(0)=-T^{-1}\sum_x \bar G_T(x)\bar\phi(x).
\end{equation}
Using (\ref{eq:append-5}), (\ref{eq:append-4}) turns into
%
%
\begin{equation}\label{eq:append-6}
\langle|\nabla(\bar\phi_T-\bar\phi)|^2 \rangle=-T^{-2}\sum
_x\bar
G_T(x)\langle\bar\phi_T(0)\bar\phi(x) \rangle.
\end{equation}
Expressing now $\bar\phi_{T,i}(0)$ and $\bar\phi_i(x)$ in terms the
Green's functions\footnote{Attention should be paid here to turn this
into a rigorous argument since $\bar G$ is not in $L^1(\Z^d)$.} $\bar
G_T$ and $\bar G$,
\begin{eqnarray*}
\bar\phi_{T,i}(0)&=&\sum_{x'}\bar G_T(x')\nabla^* \cdot(A(x')\ee
_i) \\
&=& -\sum_{x'} \nabla_i \bar G_T(x') \bigl(a_i(x')-\langle a \rangle\bigr), \\
\bar\phi_i(x)&=&\sum_{x'}\bar G(x-x')\nabla^* \cdot(A(x')\ee_i) \\
&=& -\sum_{x''} \nabla_i \bar G(x-x'')\bigl(a_i(x'')-\langle a \rangle\bigr),
\end{eqnarray*}
and using the independence of $a_i(x')$ and $a_i(x'')$ for $x'\neq
x''$, we get
\[
\langle\bar\phi_{T,i}(0)\bar\phi_i(x) \rangle = \sum_{x'}
\nabla_i
\bar
G_T(x') \nabla_i \bar G(x-x') \bigl\langle\bigl(a_i(x')-\langle a \rangle\bigr)^2
\bigr\rangle.
\]
Hence,
\begin{eqnarray*}
\sum_{i=1}^d\langle\bar\phi_{T,i}(0)\bar\phi_i(x) \rangle& =
&\operatorname{var}[a]
\sum
_{x'} \nabla\bar G_T(x')\cdot\nabla\bar G(x-x') \\
&=& \operatorname{var}[a] \bar G_T(x) ,
\end{eqnarray*}
since $-\triangle\bar G(x)=\delta(x)$.
Combined with (\ref{eq:append-6}), this proves (\ref{eq:append-2}).

\textit{Step} 2. Argument for (\ref{eq:append-1}).

Using the Green's function, one has
\begin{eqnarray*}
\bar\phi_i(x)&=&\sum_{x'}\bar G(x-x')\nabla^*\cdot\bigl((A-\langle A
\rangle)\ee
_i\bigr)(x') \\
&=&-\sum_{x'}\nabla_i \bar G(x-x')\bigl(a_i(x')-\langle a \rangle\bigr),
\end{eqnarray*}
and therefore
\[
\nabla\bar\phi_i(x) = -\sum_{x'}\nabla\nabla_i \bar
G(x-x')\bigl(a_i(x')-\langle a \rangle\bigr).
\]
Hence, denoting by $\mathcal{A}_{ij}$ the argument of the variance in
(\ref{eq:append-1}), one has
\begin{eqnarray*}
\mathcal{A}_{ij}  :\!&=& \sum_x \bigl(\ee_j\cdot\ee_i\bigl(a_i(x)-\langle a
\rangle\bigr)+2\ee_j \cdot\nabla\bar\phi_i(x)\bigr)\eta_L(x)\\
&=&\sum_x\sum_{x'}\bigl(a_i(x')-\langle a \rangle\bigr)\ee_j \cdot\bigl(\delta
(x-x')\ee
_i-2\nabla\nabla_i \bar G(x-x') \bigr)\eta_L(x).
\end{eqnarray*}
Using the independence of the $a_i$, one obtains for the variance
\begin{eqnarray*}
\operatorname{var}[\mathcal{A}_{ij}]&=&\operatorname{var}[a]\sum_x
\sum_{x'} \sum_{x''}\ee
_j\cdot
\bigl(\delta(x-x')\ee_i-2\nabla\nabla_i\bar G(x-x')\bigr) \\
&&\hspace*{73.4pt}{}\times
\ee_j\cdot\bigl(\delta(x''-x')\ee_i-2\nabla\nabla_i\bar
G(x''-x')
\bigr)\\
&&\hspace*{73.4pt}{}\times\eta_L(x)\eta_L(x'').
\end{eqnarray*}
Rearranging the terms yields
\begin{eqnarray*}
\operatorname{var}[\mathcal{A}_{ij}]&=&\operatorname{var}[a]\sum_x
\sum_{x'} \delta(j-i)
\bigl(\delta
(x-x')-4\nabla_i \nabla_i\bar G(x-x')\bigr)\eta_L(x)\eta_L(x') \\
&&{} + \operatorname{var}[a]\sum_x \sum_{x''} 4 \eta_L(x)\eta_L(x'')
\underbrace
{\sum_{x'}
\nabla_j \nabla_i\bar G(x-x')\nabla_j \nabla_i\bar G(x''-x')}_{=
-\nabla_i\nabla_i \sum_{x'}
\bar G(x-x')\nabla_j\nabla_j \bar G(x''-x') }.
\end{eqnarray*}
Summing in $j$ and using that $-\triangle G(x)=\delta(x)$, this turns into
\begin{eqnarray*}
\sum_{j=1}^d\operatorname{var}[\mathcal{A}_{ij}]&=&\operatorname
{var}[a]\sum_x \sum_{x'}
\bigl(\delta
(x-x')-4\nabla_i \nabla_i\bar G(x-x')\bigr)\eta_L(x)\eta_L(x') \\
&&{} + \operatorname{var}[a]\sum_x \sum_{x''} 4 \eta_L(x)\eta_L(x'')
\nabla
_i\nabla
_i\bar G(x-x'') \\
&=&\operatorname{var}[a]\sum_x \sum_{x'} \delta(x-x')\eta_L(x)\eta
_L(x') \\
&=&\operatorname{var}[a]\sum_x\eta_L(x)^2,
\end{eqnarray*}
from which we deduce (\ref{eq:append-1}).
\end{appendix}

\section*{Acknowledgments}

A. Gloria acknowledges full support and F. Otto acknowledges partial
support of the Hausdorff Center for Mathematics, Bonn, Germany.

%

%
\printaddresses

\end{document}